\documentclass[11pt]{amsart}
\usepackage{amsfonts}
\usepackage{amsmath,amscd}
\usepackage{amsthm}
\usepackage{amssymb}
\usepackage{latexsym}
\usepackage{dsfont}
\usepackage{soul}
\usepackage{cancel}

\usepackage{calligra}
\usepackage[T1]{fontenc}
\usepackage{aurical}
\usepackage[T1]{fontenc}
\usepackage{carolmin}
\usepackage[T1]{fontenc}

\usepackage[normalem]{ulem}

\usepackage{braket}
\usepackage{tikz}
\usepackage{enumitem}   
\usepackage{float}
\usepackage{caption}
\usetikzlibrary{matrix}
\usepackage[linktocpage=true]{hyperref}

\usepackage{frcursive}
\newenvironment{frcseries}{\fontfamily{frc}\selectfont}{}
\newcommand{\textfrc}[1]{{\frcseries#1}}
\newcommand{\mathfrc}[1]{\text{\textfrc{#1}}}

\numberwithin{equation}{section}

\usepackage{soul}
\usepackage{color}

\newtheorem{thm}{Theorem}[section]
\newtheorem{cor}[thm]{Corollary}
\newtheorem{lem}[thm]{Lemma}

\newtheorem{prop}[thm]{Proposition}

\newtheorem{example}[thm]{Example}
\newtheorem{defn}[thm]{Definition}
\newtheorem{conj}[thm]{Conjecture}
\newtheorem{rem}[thm]{Remark}

\definecolor{darkspringgreen}{rgb}{0.09, 0.45, 0.27}

\newcommand\numberthis{\addtocounter{equation}{1}\tag{\theequation}}
\numberwithin{equation}{section}
\begin{document}
	
	\newcommand{\Uq}{U_q sl_2}
	\newcommand{\Uqhat}{U_q(\widehat{sl}_2)}
	\newcommand{\Loop}{\mathcal{L} U_q sl_2}
	
	\newcommand{\beqa}{\begin{eqnarray}}
		\newcommand{\eeqa}{\end{eqnarray}}
	\newcommand{\thmref}[1]{Theorem~\ref{#1}}
	\newcommand{\secref}[1]{Sect.~\ref{#1}}
	\newcommand{\lemref}[1]{Lemma~\ref{#1}}
	\newcommand{\propref}[1]{Proposition~\ref{#1}}
	\newcommand{\corref}[1]{Corollary~\ref{#1}}
	\newcommand{\remref}[1]{Remark~\ref{#1}}
	\newcommand{\er}[1]{(\ref{#1})}
	\newcommand{\nc}{\newcommand}
	\newcommand{\rnc}{\renewcommand}
	
	\newcommand{\CWM}{\mathcal{W}_-^{(N)}}
	\newcommand{\CWP}{\mathcal{W}_+^{(N)}}
	\newcommand{\CG}{\mathcal{G}^{(N)}}
	\newcommand{\CtG}{\widetilde{\mathcal{G}}^{(N)}}
	
	\newcommand{\WMu}{\displaystyle{\sum_{k=0}^{N-1} P_{-k}^{(N)}(u) \cal W_{-k}^{(N)}}}
	\newcommand{\WPu}{\displaystyle{\sum_{k=0}^{N-1} P_{-k}^{(N)}(u) \cal W_{k+1}^{(N)}}}
	\newcommand{\Gu}{\displaystyle{\sum_{k=0}^{N-1} P_{-k}^{(N)}(u) \cal G_{k+1}^{(N)}}}
	\newcommand{\tGu}{\displaystyle{\sum_{k=0}^{N-1} P_{-k}^{(N)}(u) \tilde{\cal G}_{k+1}^{(N)}}}
	
	\nc{\cal}{\mathcal}
	\nc{\diag}{\mathrm{diag}}
	\nc{\goth}{\mathfrak}
	\rnc{\bold}{\mathbf}
	\renewcommand{\frak}{\mathfrak}
	\renewcommand{\Bbb}{\mathbb}

	\newcommand{\epsp}{\varepsilon_+}
	\newcommand{\epsm}{\varepsilon_-}
	\newcommand{\bepsp}{\overline{\varepsilon}_+}
	\newcommand{\bepsm}{\overline{\varepsilon}_-}
	
	\newcommand{\sfr}{\mbox{\mathfrc{s}}}
	
	\newcommand{\bt}{{\bf {\mathsf{T}}}}
	\newcommand{\by}{{\bf {\mathsf{Y}}}}
	\newcommand{\id}{\mathrm{id}}
	\nc{\Cal}{\mathcal}
	\nc{\Xp}[1]{X^+(#1)}
	\nc{\Xm}[1]{X^-(#1)}
	\nc{\on}{\operatorname}
	\nc{\ch}{\mbox{ch}}
	\nc{\Z}{{\bold Z}}
	\nc{\J}{{\mathcal J}}
	\nc{\C}{{\bold C}}
	\nc{\Q}{{\bold Q}}
	\nc{\oC}{{\widetilde{C}}}
	\nc{\oc}{{\tilde{c}}}
	\nc{\ocI}{ \overline{\cal I}}
	\nc{\og}{{\tilde{\gamma}}}
	\nc{\lC}{{\overline{C}}}
	\nc{\lc}{{\overline{c}}}
	\nc{\Rt}{{\tilde{R}}}
	
	\nc{\tW}{\normalfont{{\textsf{W}}}}
	\nc{\tG}{\normalfont{{\textsf{G}}}}
	\nc{\tI}{{\textsf{I}}}

	\nc{\Psig}{\pmb{\psi}}
	\nc{\Phig}{\pmb{\phi}}
	
	\nc{\tE}{{\textsf{E}}}
	\nc{\tF}{{\textsf{F}}}
	\nc{\tK}{{\textsf{K}}}

	\nc{\tx}{{\textsf{x}}}
	\nc{\tho}{{\textsf{h}}}
	\nc{\tk}{{\textsf{k}}}
	\nc{\tep}{{\bf{\cal E}}}

	\nc{\te}{{\textsf{e}}}
	\nc{\tf}{{\textsf{f}}}
	\nc{\tw}{{\textsf{w}}}
	\nc{\tgt}{{\textsf{g}}}

	\nc{\odel}{{\overline{\delta}}}
	
	\def\pr#1{\left(#1\right)_\infty}  
	
	\renewcommand{\P}{{\mathcal P}}
	\nc{\N}{{\Bbb N}}
	\nc\beq{\begin{equation}}
		\nc\enq{\end{equation}}
	\nc\lan{\langle}
	\nc\ran{\rangle}
	\nc\bsl{\backslash}
	\nc\mto{\mapsto}
	\nc\lra{\leftrightarrow}
	\nc\hra{\hookrightarrow}
	\nc\sm{\smallmatrix}
	\nc\esm{\endsmallmatrix}
	\nc\sub{\subset}
	\nc\ti{\tilde}
	\nc\nl{\newline}
	\nc\fra{\frac}
	\nc\und{\underline}
	\nc\ov{\overline}
	\nc\ot{\otimes}
	
	\nc\ochi{\overline{\chi}}
	\nc\bbq{\bar{\bq}_l}
	\nc\bcc{\thickfracwithdelims[]\thickness0}
	\nc\ad{\text{\rm ad}}
	\nc\Ad{\text{\rm Ad}}
	\nc\Hom{\text{\rm Hom}}
	\nc\End{\text{\rm End}}
	\nc\Ind{\text{\rm Ind}}
	\nc\Res{\text{\rm Res}}
	\nc\Ker{\text{\rm Ker}}
	\rnc\Im{\text{Im}}
	\nc\sgn{\text{\rm sgn}}
	\nc\tr{\text{\rm tr}}
	\nc\Tr{\text{\rm Tr}}
	\nc\supp{\text{\rm supp}}
	\nc\card{\text{\rm card}}
	\nc\bst{{}^\bigstar\!}
	\nc\he{\heartsuit}
	\nc\clu{\clubsuit}
	\nc\spa{\spadesuit}
	\nc\di{\diamond}
	\nc\cW{\cal W}
	\nc\cG{\cal G}
	\nc\cZ{\cal Z}
	\nc\ocW{\overline{\cal W}}
	\nc\ocZ{\overline{\cal Z}}
	\nc\al{\alpha}
	\nc\bet{\beta}
	\nc\ga{\gamma}
	\nc\de{\delta}
	\nc\ep{\epsilon}
	\nc\io{ç}
	\nc\om{\omega}
	\nc\si{\sigma}
	\rnc\th{\theta}
	\nc\ka{\kappa}
	\nc\la{\lambda}
	\nc\ze{\zeta}
	\nc{\Thgsr}{\pmb{\grave{\Theta}}}
	\nc{\Thacr}{\pmb{\acute{\Theta}}}
	\nc{\chicr}{\pmb{\acute{\chi}}}
	
	\nc\vp{\varpi}
	\nc\vt{\vartheta}
	\nc\vr{\varrho}
	
	\nc\odelta{\overline{\delta}}
	\nc\Ga{\Gamma}
	\nc\De{\Delta}
	\nc\Om{\Omega}
	\nc\Si{\Sigma}
	\nc\Th{\Theta}
	\nc\La{\Lambda}
	
	\nc\boa{\bold a}
	\nc\bob{\bold b}
	\nc\boc{\bold c}
	\nc\bod{\bold d}
	\nc\boe{\bold e}
	\nc\bof{\bold f}
	\nc\bog{\bold g}
	\nc\boh{\bold h}
	\nc\boi{\bold i}
	\nc\boj{\bold j}
	\nc\bok{\bold k}
	\nc\bol{\bold l}
	\nc\bom{\bold m}
	\nc\bon{\bold n}
	\nc\boo{\bold o}
	\nc\bop{\bold p}
	\nc\boq{\bold q}
	\nc\bor{\bold r}
	\nc\bos{\bold s}
	\nc\bou{\bold u}
	\nc\bov{\bold v}
	\nc\bow{\bold w}
	\nc\boz{\bold z}
	
	\nc\ba{\bold A}
	\nc\bb{\bold B}
	\nc\bc{\bold C}
	\nc\bd{\bold D}
	\nc\be{\bold E}
	\nc\bg{\bold G}
	\nc\bh{\bold H}
	\nc\bi{\bold I}
	\nc\bj{\bold J}
	\nc\bk{\bold K}
	\nc\bl{\bold L}
	\nc\bm{\bold M}
	\nc\bn{\bold N}
	\nc\bo{\bold O}
	\nc\bp{\bold P}
	\nc\bq{\bold Q}
	\nc\br{\bold R}
	\nc\bs{\bold S}
	\nc\bu{\bold U}
	\nc\bv{\bold V}
	\nc\bw{\bold W}
	\nc\bz{\bold Z}
	\nc\bx{\bold X}

	\nc\ca{\mathcal A}
	\nc\cb{\mathcal B}
	\nc\cc{\mathcal C}
	\nc\cd{\mathcal D}
	\nc\ce{\mathcal E}
	\nc\cf{\mathcal F}
	\nc\cg{\mathcal G}
	\rnc\ch{\mathcal H}
	\nc\ci{\mathcal I}
	\nc\cj{\mathcal J}
	\nc\ck{\mathcal K}
	\nc\cl{\mathcal L}
	\nc\cm{\mathcal M}
	\nc\cn{\mathcal N}
	\nc\co{\mathcal O}
	\nc\cp{\mathcal P}
	\nc\cq{\mathcal Q}
	\nc\car{\mathcal R}
	\nc\cs{\mathcal S}
	\nc\ct{\mathcal T}
	\nc\cu{\mathcal U}
	\nc\cv{\mathcal V}
	\nc\cz{\mathcal Z}
	\nc\cx{\mathcal X}
	\nc\cy{\mathcal Y}

	\nc\e[1]{E_{#1}}
	\nc\ei[1]{E_{\delta - \alpha_{#1}}}
	\nc\esi[1]{E_{s \delta - \alpha_{#1}}}
	\nc\eri[1]{E_{r \delta - \alpha_{#1}}}
	\nc\ed[2][]{E_{#1 \delta,#2}}
	\nc\ekd[1]{E_{k \delta,#1}}
	\nc\emd[1]{E_{m \delta,#1}}
	\nc\erd[1]{E_{r \delta,#1}}
	
	\nc\ef[1]{F_{#1}}
	\nc\efi[1]{F_{\delta - \alpha_{#1}}}
	\nc\efsi[1]{F_{s \delta - \alpha_{#1}}}
	\nc\efri[1]{F_{r \delta - \alpha_{#1}}}
	\nc\efd[2][]{F_{#1 \delta,#2}}
	\nc\efkd[1]{F_{k \delta,#1}}
	\nc\efmd[1]{F_{m \delta,#1}}
	\nc\efrd[1]{F_{r \delta,#1}}

	\nc\fa{\frak a}
	\nc\fb{\frak b}
	\nc\fc{\frak c}
	\nc\fd{\frak d}
	\nc\fe{\frak e}
	\nc\ff{\frak f}
	\nc\fg{\frak g}
	\nc\fh{\frak h}
	\nc\fj{\frak j}
	\nc\fk{\frak k}
	\nc\fl{\frak l}
	\nc\fm{\frak m}
	\nc\fn{\frak n}
	\nc\fo{\frak o}
	\nc\fp{\frak p}
	\nc\fq{\frak q}
	\nc\fr{\frak r}
	\nc\fs{\frak s}
	\nc\ft{\frak t}
	\nc\fv{\frak v}
	\nc\fz{\frak z}
	\nc\fx{\frak x}
	\nc\fy{\frak y}
	
	\nc\fA{\frak A}
	\nc\fB{\frak B}
	\nc\fC{\frak C}
	\nc\fD{\frak D}
	\nc\fE{\frak E}
	\nc\fF{\frak F}
	\nc\fG{\frak G}
	\nc\fH{\frak H}
	\nc\fJ{\frak J}
	\nc\fK{\frak K}
	\nc\fL{\frak L}
	\nc\fM{\frak M}
	\nc\fN{\frak N}
	\nc\fO{\frak O}
	\nc\fP{\frak P}
	\nc\fQ{\frak Q}
	\nc\fR{\frak R}
	\nc\fS{\frak S}
	\nc\fT{\frak T}
	\nc\fU{\frak U}
	\nc\fV{\frak V}
	\nc\fZ{\frak Z}
	\nc\fX{\frak X}
	\nc\fY{\frak Y}
	\nc\tfi{\ti{\Phi}}
	\nc\bF{\bold F}
	\rnc\bol{\bold 1}
	
	\nc\ua{\bold U_\A}
	
	\nc\qinti[1]{[#1]_i}
	\nc\q[1]{[#1]_q}
	\nc\xpm[2]{E_{#2 \delta \pm \alpha_#1}}  
	\nc\xmp[2]{E_{#2 \delta \mp \alpha_#1}}
	\nc\xp[2]{E_{#2 \delta + \alpha_{#1}}}
	\nc\xm[2]{E_{#2 \delta - \alpha_{#1}}}
	\nc\hik{\ed{k}{i}}
	\nc\hjl{\ed{l}{j}}
	\nc\qcoeff[3]{\left[ \begin{smallmatrix} {#1}& \\ {#2}& \end{smallmatrix}
		\negthickspace \right]_{#3}}
	\nc\qi{q}
	\nc\qj{q}
	
	\nc\ufdm{{_\ca\bu}_{\rm fd}^{\le 0}}
	
	 \newcommand{\field}{\mathbb{K}} 

	\newcommand{\Rbar}{\overline{\roots}}
	\newcommand{\gfrak}{{\mathfrak{g}}}
	\newcommand{\ogamma}{\overline{\gamma}}
	\newcommand{\ons}{O}
	\newcommand{\ima}{{\mathrm{im}}}
	\newcommand{\real}{{\mathrm{re}}}
	\newcommand{\uqg}{U_q(\gfrak)}
	\newcommand{\roots}{{\mathcal R}}

	
	\nc\isom{\cong} 
	
	\nc{\pone}{{\Bbb C}{\Bbb P}^1}
	\nc{\pa}{\partial}
	\def\H{\mathcal H}
	\def\L{\mathcal L}
	\nc{\F}{{\mathcal F}}
	\nc{\Sym}{{\goth S}}
	\nc{\A}{{\mathcal A}}
	\nc{\arr}{\rightarrow}
	\nc{\larr}{\longrightarrow}
	
	\nc{\ri}{\rangle}
	\nc{\lef}{\langle}
	\nc{\W}{{\mathcal W}}
	\nc{\uqatwoatone}{{U_{q,1}}(\su)}
	\nc{\uqtwo}{U_q(\goth{sl}_2)}
	\nc{\dij}{\delta_{ij}}
	\nc{\divei}{E_{\alpha_i}^{(n)}}
	\nc{\divfi}{F_{\alpha_i}^{(n)}}
	\nc{\Lzero}{\Lambda_0}
	\nc{\Lone}{\Lambda_1}
	\nc{\ve}{\varepsilon}
	\nc{\bepsilon}{\bar{\epsilon}}
	\nc{\bak}{\bar{k}}
	\nc{\phioneminusi}{\Phi^{(1-i,i)}}
	\nc{\phioneminusistar}{\Phi^{* (1-i,i)}}
	\nc{\phii}{\Phi^{(i,1-i)}}
	\nc{\Li}{\Lambda_i}
	\nc{\Loneminusi}{\Lambda_{1-i}}
	\nc{\vtimesz}{v_\ve \otimes z^m}
	
	\nc{\asltwo}{\widehat{\goth{sl}_2}}
	\nc\ag{\widehat{\goth{g}}}  
	\nc\teb{\tilde E_\boc}
	\nc\tebp{\tilde E_{\boc'}}
	
	\newcommand{\LR}{\bar{R}}
		\newcommand{\eeq}{\end{equation}}
	\newcommand{\ben}{\begin{eqnarray}}
		\newcommand{\een}{\end{eqnarray}}

	
	\setcounter{MaxMatrixCols}{30}
	\newcommand{\h}{\frac{1}{2}}
	\newcommand{\tha}{\frac{3}{2}}
	
	\newcommand{\bep}{\overline{\epsilon}_+}
	\newcommand{\bem}{\overline{\epsilon}_-}
	\newcommand{\bkp}{\overline{k}_+}
	\newcommand{\bkm}{\overline{k}_-}
	\newcommand{\kp}{k_+}
	\newcommand{\km}{k_-}
	
	\newcommand{\ds}{\mathds}
	
	\newcommand{\fu}{{\langle 12 \rangle} }
	\newcommand{\futt}{{\langle 23 \rangle} }
	\newcommand{\sfu}{{\langle 34 \rangle} }
	
	\newcommand{\CE}{\cal{E} }
	\newcommand{\CF}{\cal{F} }
	\newcommand{\CH}{\cal{H} }
	\newcommand{\CW}{\cal{W}}
	
	\allowdisplaybreaks

	\makeatletter
	\def\@textbottom{\vskip \z@ \@plus 1pt}
	\let\@texttop\relax
	\makeatother
	

	\title[FRT presentation of the $q$-Onsager algebra]{FRT presentation of the $q$-Onsager comodule algebra\\ and refined $q$-characters} 
		\author{Pascal Baseilhac}
	\author{Azat M. Gainutdinov}
	\address{Universit\'e de Tours, CNRS, Université d'Orl\'eans, IDP UMR 7013, 37200 Tours, FRANCE}
	\email{pascal.baseilhac@idpoisson.fr,azat.gainutdinov@idpoisson.fr}

	\begin{abstract} 
	 We introduce a concept of counital comodule algebra $(B,\delta_B,\epsilon_B)$ over a Hopf ${\mathbb K}$-algebra $H$, where $\delta_B\colon B\rightarrow B \otimes H$ and $\epsilon_B\colon B\rightarrow {\mathbb K}$ are coaction and counit homomorphisms, respectively. This structure gives a canonical (counital) comodule algebra map $\iota_B\colon B \rightarrow H$. In this paper, we study in details the case $B=O_q$, and $H$ is the quantum loop algebra $\Loop$ for ${\mathbb K}={\mathbb C}$. Here, $O_q$ is the $q$-Onsager algebra which admits four presentations, respectively referred to as (i) the (original) {\it $q$-Dolan-Grady} presentation with two generators $W_0,W_1$; (ii) the {\it alternating} presentation;	(iii) the {\it root vectors} presentation;	(iv) the {\it Drinfeld} type presentation.	For $(O_q,\delta_{O_q},\epsilon_{O_q})$, in fact a 3-parameter family of counital comodule algebras, we establish a Faddeev-Reshetikhin-Takhtajan (FRT)  type presentation thanks to a generalization of  Terwilliger's correspondence between (ii) and (iii). This FRT type presentation is written in terms of a K-operator satisfying a reflection algebra relation, and the coaction $\delta_{O_q}$ sends the K-operator to a `dressed' K-operator using Ding-Frenkel type L-operators of $\Loop$; the counit $\epsilon_{O_q}$ sends the K-operator to the most general (non-diagonal) solution of the reflection equation; the map $\iota_{O_q}$ sends the K-operator to a `dressed' K-matrix with entries in $\Loop$. The FRT type presentation allows to compute the image via $\delta_{O_q}$ and $\iota_{O_q}$ of the generating functions of $O_q$ for each presentation  (ii)-(iv). Moreover, using Beck's correspondence between root vectors of $\Loop$ and Drinfeld generators of ${}^{Dr}\!\Loop$ -- the so-called Drinfeld second realization --   the images in ${}^{Dr}\!\Loop$  of Lu-Wang's generating functions of $O_q$ are obtained,  as well as their images via the (ordinary) coproduct of ${}^{Dr}\!\Loop$, which confirm and enhance two theorems by Prze\'zdziecki.  Finally, we apply our FRT presentation to the construction of new type $q$-characters for $\Loop$ in terms of imaginary root vectors of $O_q$, refining those of Frenkel-Reshetikhin, and we explain representation-theoretic meaning of such refining.
    \end{abstract}

	\maketitle
	
	\vskip -0.7cm
	
	{\small MSC:\  81R50; 81R10;  81U15.}

	{{\small  {\it \bf Keywords}: Reflection equation; FRT presentation; $q$-Onsager algebra;  Drinfeld second realization; $q$-characters}}

	\vspace{-0.3cm}
    \setcounter{tocdepth}{1}
	\tableofcontents
	
	\section{Introduction}
\subsection{Background} Hopf algebras associated with
quantum affine algebras\footnote{Here, note that the so-called derivation generator is omitted (see \cite[Rem.\,2, p. 393]{CP95}).}  $U_q(\widehat{g})$ are known to admit at least three presentations\footnote{For $U_q(\widehat{sl}_2)$ note that a fourth presentation called `equitable' has been introduced in \cite{IT03}.}: the first two are, respectively, the  {\it Drinfeld-Jimbo} presentation originally introduced in \cite{Jim,Dr0}
and the so-called {\it Drinfeld second realization} \cite{Dr88}. The Hopf algebra structure is ensured by the existence of coproduct $\Delta$ and counit $\epsilon$ homomorphisms, and antipode $S$ anti-automorphism.
Importantly, the root vectors construction due to Damiani and Beck \cite{Dam93,Beck} using Lusztig's braid group action was a turning point in establishing the explicit isomorphisms relating the two presentations \cite{Beck,Ji,Da14}. A third presentation for the Hopf algebra $(U_q(\widehat{g}),\Delta,\epsilon,S)$ is the so-called {\it FRT} presentation \cite{RS90}, built on previous works  by Faddeev-Reshetikhin-Takhtajan
\cite{FRT89}.  It is given in terms of L-operators satisfying relations of a Yang-Baxter algebra (that usually refers to the `R-matrix presentation' of the quantum algebra), with explicit coproduct and counit maps in terms of the L-operators' entries. For the simplest example $U_q(\widehat{sl}_2)$, the explicit FRT presentation is due to Ding-Frenkel \cite{DF93}. For the extension to higher rank case $U_q(\widehat{g})$, see~\cite{HM97,JLM,JLMBD,LP21}.  Importantly, FRT presentations brought a major contribution to the fields of quantum affine algebras, representation theory and mathematical physics. Among the well-known applications, one finds the definition of $q$-characters introduced by Frenkel and Reshetikhin~\cite{FR99} (see also~\cite{FMu99}) in the context of representation theory.  These are ring homomorphisms from the Grothendieck ring $[\mathrm{Rep}\,U_q(\widehat{g})]$ to a commutative subalgebra in $U_q(\widehat{g})$, so they help to understand algebraic relations in $[\mathrm{Rep}\,U_q(\widehat{g})]$.

Comodule algebras over the quantum affine algebras and related coideal subalgebras have attracted growing interest in recent years~\cite{MRS02,BB09,Ko12,LBG,Lu23,AV24,LP26}. A (right) comodule algebra $B$ over a bialgebra $H$ (over a field $\mathbb K$) is equipped with a coaction map $\delta_B: B \rightarrow B \otimes H$ satisfying certain axioms. If in addition it is equipped with an algebra map $\epsilon_B: B \rightarrow {\mathbb K}$, we call it a counital comodule algebra, denoted $(B,\delta_B,\epsilon_B)$; see Definition~\ref{def:right-counit-com}. We explain by Corollary~\ref{cor:im-iota}  in what sense the counital comodule algebras are generalizations of coideal subalgebras in $H$. In particular, by Lemma~\ref{lem:Bcoid}, their structure provides a canonical algebra map $\iota_B\colon B\rightarrow H$ whose image is a coideal subalgebra in $H$.

Given a counital comodule algebra $B$ in terms of generators and defining relations, it is natural to ask for the corresponding FRT presentation. This amounts to provide the so-called K-operator $K\in  B \otimes \mathrm{End}_{\mathbb{K}}(V)$, for an appropriate $H$-module $V$, satisfying relations of a reflection type algebra~\cite{Skly88}, and which are equivalent to the original defining relations.
Moreover, the coaction map $\delta_B$ should be given in terms of entries of the K- and L-operators, while the counits $\epsilon_B$ applied to the K-operator give K-matrix reflection-equation solutions on $V$. And a further question is what such FRT presentation tells us about  the image via $\iota_B$ of $B$ in $H$, and about the representation theory of the algebras $B$ and $H$.

 In this paper, we focus on the  
\textit{$q$-Onsager algebra} $O_q$, which is a paradigmatic example of a comodule algebra over the quantum loop algebra $H=\Loop$.
 The $q$-Onsager algebra has been introduced  by Terwilliger \cite{Ter03}, and  appeared independently later on in the context of mathematical physics \cite{B05}.  There are four presentations of $O_q$ with generators and relations that can be summarized as follows. 
 \begin{itemize}
 	\item[(i)] The (original) {\it $q$-Dolan-Grady} type presentation  by two generators $\{W_0,W_1\}$ and $q$-Dolan-Grady relations, see Definition \ref{def:Oq} and \cite{Ter03,B05} for details.\vspace{1mm}
 	
 	\item[(ii)] The {\it alternating} type presentation $^{\{a\}}O_q$ by generators $\{W_{-k},W_{k+1},G_{k+1},{\tilde G}_{k+1}|k\in {\mathbb N}\}$ and relations, see Definition \ref{def:Oqalt} and  \cite{BasBel,Ter21}  for details.\vspace{1mm}
 	
 	\item[(iii)] The {\it root vectors} type presentation $^{\{r\}}O_q$ by root vectors $\{B_{n\delta+\alpha_1}, B_{m \delta} \, | \, n\in {\mathbb Z}, m\in \N^*\}$ and relations, see Theorem \ref{thm:BK17} and \cite{BK17} for details.\vspace{1mm}
 	
 	\item[(iv)] The Drinfeld type presentation $^{Dr}O_q$ by generators $\{B_{n\delta+\alpha_1}, {\acute H}_{m}\, | \, n\in {\mathbb Z}, m\in \N^*\}$ and relations, where the generators are the coefficients of Lu-Wang's generating functions ${\bf B}(z)$, ${\bf \Theta}(z)$ in Proposition~\ref{prop:DrOq}; see \cite{LW20} for details.\vspace{1mm}
 \end{itemize}

 Each presentation has its own field of interest. For instance, the presentation (i) has been extensively studied in the context of representation theory: irreducible finite-dimensional representations of $O_q$ (for $q$ not a root of unity) are classified using the theory of tridiagonal pairs; see \cite{IT10} and e.g.\ \cite{INT10,NT17} for a review. Also, $O_q$ occurs in the context of coideal subalgebras~\cite{BB09,BB12}: it is  one of the most basic examples of quantum affine symmetric pairs \cite{Ko12},  also referred as the split affine	$\imath$quantum group of rank one \cite{LW20}.
The presentation (ii) is relevant in studying spectrum problem in such quantum integrable models as the open XXZ spin chain with generic boundary conditions \cite{BK05,BB12,BB16,BGL25}. The presentations (iii) and (iv) relate to  the braid group action for the $q$-Onsager algebra \cite{BK17} and to the $\imath$Hall algebra of the projective line \cite{LW20b,LW22}. Importantly, (iv) is used to build a theory of `boundary' $q$-characters for studying finite-dimensional representations of quantum affine symmetric pairs~\cite{P24,LPr25},  independently of Ito-Terwilliger's result in the case of $O_q$~\cite{IT10}.

\subsection{Goal \& Main results}
 The goal of this paper is firstly to provide an FRT presentation of $O_q$ from the perspective of counital comodule algebras, building on previous works~\cite{BSh1,BasBel,Ter21b,L23,LBG,BGL25} on its alternating central extension $\cal A_q$. Secondly, to use this FRT presentation in order to establish the precise relationship between $O_q$ and Drinfeld second realization of $\Loop$, beyond the partial results in~\cite{P24}. Thirdly, to propose a refinement of Frenkel-Reshetikhin's $q$-characters for  $\Loop$.
More precisely, the main results are: 
\begin{enumerate}[label=(\alph*)]
    \item 
    In Theorem~\ref{thm:FRTOqDr}, we formulate FRT type presentations for a 3-parameter family of counital comodule  $q$-Onsager algebras, in terms of K-operator $\bold K^{\{r\}}(u)$ with twisted-LDU decomposition. The corresponding defining relations of $O_q$ then form a single reflection equation in terms of~$\bold K^{\{r\}}(u)$. Importantly, this presentation is given in terms of generating functions $B^\pm(u^{-2})$ and $\bar{B}(u^{-2})$ of real and imaginary root vectors of $O_q$, where $u$ is a formal parameter. Furthermore, the coaction  $\delta_{O_q}$ and the counit $\epsilon_{O_q}$ maps are given explicitly using Ding-Frenkel's L-operators of $\Loop$ and the most general K-matrix solution~\cite{dVG,GZ94}, respectively.
    
    \item For each of the presentations (ii)-(iv), the
     images of {\it all} the generators of $O_q$ via $\delta_{O_q}$ and the injective homomorphism $\iota_{O_q}: O_q \rightarrow \Loop$ are derived in terms of the corresponding generating functions. For the cases (ii) and (iii) we refer to Corollaries~\ref{cor:DelWG},~\ref{cor:DeltBinKij},~\ref{cor:WKdress},~\ref{cor:BinKij}. For the presentation (iv), the images of ${\bf B}(z)$ and ${\bf \Theta}(z)$ in Drinfeld's second realization ${}^{Dr}\!\Loop$ and their ordinary coproducts in ${}^{Dr}\!\Loop \otimes {}^{Dr}\!\Loop$ are given in Section~\ref{sec:DrOq};  see Theorem~\ref{thm:imOqUqnc} and~Proposition~\ref{prop:imDelOqUqnc}. 
          
\item  The concept of `refined' $q$-character map is introduced and studied in Section~\ref{sec:refqchar}. It is a ring homomorphism $\acute{\chi}_q\colon [\mathrm{Rep}\Loop]\rightarrow  O_q\tilde{H}[[z]]$, where $O_q\widetilde{H}\subset O_q$ is the commutative subalgebra
generated by  Lu-Wang's generators $\{\acute{H}_m|m\in{\mathbb N}^*\}$ which are polynomials in the imaginary root vectors $\{B_{m\delta}|m\in{\mathbb N}^*\}$; see Definition~\ref{def:chicute}, Theorem~\ref{thm:homhcute} and~\eqref{def:refchar}.
Examples of the refined $q$-characters, that depend on two parameters, are computed explicitly in terms of the generating functions $\acute{Y}_a^{\pm 1}$, see Proposition~\ref{prop:acutechi-half}, Example~\ref{valchi1} and Appendix~\ref{app:refchar}. 
At the end of Section~\ref{sec:refqchar} we give a representation-theoretic meaning of the refined $q$-characters, via an algorithm of extracting relations in $[\mathrm{Rep}\Loop]$. 
\end{enumerate}

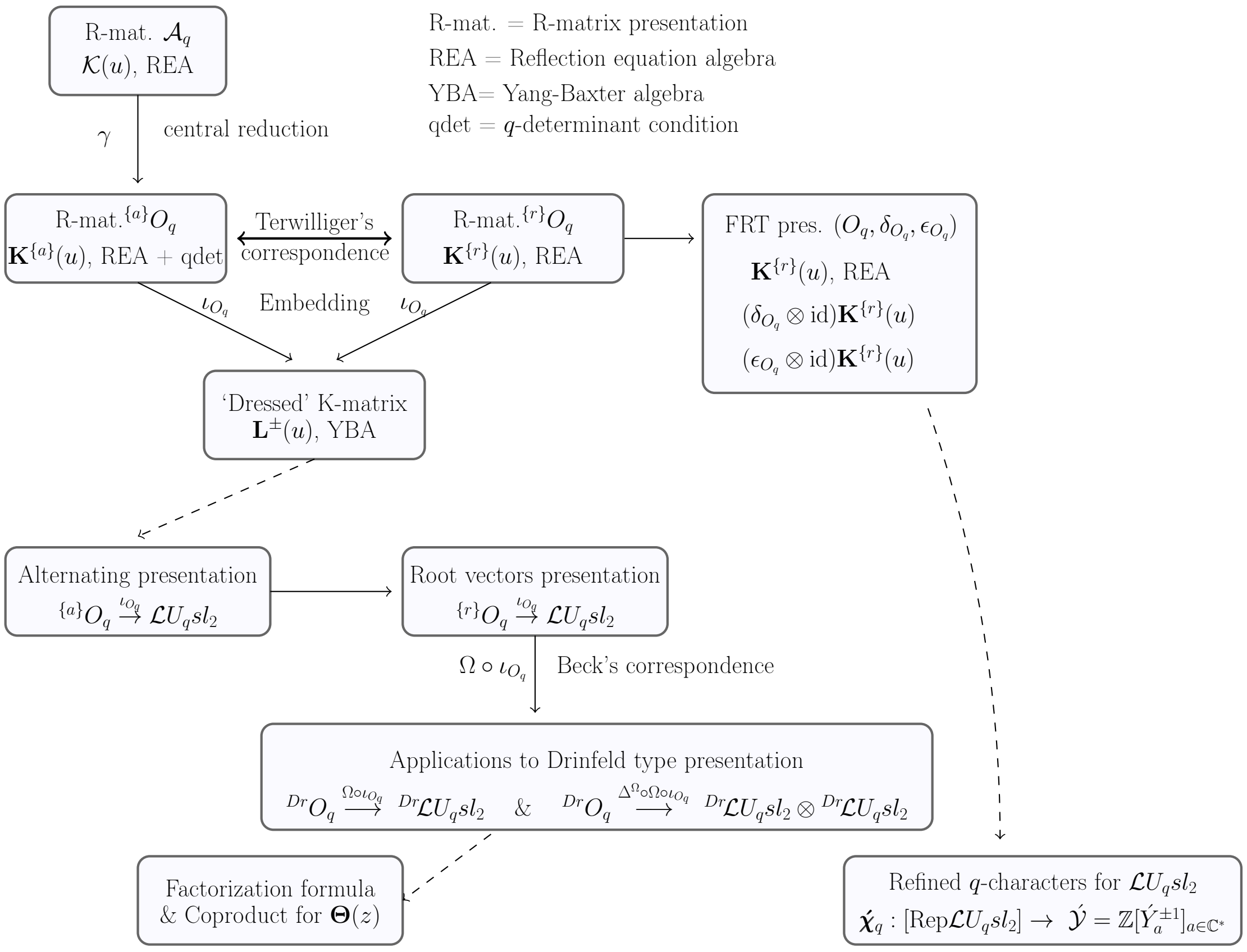
\begin{figure}
		\begin{tikzpicture}[scale=
			0.52, every node/.style={scale=
			0.33}]
			\tikzset{every node/.append style={font=\fontsize{24}{29}\selectfont}}
			
			\def\x{1}
			\def\y{-0.25}

			\draw[color=white!40!black, thick,rounded corners,fill=white!98!blue] (-6, 0) rectangle (-2, -2) node[midway,color=black,above] {R-mat. $\mathcal{A}_q$} node[midway, color=black,below] {$\mathcal{K}(u)$, REA};
			
			\draw[color=white!40!black, thick, rounded corners, fill=white!98!blue] (-7, -4+\y) rectangle (-2, -6+\y) node[midway, color=black,above] {R-mat.$^{\{a\}}O_q$} node[midway, color=black,below] {${\bf K}^{\{a\}}(u)$, REA + qdet};

			\draw[->] (7, -5+\y)  -- (8.6, -5+\y) node[midway,above]{
			};
			
			\draw[color=white!40!black, thick,rounded corners,fill=white!98!blue] (8.8, -4+\y) rectangle (15, -8.5+\y)  	node[black, anchor=west ] (r) at (9.2,-5) {FRT pres. $(O_q,\delta_{O_q},\epsilon_{O_q})$};

		\node[black, anchor=west ] (r) at (9.8,-6) {${\bf K}^{\{r\}}(u)$, REA};
			
	    \node[black, anchor=west ] (r) at (9.6,-7) {$(\delta_{O_q}\otimes \id){\bf K}^{\{r\}}(u)$};
				
		    \node[black, anchor=west ] (r) at (9.6,-8) {$(\epsilon_{O_q}\otimes \id){\bf K}^{\{r\}}(u)$};

			\draw[color=white!40!black, thick,rounded corners,fill=white!98!blue] (-2.5, -8+\y) rectangle (2.5, -10+\y) node[midway, color=black,above] {`Dressed' K-matrix} node[midway,color=black,below]{ ${\bf L}^\pm(u)$, YBA};
			
			\draw[color=white!40!black, thick,rounded corners,fill=white!98!blue] (2, -4+\y) rectangle (7, -6+\y) node[midway, color=black,above] {R-mat.$^{\{r\}}O_q$} node[midway, color=black,below] {${\bf K}^{\{r\}}(u)$, REA};

			\draw[color=white!40!black, thick,rounded corners,fill=white!98!blue] (-7, -12+\y) rectangle (-1, -14+\y) node[midway, color=black,above] {Alternating presentation} node[midway,color=black,below]{$^{\{a\}}O_q \overset{\iota_{O_q}}{\rightarrow} \Loop$};
			
			\draw[color=white!40!black, thick,rounded corners,fill=white!98!blue] (2, -12+\y) rectangle (8, -14+\y) node[midway, color=black,above] {Root vectors presentation} node[midway,color=black,below]{$^{\{r\}}O_q \overset{\iota_{O_q}}{\rightarrow} \Loop$};
			
			\node[black, anchor=west ] (r) at (3.2,-15) {$\Omega \circ \iota_{O_q}$ \quad Beck's correspondence};
			
			\draw[color=white!40!black, thick,rounded corners,fill=white!98!blue] (-1.2, -16+\y) rectangle (14, -18.4+\y) node[midway, color=black,above] {Applications to Drinfeld type presentation} node[midway,color=black,below]{$^{Dr}O_q \overset{\Omega \circ\iota_{O_q}}{\longrightarrow} \ {}^{Dr}\!\Loop$ \quad \& \quad $^{Dr}O_q \overset{\Delta^\Omega \circ \Omega \circ\iota_{O_q}}{\longrightarrow} \ {}^{Dr}\!\Loop\otimes{}^{Dr}\!\Loop$ };
		
		 \draw[color=white!40!black, thick,rounded corners,fill=white!98!blue] (-4, -19+\y) rectangle (2, -21+\y) 
			node[midway, color=black,above] {Factorization formula} 
			node[midway, color=black,below] {\& Coproduct for ${\bf \Theta}(z)$ };

		\draw[color=white!40!black, thick,rounded corners,fill=white!98!blue] (12, -19+\y) rectangle (21, -21+\y) 
			node[midway, color=black,above] {Refined $q$-characters for $\Loop$} 
			node[midway, color=black,below] {$\chicr_q: [\mathrm{Rep}\Loop]  \rightarrow\ \acute{\cal Y} = {\mathbb Z}[\acute{Y}_{a}^{\pm 1}]_{a\in{\mathbb C}^*}$};

			\node[black, anchor=west ] (r) at (2.5,-0.4) {R-mat. =  R-matrix presentation};
			
			\node[black, anchor=west ] (r) at (2.5,-1.2) {REA = Reflection equation algebra};
			
			\node[black, anchor=west ] (r) at (2.5,-2) {YBA= Yang-Baxter algebra};
			
			\node[black, anchor=west ] (r) at (2.5,-2.7) {qdet = $q$-determinant condition};

			\node[black, anchor=west ] (r) at (-1.3,-7.35) {};
			
			\node[black] (r) at (0, -6.5+\y) {Embedding};

			\draw[->] (-4,-2) -- (-4,-4);
			
			\node[black, anchor=west ] (r) at (-5,-3) {$\gamma$};
			\node[black, anchor=west ] (r) at (-3.5,-2.75) {central reduction};
			\node[black, anchor=west ] (r) at (-3.,-3.25) {};
			
			\draw[<->,thick] (-1.75, -5+\y)  -- (1.75, -5+\y) node[midway,above]{Terwilliger's};
			
			\node[black, anchor=west ] (r) at (-1.75,-5.6) {correspondence};
			\node[black, anchor=west ] (r) at (-3.,-3.25) {};
			
			\draw[->] (-4,-6.25) -- (-0.5,-8) node[midway, above] {$\iota_{O_q}$};
			
			\draw[->] (4,-6.25) -- (0.5,-8) node[midway, above] {$\iota_{O_q}$};
			
			
			\draw[->,dashed] (0, -10+\y) -- (-4, -11.75+\y) ;
			

			\draw[->] (-1, -13+\y) --   (1.75, -13+\y) ;
			
			\draw [->] (5, -14+\y) -- (5, -15.75+\y) ;

			\draw[->,dashed] (4, -18.5+\y) -- (2, -20+\y) ;
			


    \draw[dashed 
		, ->] (13.9,-9.1) arc (19:0:30);
            
		\end{tikzpicture}
	\caption{The path followed starting from $\cal A_q$, the alternating central extension of the $q$-Onsager algebra. Here, $\Omega$ denotes an anti-automorphism of $\Loop$.}
		\label{fig:overview}
	\end{figure}

The path leading to the above listed results is the following,  summarized in Figure \ref{fig:overview}. 
We start from the R-matrix presentation of ${\cal A}_q$ \cite{BSh1}. Through a central reduction via a surjective homomorphism 
$\gamma: {\cal A}_q \rightarrow O_q$ we get the R-matrix presentation of $^{(a)}O_q$ from presentation (ii). 
Generalizing Terwilliger's correspondence \cite{Ter21b}  relating $^{(a)}O_q$ and $^{(r)}O_q$ in presentation (iii), we obtain
 the R-matrix presentation of $^{(r)}O_q$, and consequently for $^{Dr}O_q$.
Based on previous results~\cite[Sec.\,6.2]{LBG} and~\cite[Sec.\,2]{BGL25}, we formulate the coaction map $\delta_{O_q}$ and counit $\epsilon_{O_q}$, and as a result we get the `dressed' K-matrix form for the algebra map $\iota_{O_q}$. Altogether, this gives the FRT presentation for $(O_q,\delta_{O_q},\epsilon_{O_q})$ which is displayed in the upper right box on the figure. The lower half of the figure connected by dashed arrows highlights the applications (b) and (c) of our FRT presentation.

Concerning the above listed  main results, some comments are in order. 
First about (a), we introduce in Proposition~\ref{prop:coact} a family of counital $\Loop$-comodule algebra structures on $O_q$ with 4 parameters $(k_\pm,\varepsilon_\pm)$. Then, according to the equivalence relations between counital comodule algebras in Definition~\ref{def:equiv-rel}, we reduce it in Proposition~\ref{prop:isocomod} to a  3-parameter family. 
 We further conjecture in Conjecture~\ref{conj:IsoCom} that $O_q$ admits a 3-parameter family of \textit{inequivalent} counital right comodule $\Loop$-algebra structures.
We would like to point out that in the context of twisted $q$-Yangians~\cite{MRS02}, R-matrix presentation of $O_q$  was studied in~\cite{Lu23} for one particular choice (with $\varepsilon_\pm=0$) among this 3-parameter family, see also Remark~\ref{rem:qtwist}.

Concerning (b), our results give a complete answer  to the question addressed in \cite[Sec.\,1.5]{P24}: ``What is the relation between the Drinfeld presentation of a quantum affine algebra and Lu and Wang’s Drinfeld-type presentation of a QSP coideal subalgebra?''. As corollaries of our results, independent proofs of two theorems by Prze\'zdziecki follow; see Corollaries~\ref{cor:finnc} and~\ref{cor:DelTh}.
 With respect to the images via $\delta_{O_q}$ of the generating functions of root vectors, let us point out that  the $n$-fold coaction $\delta_{O_q}^{(n)} = (\delta_{O_q}^{(n-1)}\otimes \id)\circ \delta_{O_q}$ of any generating function -- starting from Corollary~\ref{cor:DeltBinKij} with Lemma~\ref{copKdressij} -- becomes increasingly complicated with $n$. This situation is similar to the case of quantum affine Hopf algebras with `ordinary' coproduct $\Delta$. To circumvent the increasing complexity, in~\cite{Dr88} Drinfeld introduced the so-called `Drinfeld coproduct' defined in an appropriate
completion $\Loop \widehat{\otimes} \Loop$, see e.g.~\cite{He06,Da24} for more details. Clearly, our coaction map $\delta_{O_q}$ can be viewed as `ordinary' in similar sense. 
A coaction analogue of the Drinfeld coproduct for $O_q$ was obtained in~\cite[Prop.\,4.15]{LP26}.

Concerning the above part (c), key ingredients of the construction are the universal boundary transfer matrices for $O_q$ and their corresponding TT-relations previously obtained in~\cite{L23,AV24b,BGL25}, and the introduction of an homomorphism $\acute{h}_q: O_q \rightarrow O_q\tilde{H}$ generalizing the Harish-Chandra map of Frenkel-Reshitikhin \cite{FR99}. Furthermore, in Section~\ref{sec:defchi} we discuss how  the refined $q$-characters are related to a conjectured universal K-matrix for $O_q$, culminating in Conjecture~\ref{conj:chiunivK}.  For related comments on this latter point, see also~\cite[Sec.\,1.6.1]{AV24b}.

\subsection{Perspectives} The following directions of research may be highlighted.

\begin{itemize}
\item {\bf Extension to higher rank}. Generalized $q$-Onsager algebras  $O_q(\widehat{g})$ introduced in~\cite{BB09} and their classical analogs  in~\cite{UI,S18} are higher rank generalizations of the $q$-Onsager algebra $O_q\equiv O_q(\widehat{sl}_2)$ and its classical analog $O$~\cite{Ons44}, respectively.  From the perspective of quantum affine symmetric pairs, generalized $q$-Onsager algebras are isomorphic to affine $\imath$quantum groups of split type for which Drinfeld type presentations are known~\cite{LW20,Z21}. 
Examples of $K$-operators for $O_q(\widehat{g})$ have already been considered in the literature, with applications in quantum integrable models, see~\cite{BF11,KP,KOY19}. The results presented here suggest to investigate counital comodule algebras associated with $O_q(\widehat{g})$ and their relation with its Drinfeld type presentation~\cite{LW20}. In this direction, the R-matrix presentations of $U_q(\widehat{g})$ completed in \cite{HM97,JLM,JLMBD,LP21} are crucial.

\item {\bf Towards a universal K-matrix.} As pointed out in Section~\ref{sec:defchi}, by analogy with the case of $\cal A_q$~\cite{LBG}, it is natural to ask for a universal K-matrix\footnote{The existence of (tensor) universal K-matrices is proven in~\cite{AV24} for a choice of so-called twist pair that differs from the one considered in~\cite{LBG};}   $\mathfrak{K}\in O_q \widehat{\otimes} \Loop$ (in an appropriate completion of $O_q \otimes \Loop$) that satisfies the axioms~\cite[Def.\,2.7 for $B=O_q$, $H=\Loop$]{LBG}, and such that it specializes to our K-operator $\bold K^{\{r\}}(u)$ given by~\eqref{LDU-Kop}.  An explicit expression for $\mathfrak{K}$ in terms of root vectors of $O_q$ and $\Loop$ is clearly desirable, especially in order to  prove and investigate further our refined $q$-character Conjecture~\ref{conj:chiunivK}.
Results in~\cite{LP26}, namely the derivation of an analog of the Drinfeld coproduct for quantum affine symmetric pairs should help on this problem.
\end{itemize}

\subsection{Text organization} In Section~\ref{sec:councom}, we introduce the concept of counital comodule algebra. Our main example, the $q$-Onsager algebra, is introduced in this context. The corresponding counital comodule algebra $(O_q,\delta_{O_q},\epsilon_{O_q})$ in its original presentation (i), with coaction $\delta_{O_q}$, counit $\epsilon_{O_q}$ and injection $\iota_{O_q}=(\epsilon_{O_q}\otimes \id)\circ \delta_{O_q}$ is given; see Definition~\ref{def:Oq}, Proposition~\ref{prop:coact} and Remark~\ref{rem:iotadef}, respectively.
In Section~\ref{sec:Oqcomod}, the three other known presentations (ii), (iii) and (iv) of $O_q$ that use generating functions are briefly recalled, as well as how they relate to each other using a generalized version of the correspondence between (ii) and (iii) due to Terwilliger. In particular, we introduce an auxiliary generating function $\bar B(z)$, see Definition~\ref{def:Bbar}, that plays a central role in further sections.
Section~\ref{sec:FRTOq} is devoted to the construction of the FRT presentation of $(O_q,\delta_{O_q},\epsilon_{O_q})$, whereas the images via $\delta_{O_q}$ and $\iota_{O_q}$ of the generating functions of $O_q$ are computed in Section~\ref{sec:actgen}. In Section~\ref{sec:DrOq}, after briefly reviewing Drinfeld second realization ${}^{Dr}\!\Loop$,  we derive the images of Lu-Wang's Drinfeld type generating functions ${\bf B}(z)$, ${\bf \Theta}(z)$  in ${}^{Dr}\!\Loop$ and their ordinary coproducts. This is done combining an anti-automorphism of $\Loop$ denoted $\Omega$ and Beck's correspondence (see Table~\ref{tab:Omx}) between root vectors for $\Loop$ and Drinfeld generators of ${}^{Dr}\!\Loop$.  Finally, in Section~\ref{sec:refqchar} we revisit the construction of standard $q$-characters using the FRT presentation of $\Loop$~\cite{FR99,FMu99}, emphasizing an approach based on the standard TT-relations. Then, we apply our FRT presentation to the construction of `refined' $q$-characters for $\Loop$ based on our K-operator~\eqref{LDU-Kop}.
 In Appendix~\ref{ap:Loop}, we review the Drinfeld-Jimbo presentation of $\Loop$, its root vectors and evaluation representations.
In Appendix~\ref{ap:AqRmat}, the alternating central extension of the $q$-Onsager algebra $\cal A_q$ is recalled.  In Appendix~\ref{ap:Bbarsol}, the generating function $\bar B(z)$ is studied in details for the standard and non-standard cases. In Appendix~\ref{app:refchar}, explicit examples of refined $q$-characters are given.

\smallskip

{\bf Notations}
	We denote  the set of natural numbers ${\mathbb N} = \{0$, $1$, $2$, $\ldots\}$,  integers ${\mathbb Z} = \{0, \pm 1, \pm 2, \ldots\}$ and complex numbers ${\mathbb C}$. We denote ${\mathbb N}^* = {\mathbb N}\setminus\{0\}$, ${\mathbb Z}^* = {\mathbb Z}\setminus\{0\}$, ${\mathbb C}^* = {\mathbb C}\setminus\{0\}$.
	$\field$ is any field and we assume $q\in\field^*$ and  $q$ is not a root of unity. The $q$-commutator  is defined by
	\beqa
	\big[X,Y\big]_q=qXY-q^{-1}YX \label{defcom}
	\eeqa
	and $\big[X,Y\big]=\big[X,Y\big]_1=XY-YX$ and we denote the $q$-number $[n]_q= (q^n-q^{-n})/(q-q^{-1})$. 
	
	In this paper, all algebras are unital, all scalars are in ${\mathbb K}$. We use the following conventions for algebras: $A = \bigl\langle x_i \, , \,  i\in I \big|\, f_k({x_i}_{i\in I})=0\,\bigr\rangle$ means a unital associative algebra over $\mathbb{C}$ generated by $x_i$ with the defining relations $f_k$. Similarly,  $B = \bigl\langle x_i \, \big| \,  i\in J\subset I \,\bigr\rangle$
	means that it is  a unital subalgebra $B\subset A$ generated by $x_i$'s for all $i\in J$.
 
	Let 
    ${\cal F}(u)=\sum_{n\in \mathbb{Z}} f_n u^{-n}$
    be a formal Laurent series with coefficients in an algebra (like $O_q$), for this we use the notation $O_q[[u^{-1}]]$.
	We use the convention that every rational function of the form $1/p(u)$, where $p(u)$ is a Laurent polynomial, is expanded in $u^{-1}$.

In a few equations, e.g.~\eqref{qdetOq}, we set
			\beqa
			{\cal P}^{-}_{12}=(1-P)/2 \label{eq:P12}
			\eeqa
			where $P=R(1)/(q-q^{-1})$,
		with R-matrix (\ref{R}) and $P$ is the flip operator, i.e.\ $P(a\otimes b) = b\otimes a$.
	
	The subscript $[1]$ attached to an algebra-valued matrix $T$ indicates on which factor in the tensor
product of modules the entries of $T$ act nontrivially, and we use the convention 
\beqa
((T)_{[\mathsf 2]}(T')_{[\mathsf 1]}(T'')_{[\mathsf 2]})_{ij}
		=\sum_{k,\ell=1}^2 (T')_{k\ell} \otimes  (T)_{ik}(T'')_{\ell j}\ .\label{eq:TTT}
\eeqa	

{\bf AI disclosure.} The initial draft of this paper, including most of computations, was completed in
summer 2026, without the use of AI. Furthermore, AI was not used in any capacity (generating ideas, proof formalization or verification, proofreading, writing or otherwise) in the preparation of this manuscript.

\section{On counital comodule algebras and $O_q$}\label{sec:councom}
We first recall  the standard definition of a comodule algebra over a bialgebra, and then introduce a new notion of counital comodule algebras together with equivalence relations among them in Section~\ref{sec:comod-alg}. As a first example, we provide an infinite discrete family of inequivalent counital comodule algebras over Sweedler's bialgebra. We also briefly discuss their categorical interpretation in terms of module categories. In Section~\ref{sec:comod-Oq}, we give as our main example a 3-parameter continuous family of counital comodule algebras based on the $q$-Onsager algebra $O_q$.

\subsection{Counital comodule algebras}\label{sec:comod-alg}
	Let $H$ be a bialgebra over a field $\field$. This means that $H$ is an associative $\field$-algebra with the coproduct $\Delta\colon H\rightarrow H \otimes H$ and the counit  $\epsilon\colon H \to \field$, which are both  algebra homomorphisms, and additionally $(\epsilon\otimes\id) \circ \Delta = \id =(\id \otimes \epsilon) \circ \Delta$ and $\Delta$ is co-associative, i.e.\ $(\Delta\otimes\id) \circ \Delta = (\id \otimes\Delta) \circ \Delta$. 
	
 \begin{defn} \label{def:right-counit-com}
 A right comodule algebra $B$ over $H$ is an associative $\field$-algebra with an algebra map 
 $$
 \delta_B \colon B \rightarrow B  \otimes H\ ,
 $$
 which we call the right coaction of $B$, such that the following conditions hold:
\begin{itemize}
	\item Co-associativity:
	\begin{equation} \label{coassoc}
	(\id_B \otimes \Delta) \circ \delta_B = (\delta_B \otimes id_H) \circ \delta_B \ ,
	\end{equation}
	\item Right counital:
	\begin{equation} \label{right-counital}
	(\id \otimes \epsilon) \circ \delta_B = \id\ ,
	\end{equation}
\end{itemize}

A right comodule algebra $B$  is called counital if  we have an algebra map  $\epsilon_B \colon B \to \field$. 

For brevity, we often write  `counital comodule algebra' instead of `counital right comodule algebra'.
 \end{defn}
 
 We note that every counital right $H$-comodule algebra has automatically a left counitality property, generalizing~\eqref{right-counital} in the following sense:

\begin{lem}\label{lem:Bcoid}
    Let $B$ be a counital right $H$-comodule algebra. 
Then the kernel of the induced algebra map 
	\begin{equation} \label{iota-B}
	\iota_B := (\epsilon_B \otimes \id ) \circ \delta_B \colon B \to H \; ,
	\end{equation}
	which is a two-sided ideal in $B$, is  a right comodule over $H$, i.e.\ 
	\begin{equation} \label{left-counital}
	\delta_B\bigl(\ker(\iota_B)\bigr) \subset \ker(\iota_B) \otimes H\ .
	\end{equation}
\end{lem}
\begin{proof}
Applying the counit $\epsilon_B$ on the first tensor factor to both sides of the co-associativity condition~\eqref{coassoc}, we
get that the $\iota_B$  commutes with both coactions in the following sense:
	\begin{equation}\label{eq:Delta-delta}
	 \Delta \circ \iota_B  =
	(\iota_B\otimes \id)\circ\delta_B \ .
	\end{equation}
	From this equation we  conclude $(\iota_B\otimes \id)\circ\delta_B (\ker(\iota_B)) = 0$ or equivalently $\delta_B (\ker(\iota_B)) \subset \ker(\iota_B\otimes \id)$. But because $ \ker(\iota_B\otimes \id) = \ker(\iota_B) \otimes H$, by left exactness of $\otimes$ bi-functor in the first argument, and so it preserves the kernels, we get~\eqref{left-counital}.
\end{proof}

\begin{rem}
In the above Definition~\ref{def:right-counit-com}, presence of the  counit $\epsilon_B$ tells that the underlying algebra $B$ is the so-called algebra with augmentation (it possess a 1-dimensional representation), while the  property~\eqref{left-counital} guarantees that the quotient of $B$ by the two-sided ideal $\ker(\iota_B)$ is again a right comodule algebra over $H$. 
 This quotient is additionally counital if $\epsilon_B(\ker(\iota_B))=0$.
\end{rem}

It is clear that right coideal subalgebras $B$ in $H$ are examples of counital right comodule algebras over~$H$: in this case $\delta_B := \Delta|_{B}$ and $\epsilon_B := \epsilon|_{B}$, and the kernel $\ker(\iota_B)$ is trivial. However, in general  counital comodule algebras might have non-trivial $\ker(\iota_B)$ as in the following example.

\begin{example}\label{ex:Sw}
   Consider Sweedler's 4-dimensional bialgebra
\[ \mathsf{Sw} = \bigl\langle x,g \, \big|\,  gx = -xg\ , \quad g^2 = 1\ , \quad x^2 = 0 \bigr\rangle\ , \]\  
with the coproduct and counit
defined by
\begin{equation*}
\Delta(g) = g \otimes g\ , \quad 
\Delta(x) = x \otimes 1 + g \otimes x\ ,\qquad
\epsilon(g) = 1\ , \quad \epsilon(x) = 0\ .
\end{equation*}
For any pair of positive integers $m,n \geq 1$ consider the $4mn$-dimensional algebra
\[ B_{m,n} = \bigl\langle \tilde{g}\ , \tilde{x} \, \big|\, \tilde{g}\tilde{x} = -\tilde{x}\tilde{g}\ , \: \tilde{g}^{2m} = 1\ , \: \tilde{x}^{2n} = 0 \bigr\rangle \]
with the coaction
$$
\delta_B(\tilde{g}) = \tilde{g} \otimes g\ ,\qquad
\delta_B(\tilde{x}) = \tilde{x}\otimes1 + \tilde{g}\otimes x\ ,
$$
which makes it a right comodule algebra over $H=\mathsf{Sw}$.
Fix such a pair $(m,n)$ and denote $B=B_{m,n}$. For each $2m$-th root of unity $\zeta_k$,  for $1\leq k \leq 2m$, the algebra $B_{m,n}$ admits the counits  $\epsilon_B$ such that
$$
\epsilon_B(\tilde{g}) = \zeta_k\ , \qquad \epsilon_B(\tilde{x})=0\ .
$$
It is readily seen that the induced algebra map $\iota_B$ from~\eqref{iota-B} acts as follows:
$$
\iota_B(\tilde{g}) = \zeta_k g\ , \qquad \iota_B(\tilde{x})= \zeta_k x\ ,
\qquad \iota_B(\tilde{x}^2)=0\ .
$$ 
 In fact,  $\ker(\iota_B)$ is generated, as a two-sided ideal, by $\tilde{x}^2$:
$$
\ker(\iota_B) = \bigl\langle \tilde{g}^i \tilde{x}^j \, \big|\, 1\leq i \leq 2m, 1 < j \leq 2n \bigr\rangle \ .
$$
It can be verified directly that the property~\eqref{left-counital} holds indeed, in particular $\delta_B(\tilde{x}^2) = \tilde{x}^2\otimes 1$. We thus get a discrete infinite family of counital right comodule algebras $B_{m,n}$ over $\mathsf{Sw}$ with $2m$ distinct counital maps $\epsilon_B$ for each $1\leq k\leq 2m$, and the image of $\iota_B$ in all these cases is the whole $\mathsf{Sw}$.  

Consider furthermore the unital $2n$-dimensional subalgebra of $B_{2,n}$ generated by $y:=\tilde{g}\tilde{x}$, it is in fact isomorphic to the algebra of `truncated' polynomials $\mathbb{C}[y]/\langle y^{2n}\rangle$. It has a structure of a counital right comodule algebra over $\mathsf{Sw}$. Indeed, for $1\leq k\leq n$, $\delta_B(y^k) = y^k\otimes g + y^{k-1} \otimes gx$ for odd $k$ and simply $\delta_B(y^k) = y^k\otimes 1$ for even $k$, verifying~\eqref{coassoc} and~\eqref{right-counital}. Its counit $\epsilon_B$ sends $y$ to zero. The kernel $\ker(\iota_B)$ is generated by $y^2$, verifying~\eqref{left-counital}. The corresponding image under $\iota_B$ is the 2-dimensional  coideal subalgebra of $\mathsf{Sw}$ generated by $gx$, for any $n\geq1$.
\end{example}

 The following corollary shows that the counital comodule algebras are generalisations of coideal subalgebras in $H$ in the sense that the images of $\iota_B$ are always coideal subalgebras, as we saw in the previous example.

\begin{cor}\label{cor:im-iota}
A counital comodule algebra $B$ provides via the image of the algebra map $\iota_B$ from~\eqref{iota-B} a right coideal subalgebra in $H$. Moreover, both the counits $\epsilon$ and $\epsilon_B$ are compatible in the following sense:
\begin{equation}\label{eq:iota-B-iota}
    \epsilon\bigl(\iota_B(b)\bigr) = \epsilon_B(b)\ , \qquad \text{for any}\quad b\in B\ . 
\end{equation}
\end{cor}
\begin{proof}
	The equation~\eqref{eq:Delta-delta} shows that $\Delta\bigl(\iota_B(B)\bigr) \subset \iota_B(B) \otimes H$, i.e.\ the image $\iota_B(B)$ is a right coideal. 
	Moreover, using the right counital condition~\eqref{right-counital} we see from~\eqref{iota-B} that~\eqref{eq:iota-B-iota} holds.
\end{proof}

\begin{defn}\label{def:equiv-rel}
Two counital right $H$-comodule algebras $B$ and $\tilde{B}$ are called isomorphic, or equivalent, if there is an algebra isomorphism between them such that it commutes with all the structures: the coaction and the counit. More precisely,
an algebra isomorphism $\psi: B \to \tilde{B}$ should satisfy:
\begin{equation}\label{eq:comod-alg-equiv}
\begin{split}
(\psi\otimes \id)\circ \delta_B(b) &= \delta_{\tilde{B}}\bigl(\psi(b)\bigr),\\
 \epsilon_{B}(b) &= \epsilon_{\tilde{B}}\bigl( \psi(b) \bigr)
\end{split}
\qquad \text{for any}\quad b\in B\ .
\end{equation}
More generally, we say that an algebra homomorphism $\psi: B \to \tilde{B}$ is the one of counital right $H$-comodule algebras if it satisfes~\eqref{eq:comod-alg-equiv}.
\end{defn}

\begin{example}
    Consider $H$ as a counital right comodule algebra over $H$ with the coaction $\delta_H=\Delta$ and the counit $\epsilon_H=\epsilon$. Then,
combining~\eqref{eq:Delta-delta} and~\eqref{eq:iota-B-iota}, we see for any $B$ that $\iota_B\colon B \to H$ is a homomorphism of counital right comodule algebras.
\end{example}

\begin{example}
Continuing Example~\ref{ex:Sw},
for each $B_{m,n}$, we have $2m$ distinct counits $\epsilon_B$ parametrized by $2m$-th roots of unity $\zeta_k$, for $1\leq k \leq 2m$. However, for a fixed $(m,n)$, the corresponding counital right comodule algebras are all pairwise equivalent. Indeed, let $\tilde{B}$ be $B_{m,n}$ as a right comodule algebra and with $\epsilon_{\tilde{B}}(\tilde{g})=1$. The algebra automorphisms $\psi_k$ on $B_{m,n}$ defined by
$$
\psi_k(\tilde{g}) = \zeta_k\tilde{g}\ , \qquad
\psi_k(\tilde{x}) = \tilde{x}\ ,
$$
provide equivalence of $B_{m,n}$ with $\epsilon_{\tilde{B}}(\tilde{g})=\zeta_k$ to $\tilde{B}$, as counital right $\mathsf{Sw}$-comodule algebras, e.g.\ both equations in~\eqref{eq:comod-alg-equiv} are satisfied.

Furthermore, for a fixed $m$, the canonical quotient from $B_{m,n+2l}$ to $B_{m,n}$ for any positive $l$ is a homomorphism of counital right $\mathsf{Sw}$-comodule algebras.
\end{example}

\begin{rem}
In the context of category theory, Definition~\ref{def:right-counit-com} means that the category $\mathcal{M}=B\mathrm{-mod}$ of $B$-modules is a right module category over the tensor category $\mathcal{C}=H\mathrm{-mod}$, and is equipped with a distinguished  object of dimension one. The right action of an $H$-module $V$ on a $B$-module $M$ is $M\otimes V$ with $B$-action given by the coaction $\delta_B$. The pullback along the counital map $\iota_B\colon B \to H$ corresponds to a $\mathcal{C}$-module functor from the regular module category $\mathcal{C}$ to the 
module category $\mathcal{M}$ such that it sends the tensor unit of $\mathcal{C}$ to the distinguished object of $\mathcal{M}$.
In this context, the main idea of  Definition~\ref{def:equiv-rel} is that the pullback along $\psi$ gives rise to a $\mathcal{C}$-module functor preserving the distinguished objects, and inequivalent counital comodule algebras give inequivalent $\mathcal{C}$-module categories with a distinguished object.
\end{rem}
	
\subsection{Counital comodule algebra structures on $O_q$}
\label{sec:comod-Oq}
As our main example, we now introduce a family of counital right comodule algebras over the quantum loop algebra $\Loop$. We  refer to Appendix~\ref{ap:Loop} for the definition of the bialgebra structure on $H=\Loop$.

The $q$-Onsager algebra $O_q$ introduced in~\cite{Ter93,Ter03,BK05} admits a presentation with two generators that are subject to $q$-deformed analogs of the Dolan-Grady relations \cite{DG}.
\begin{defn}\label{def:Oq} The $q$-Onsager algebra ${O}_q$ is the associative algebra with scalar $\rho\in{\field}^*$ and generators $W_0,W_1$ that satisfy:
\begin{align} \label{qDG1}
	\lbrack  {\normalfont W}_0, \lbrack {\normalfont W}_0, \lbrack {\normalfont W}_0, {\normalfont W}_1\rbrack_q \rbrack_{q^{-1}} \rbrack &=\rho \lbrack {\normalfont W}_0, {\normalfont W}_1 \rbrack\ ,\\ \label{qDG2}
	\lbrack {\normalfont W}_1, \lbrack {\normalfont W}_1, \lbrack {\normalfont W}_1, {\normalfont W}_0\rbrack_q \rbrack_{q^{-1}}\rbrack &= \rho  \lbrack {\normalfont W}_1, {\normalfont W}_0 \rbrack\ .
\end{align}   
\end{defn}
\begin{rem} Let $\tilde{O}_q$ be the algebra with generators $\tilde{W}_0,\tilde{W}_1$ and relations~\eqref{qDG1}, \eqref{qDG2} with the substitution $\rho \rightarrow \tilde\rho\in{\field}^*$.
The algebras $O_q$ and $\tilde O_q$  are isomorphic. The isomorphism $\Psi\colon O_q \rightarrow \tilde{O}_q$ is such that 
\beqa\label{eq:Psi}
\Psi(W_0)= \sqrt{\frac{\rho}{\tilde\rho}}\tilde{W}_0, \qquad
\Psi(W_1)= \sqrt{\frac{\rho}{\tilde\rho}}\tilde{W}_1\ .
\eeqa
\end{rem}

We now introduce  
a counital $\Loop$-comodule algebra structure on $O_q$ with  parameters $(k_\pm,\varepsilon_\pm)$.

\begin{prop}\label{prop:coact} Let $k_\pm\in\field^*$ and $\varepsilon_\pm\in\field$. 
The following formulas
\beqa
	\delta_{O_q}(W_0) &=& W_0 \otimes \tk_1 + 1 \otimes (k_+q\te_1 \tk_1 + k_-\tf_1)\ ,\label{eq:del1}\\
	\delta_{O_q}(W_1)&=& W_1 \otimes \tk_0 +  1 \otimes(k_-q\te_0 \tk_0 + k_+\tf_0) \ ,\label{eq:del2}
\eeqa
define coaction $\delta_{O_q}$ if and only if
$\rho=k_+ k_-(q+q^{-1})^2$. 
The counit $\epsilon_{O_q}$ is given by:
\beqa
	\epsilon_{O_q}(W_0) = \varepsilon_+ \ ,\quad \epsilon_{O_q}(W_1) = \varepsilon_-\ .\label{eq:coun1}
\eeqa
The triple $(O_q,\delta_{O_q},\epsilon_{O_q})$ is thus a counital right comodule algebra over $\Loop$.
\end{prop}
\begin{proof}
 Firstly, for the proof that $\delta_{O_q}$ preserves the $q$-Dolan-Grady relations, see~\cite[Prop.\,2.2]{BB09}. Then, using the definitions of the coproduct~\eqref{coprod} and counit~\eqref{counit} of $\Loop$, the conditions~\eqref{coassoc} and ~\eqref{right-counital} are readily checked. Finally,  the homomorphism property of $\epsilon_{O_q}$ is straightforward.
\end{proof}

\begin{rem}\label{rem:iotadef}
 We now compute the image of $W_0,W_1$ by the induced algebra map $\iota_{O_q}$  defined in~\eqref{iota-B}:
\beqa
\iota_{O_q}(W_0)&=&  k_+q\te_1 \tk_1 + k_-\tf_1 +  \varepsilon_+\tk_1 \in \Loop \ ,\label{eq:iot1}\\
	\iota_{O_q}(W_1)&=&  k_-q\te_0 \tk_0 + k_+\tf_0 + \varepsilon_-\tk_0 \in \Loop\ .\label{eq:iot2}\
\eeqa
This image generates a right coideal subalgebra in $\Loop$, according to Corollary~\ref{cor:im-iota}, and this is indeed a well known result, and moreover $\iota_{O_q}$ is injective in this case~\cite[Thm.\,7.4\,\&\,Ex.\,7.6]{Ko12}.
\end{rem}

\begin{defn}\label{def:stnst} The counital right comodule algebra $(O_q,\delta_{O_q},\epsilon_{O_q})$ associated with $(k_+,0)$ is called standard, whereas the one associated with $(k_+,\varepsilon_\pm)$ for $\varepsilon_\pm\neq 0$ is called non-standard.
\end{defn}

We now identify certain equivalent classes among the counital comodule algebras associated with $O_q$.

\begin{prop}\label{prop:isocomod} 
Let 
$(\tilde{O}_q,\tilde{\delta}_{\tilde{O}_q},\tilde{\epsilon}_{\tilde{O}_q})$ be the counital comodule algebra with coaction $\tilde{\delta}_{O_q}$ given by~\eqref{eq:del1} and~\eqref{eq:del2}
 with the substitutions $\rho \rightarrow \tilde{\rho}$, $k_\pm \rightarrow \tilde{k}_\pm$ and the counit $\tilde{\epsilon}_{O_q}$  given by~\eqref{eq:coun1} with $\varepsilon_\pm \rightarrow \tilde{\varepsilon}_\pm$. Then
$(O_q,\delta_{O_q},\epsilon_{O_q})$ and  $(\tilde{O}_q,\tilde{\delta}_{\tilde{O}_q},\tilde{\epsilon}_{\tilde{O}_q})$ are isomorphic counital comodule  algebras if 
\beqa
\tilde{k}_+ =\sqrt{\frac{\tilde{\rho}}{\rho}}k_+\ ,\quad \tilde{\varepsilon}_\pm =\sqrt{\frac{\tilde{\rho}}{\rho}}\varepsilon_\pm\ . \label{eq:idcomod}
\eeqa
\end{prop}
\begin{proof}
Using the isomorphism $\Psi\colon O_q \rightarrow \tilde{O}_q$ introduced in~\eqref{eq:Psi} together with the respective definitions of $\delta_{O_q}$, $\epsilon_{O_q}$, $\tilde{\delta}_{\tilde{O}_q}$, and $\tilde{\epsilon}_{\tilde{O}_q}$, equations~\eqref{eq:comod-alg-equiv} are checked provided the identification~\eqref{eq:idcomod} holds.
\end{proof}

\begin{conj}\label{conj:IsoCom}
With the notations from Proposition~\ref{prop:isocomod},
we conjecture a stronger statement: the counital right comodule algebras $(O_q,\delta_{O_q},\epsilon_{O_q})$ and  $(\tilde{O}_q,\tilde{\delta}_{\tilde{O}_q},\tilde{\epsilon}_{\tilde{O}_q})$ are isomorphic  if and only if \eqref{eq:idcomod} holds.
\end{conj}

\begin{rem}\label{rem:equiv-classes}
By Proposition~\ref{prop:isocomod}, the four-parameter family of counital comodule algebras associated with  $(\tilde{k}_+,\tilde{\varepsilon}_\pm)$ and arbitrary $\tilde{\rho}$ is in the same equivalent class as the  3-parameter family of counital comodule algebras associated with  $(k_+,\varepsilon_\pm)$ and $\rho$ fixed. It is thus sufficient to study the counital comodule algebra $O_q$ associated with $(k_+,\varepsilon_\pm)$ for $\rho$ fixed without loss of generality. 
\end{rem}

\section{Three presentations of $O_q$}\label{sec:Oqcomod}
From now on, and for the rest of the paper, we fix (recall Remark~\ref{rem:equiv-classes})
\begin{equation}
	{\mathbb K}={\mathbb C} \quad \mbox{and}\quad  \rho= - (q^2-q^{-2})^2 \ . 
	\label{rhofix}
\end{equation}
Besides Definition~\ref{def:Oq}, three other presentations for $O_q$ are known. They are given in terms of an infinite number of generators and relations, respectively called the {\it alternating} presentation of $O_q$ \cite{BasBel,Ter21b,LBG}, the {\it root vectors} presentation of $O_q$ \cite{BK17} and the {\it Drinfeld} type presentation of $O_q$ \cite{LW20}, mimicking Drinfeld's second realization of $U_q(\widehat{sl}_2)$ \cite{Dr0,Dr88}. Below, the three presentations are
respectively denoted $^{\{a\}}O_q$, $^{\{r\}}O_q$ and $^{Dr}O_q$, and briefly reviewed in Sections~\ref{subalt},~\ref{sec:rootpres},~\ref{sec:Drinpres}. The correspondence between the three presentations is discussed in Section~\ref{sec:Tercor}.

\subsection{The alternating presentation $^{\{a\}}O_q$}
\label{subalt}
This presentation is obtained from the alternating central extension of the $q$-Onsager algebra ${\cal A}_q$, which has been studied in details in~\cite{BasBel,Ter21,Ter21b,Ter21d}. We give the necessary material on $\cal A_q$ in Appendix~\ref{ap:AqRmat}. For example, see the defining relations in Definition~\ref{DefAq}  in terms of infinitely many generators $\{{\tW}_{-k},{\tW}_{k+1},$
		${\tG}_{k+1},{\tilde{\tG}}_{k+1}|k\in {\mathbb N}\}$, and its center ${\cal Z}$ in Proposition~\ref{DefZ}. We have an isomorphism of algebras \cite{Ter21,Ter21b}:
\beqa
{\cal A}_q \cong O_q \otimes {\cal Z}\ , \label{isoAqOqZ}
\eeqa
where ${\cal Z}$ is  generated by $\{\Delta_{k}|k\in {\mathbb N}^*\}$ defined via~\eqref{eq:ckp1} and~\eqref{deltau}. Therefore, $O_q$ is obtained from ${\cal A}_q$ by taking the image of  $\Delta_{k}$'s to some scalars denoted  $\delta_k$, i.e.\ a central reduction is applied to~${\cal A}_q$, as in \cite{BasBel}. We thus have  the surjective homomorphism
$\gamma: {\cal A}_q \rightarrow O_q$ 
such that
\begin{equation}\label{mapDeldel}
\begin{split}
&\gamma(\Delta_{k}) = 2\delta_{k}\ ,\qquad \delta_{k} \in {\mathbb C}\ ,\\
&W_{-k} := \gamma(\tW_{-k}) \ , \quad W_{k+1} :=\gamma(\tW_{k+1})\ , 
\quad G_{k+1}:= \gamma(\tG_{k+1})\ ,\quad {\tilde G}_{k+1}:=\gamma({\tilde \tG}_{k+1})\ ,
\end{split}
\end{equation}
where we also introduced notations for the image of generators.
We note  that the $\Delta_{k}$'s are algebraically independent \cite[Rem.\,8.17, Thm.\,10.2(iii)]{Ter21} and so all $\delta_{k}$ can be fixed independently.  Due to the isomorphism (\ref{isoAqOqZ}), any choice of scalars $\delta_{k}$  gives isomorphic quotients $O_q$. For our purpose, following~\cite[Ex.\,3.14]{BGL25} it will be sufficient to set the scalars as:
\beqa
\delta_{2p+2} &=&  (q+q^{-1})^{2}(q-q^{-1})(p+1) q^{-4p-4}c_{2p+2}^{-1}(\epsp^2+  \epsm^2) \ ,\qquad\label{delpair}\\
 \delta_{2p+1} &=&  (q+q^{-1})^{2}(q-q^{-1})(2p+1) q^{-4p-2}c_{2p+1}^{-1} \epsp\epsm \ , \quad p\in {\mathbb N}\ ,\label{delimp} 
\eeqa
where the $c_{k+1}$'s are given in (\ref{eq:ckp1}). 

\begin{defn}\label{def:altOq}
The alternating presentation of $O_q$, denoted by $^{\{a\}}O_q$, is defined via the image of $\gamma$, i.e.\ by the alternating generators 
$$\{W_{-k},W_{k+1},G_{k+1},{\tilde G}_{k+1}|k\in {\mathbb N}\}$$
and with $\delta_k$'s fixed as in~\eqref{delpair}~ and~\eqref{delimp}.
\end{defn}

\begin{prop}\label{def:Oqalt} 
Define the generating functions:
	\begin{align}
		{W}_+(u)=\sum_{k\in {\mathbb N}}{\normalfont W}_{-k}U^{-k-1} \ , \quad {W}_-(u)=\sum_{k\in  {\mathbb N}}{\normalfont W}_{k+1}U^{-k-1} \ ,\label{c1oq}\\
		\quad {G}_+(u)=\sum_{k\in {\mathbb N}}{\normalfont G}_{k+1}U^{-k-1} \ , \; \quad {G}_-(u)=\sum_{k\in {\mathbb N}}{\normalfont \tilde{{G}}_{k+1}}U^{-k-1} \label{c2oq} \  ,
	\end{align}
with 
$$
U=\frac{qu^2+q^{-1}u^{-2}}{q+q^{-1}}\ ,
$$
and with the convention that  rational functions in $U$ are considered  as power series in $u^{-2}$.
Then the defining relations of $^{\{a\}}O_q$ are given by (\ref{eq1-Aq})-(\ref{eq9-Aq}) for the substitutions ${\tW}_\pm(u) \mapsto W_\pm(u)$, ${\tG}_\pm(u) \mapsto G_\pm(u)$  and
	\beqa
	\gamma(\Delta(u))=  	2\frac{(q+q^{-1})^2(q-q^{-1})}{(q^2u^2-q^{-2}u^{-2})^2}\left( \epsp^2 + \epsm^2 + \epsp\epsm(q^2u^2+q^{-2}u^{-2}) \right)\ \label{2delust}
	\eeqa
	with (\ref{deltau}). 
\end{prop}	
\begin{proof}
We apply $\gamma$  to~\eqref{eq1-Aq}-\eqref{eq9-Aq} and \eqref{deltau} in order to get the defining relations of $^{\{a\}}O_q$. In particular, using~\eqref{mapDeldel} and the definitions~\eqref{delpair},~\eqref{delimp},  by straightforward computation we get the r.h.s.\,of~\eqref{2delust}.
\end{proof}
\begin{defn}\label{def:stalt}
The alternating presentation $^{\{a\}}O_q$ is called standard if $\varepsilon_\pm=0$, and non-standard otherwise. 
\end{defn}
\begin{rem}\label{rem:Oq-auto}
Recall the automorphism $\sigma$ and anti-automorphisms $\dagger$, $\tau$ of $\cal A_q$ introduced in Remark \ref{rem:auto} \cite[Lems.\,11.7, 11.8, 11.9]{Ter21b}. Because they all preserve the kernel of $\gamma$, $O_q$ admits analogous (anti-)automorphisms given by the same formulas in terms of the alternating generators, and we thus denote them by the same symbols. 
\end{rem}

Note that the expressions of the alternating generators of $^{\{a\}}O_q$ in terms of the generators $W_0,W_1$ of $O_q$ (see Definition \ref{def:Oq})  can be derived recursively using \cite[Prop.\,3.1]{BasBel}.
\begin{example}\label{ex1} For the standard case $\varepsilon_\pm=0$, besides $W_0,W_1$, the first generators read:
	\beqa
	\quad && \ {\tilde G}_{1} = \big[{W_0},W_1\big]_q \ ,\qquad 
	{W}_{-1} = \frac{1}{\rho}\left( (q^2+q^{-2})W_0W_1W_0 -W_0^2W_1 - W_1 W_0^2\right) + W_1 
	\nonumber\ 
	\eeqa
	and $G_1=\sigma(\tilde{G}_{1})$, ${W_2}=\sigma(W_{-1})$.
	\end{example}

\begin{rem}\label{rem:q-Ons-limit}
	For $q\rightarrow 1$ and $\varepsilon_\pm=0$, the relations of Definition \ref{def:Oqalt} specialize to the defining relations of the alternating presentation of the Onsager algebra \cite[Def.\,4.1 \& Thm.\,2]{BC18}.  Indeed, first  let us fix $a=\textsf{i} (q^2-q^{-2})/4$ and renormalize the currents as $W_{\pm}(u)\mapsto a W_{\pm}(u)$, $G_{\pm}(u)\mapsto a^2 G_{\pm}(u)$. Then, from the relation \eqref{2delust} with $\varepsilon_\pm=0$ we obtain $G_+(u) = -G_-(u)$ in the limit, and secondly, the leading terms of  the relations (\ref{eq1-Aq})-(\ref{eq9-Aq}) give in the limit $q\to 1$ the classical Onsager current relations established in the proof of~\cite[Prop.\,4.1]{BC18}.
\end{rem}

In Appendix~\ref{ap:AqRmat}, we recall a PBW basis for ${\cal A}_q$ associated with the ordering~\eqref{eq:order}.  Using the surjective map~$\gamma$ with (\ref{mapDeldel}), and the fact~\cite[Prop.\,3.1]{BasBel} that in $O_q$, all $G_{k+1}$'s are expressed in terms of the alternating generators 
\beqa\label{eq:alt-gen-Oq}
\{W_{-\ell},W_{\ell+1},\tilde G_{\ell+1}|\ell \in{\mathbb N}\}\ ,
\eeqa
 a PBW basis for $O_q$ has been conjectured~\cite[Conj.\,1]{BasBel}, \cite[Conj.\,16.2]{Ter21b} and proven in~\cite[Thm.\,7.1]{Zh26b}
in any linear order $<$ such that
\begin{align}
	W_{-i} < W_{j+1} < {\tilde G}_{k+1}\ , \qquad \qquad i,j,k \in \mathbb N.
	\label{eq:orderOq}
\end{align}

\subsection{The root vectors presentation $^{\{r\}}O_q$}\label{sec:rootpres}
The 
Onsager algebra  has been introduced in~\cite{Ons44}, in terms of  generators $\{A_n,G_m\,|\,n\in {\mathbb Z}, m\in\N^*\}$.  
 In \cite{BK17}, their quantum group analogs -- the so-called root vectors of $O_q$ -- denoted $\{B_{n\delta+\alpha_1}, B_{m \delta} \, | \, n\in {\mathbb Z}, m\in \N^*\}$ are constructed via a Lusztig's like braid group action, see~\cite[Sec.\,2.3 \&\,3]{BK17} for details. 
 
We now briefly recall the construction of  $\{B_{n\delta+\alpha_1}, B_{m \delta} \, | \, n\in {\mathbb Z}, m\in \N^*\}$. Let $\alpha_0$ and $\alpha_1$ denote the simple roots for $\widehat{sl}_2$ and set $\delta=\alpha_0+\alpha_1$. Let us denote the root system of $\widehat{sl}_2$ by
 \beqa
 \roots=\{n\delta+\alpha_0,n\delta+\alpha_1, m\delta\,|\,n\in {\mathbb Z}, m\in {\mathbb Z}^*\} \ ,\label{rootsys}
 \eeqa
 and collect all the real roots into
 $$
 \roots^{\real}=\{n\delta\pm\alpha_1\,|\,n\in {\mathbb Z}\}\ .
 $$
Define 
 $B_{\alpha_0},B_{\alpha_1}$ that satisfy the $q$-Dolan-Grady relations (\ref{qDG1}), (\ref{qDG2}) with the identification (\ref{rhofix}) and
 	\beqa
 	B_{\alpha_0}\equiv W_0\ ,\qquad B_{\alpha_1}\equiv W_1\ . \nonumber
 	\eeqa
 Let $\Phi: O_q \rightarrow O_q$ denote the algebra automorphism defined by $\Phi(B_{\alpha_0})=B_{\alpha_1}$ and $\Phi(B_{\alpha_1})=B_{\alpha_0}$. Consider the braid group with two (Lusztig's type) algebra automorphisms $T_0,T_1$ with $T_1=\Phi \circ T_0 \circ \Phi$. The action of $T_0:O_q \rightarrow O_q$ is given in~\cite[Prop.\,2.2]{BK17}. The root vectors 
 of $O_q$ are denoted $B_\gamma$ with $\gamma\in \roots^{\real}$
 for all real roots and  $\gamma\in\{n\delta|n\in \N^*\}$ for all imaginary roots. The real root vectors are defined as~\cite[Eq.\,(3.7)]{BK17}:
\beqa
B_{n\delta + \alpha_0} = (T_0\Phi)^n(B_{\alpha_0}) \qquad \mbox{and}\qquad B_{n\delta + \alpha_1} = (T_0\Phi)^{-n}(B_{\alpha_1})\ ,\label{def:Be}
\eeqa
and the imaginary root vectors are~\cite[Eq.\,(3.11)]{BK17}:
\beqa
 B_{m\delta} = -q^{-1}[B_{\alpha_0}, B_{(m-1)\delta+\alpha_1}]_q + (q^{-2}-1)\sum_{\ell=0}^{m-2} B_{\ell\delta+\alpha_1} B_{(m-\ell-2)\ell\delta+\alpha_1}\ .\label{def:Bim}
\eeqa

Similarly to Example \ref{ex1}, the root vectors can be expressed as polynomials of the fundamental generators $B_{\alpha_0},B_{\alpha_1}$ using~\eqref{def:Bim} and
the recursion relation~\cite[Eq.\,(3.10)]{BK17}:
\beqa
B_{(n+1)\delta+\alpha_1} = B_{(n-1)\delta+\alpha_1} + \frac{q(q+q^{-1})}{\rho}[B_\delta,B_{n\delta+\alpha_1}]\ .\label{eq:recBe}
\eeqa
\begin{example} The first root vectors are given by: 
	\beqa
B_\delta&=& q^{-2}B_{\alpha_1}B_{\alpha_0} - B_{\alpha_0}B_{\alpha_1}\ ,\nonumber\\
B_{\delta+\alpha_1}&=& \frac{q(q+q^{-1})}{\rho}\left((1+q^{-2}) B_{\alpha_1}B_{\alpha_0}B_{\alpha_1} - B_{\alpha_0}B_{\alpha_1}^2 -q^{-2}B_{\alpha_1}^2B_{\alpha_0} \right) +B_{\alpha_0} \
 \nonumber 
\eeqa
and $B_{\delta+\alpha_0}=\tau(B_{\delta+\alpha_1})$, recall the anti-automorphism $\tau$ from Remark~\ref{rem:Oq-auto}.
They clearly differ from the expressions for the  alternating generators of Example~\ref{ex1}.
\end{example}

By construction, from~\eqref{def:Be} we have:
\beqa 
B_{-(n+1)\delta+\alpha_1}=B_{n\delta+\alpha_0}\ .\label{eq:BmBn}
\eeqa
Thus, $B_\gamma=B_{-\gamma}$ for all $\gamma\in\roots^{\mathrm{re}}$.  Now, recall the root system of $\widehat{sl}_2$ given by~\eqref{rootsys} 
and set 
\beqa
\Rbar=\roots/\sim \label{eq:Rtilde}
\eeqa
 where $\sim$ is the equivalence relation on  $\roots$ given by $\alpha \sim \beta \Leftrightarrow \alpha=\pm \beta$.
It follows that 
 $({\mathbb Z}\delta+\alpha_1)\cup \N^*\delta$  is a representative set of $\Rbar$. With respect to the uniquely determined total ordering $<$ on $\Rbar$ defined by
 \begin{align}\label{eq:order-def}
 	(m+1)\delta+\alpha_1<m\delta+\alpha_1<k\delta<(k+1)\delta \quad \mbox{for all $m\in {\mathbb Z}$, $k\in\N^*$},
 \end{align}  
a PBW basis for $O_q$ is obtained \cite[Thm.\,I]{BK17}. We note that other orderings may be chosen in constructing PBW bases, as we do below in the proofs of Propositions~\ref{prop:FRTOqrv} and~\ref{prop:Kr-inv},  see the discussion in \cite[Sec.\,5.1]{BK17} and the proof in~\cite[Thm.\,5.3]{Zh26b}.
According to the ordering~\eqref{eq:order-def}, the $q$-analog of Onsager's relations are given in~\cite[Thm.\,II]{BK17}:\footnote{The $q$-commutator $[X,Y]_{q^{-2}}$ from \cite{BK17} equals $q^{-1}[X,Y]_q$ defined by (\ref{defcom}); and their scalar $c = -\rho q^{-1}/(q+q^{-1})^2$.}

\begin{thm}[\cite{BK17}]\label{thm:BK17}
  The q-Onsager algebra admits a root vectors presentation  $^{\{r\}}O_q$ in the sense that its generators are the root vectors $\{B_{n\delta+\alpha_1}, B_{m \delta} \, | \, n\in {\mathbb Z}, m\in \N^*\}$ satisfying the defining relations:
		\begin{align}
			& [B_{m\delta},B_{n\delta}]=0\ ,   \label{eq:qOnsCom1}\\
			& [B_{r\delta+\alpha_1}, B_{(r+m)\delta+\alpha_1}]_{q}=-qB_{m\delta}+ qC^{\real}_{r,m} \ , \label{eq:qOnsCom2}\\
			& [B_{m\delta},B_{r\delta+\alpha_1}] = -\frac{\rho q^{-1}}{q+q^{-1}} \big(q^{-2(m-1)} B_{(r-m)\delta+\alpha_1} 
			- q^{2(m-1)}B_{(r+m)\delta+\alpha_1}\big) + C^{\ima}_{r,m}\ , \label{eq:qOnsCom3}
		\end{align}
		where $r\in {\mathbb Z}$ and $C^{\real}_{r,m}$ and $C^{\ima}_{r,m}$  are  linear combinations of elements of the PBW basis for the ordering (\ref{eq:order-def}), with coefficients in $(q-1)\big({\mathbb Z}[q,q^{-1}]+ {\mathbb Z}[q,q^{-1}]\rho\big)$.  
\end{thm}
The
elements $C^{\real}_{r,m}$ and $C^{\ima}_{r,m}$ are given explicitly in~\cite[Props.\,5.5 \&~5.8]{BK17}.

\begin{rem}
	For $q\rightarrow 1$, up to a rescaling of the root vectors, the relations \eqref{eq:qOnsCom1},  \eqref{eq:qOnsCom2}, and \eqref{eq:qOnsCom3} specialize to the Onsager relations \cite[Sec.\,3.3]{BK17}.  
\end{rem}

\begin{rem}\label{rem:Oq-grading}
    By Theorem~\ref{thm:BK17}, using the explicit expressions for $C^{\real}_{r,m}$ and $C^{\ima}_{r,m}$ in~\cite[Props.\,5.5 \&~5.8]{BK17}, we note that there exists a natural grading on $O_q$ given by $\deg^{O_q}: O_q \rightarrow {\mathbb Z}_2$ 
    on the root vectors:
\beqa
\deg^{O_q}(B_{n\delta+\alpha_i})=1 \quad \mbox{and}\quad  \deg^{O_q}(B_{m\delta})=0\ , \quad i=0,1,\ n\in{\mathbb N},\ m\in{\mathbb N}^*\ .\label{eq:degOq}
\eeqa
\end{rem}

With respect to the root vectors presentation of $O_q$ recalled in Theorem \ref{thm:BK17}, the defining relations (\ref{eq:qOnsCom1})-(\ref{eq:qOnsCom3}) can be presented in terms of generating functions. Following~\cite[Def.\,5.2]{Ter21b}, define:
\beqa
B^+(z) = \sum_{n\in{\mathbb N}} B_{n\delta + \alpha_1}z^{n}\ , \quad B^-(z) = \sum_{n\in{\mathbb N}} B_{n\delta + \alpha_0} z^{n}\ , \quad B(z) = \sum_{n\in{\mathbb N}} B_{n\delta} z^{n}\ , \label{rvgen1}
\eeqa
with the convention $B_{0\delta}:=q^{-2}-1$ and indeterminate $z$. The relations take the form~\cite[Props.\,5.8 \&~5.9]{Ter21b}:
\begin{align}
0&=\big[B(z),B(w)\big]\ ,\label{eq:cc0}\\
	\begin{split}
		0&=  (qz-q^{-1} w)B^\pm(z) B^\pm(w) +(qw-q^{-1}z)B^\pm(w) B^\pm(z)
		-(q-q^{-1}) w \bigl( B^\pm(w)\bigr)^2
		\\
		&\quad -(q-q^{-1}) z \bigl( B^\pm(z)\bigr)^2
		+ \frac{q(z-w)}{1-zw} B(w) + \frac{q(w-z)}{1-zw} B(z)\ ,
	\end{split} \label{eq:cc1}
	\\
	\begin{split}
		0 &= (1-q^{-2} z w)B^-(z) B^+(w) + (zw-q^{-2})B^+(w) B^-(z)
		+(1-q^{-2}) w\bigl(B^+(w)\bigr)^2 
		\\
		&\quad + (1-q^{-2}) z\bigl(B^-(z)\bigr)^2  + \frac{1-zw}{z-w} B(z) - \frac{1-zw}{z-w} B(w)\ ,
	\end{split}\label{eq:cc2}
	\\
\begin{split}
	0 &= 
	B(z)B^+(w) (1-q^{-2}zw)(w-q^2z)
	-
	B^+(w)B(z)(1-q^2zw)(w-q^{-2}z)
	\\
	&\quad +
	B(z)B^+(q^{-2}z)(q^2
	-q^{-2})(1-q^{-2}zw)z
	+
	B^-(q^{-2}z)B(z)(q^2 -q^{-2})(w-q^{-2}z)z\ ,
	\end{split}  \label{eq:cc4}\\
	\begin{split}
	0 &= 
	B^-(w)B(z)(1-q^{-2}zw)(w-q^2z)
	-
	B(z)B^-(w)(1-q^2zw)(w-q^{-2}z) 
	\\
	&\quad +
	B^-(q^{-2}z)B(z)(q^2
	-q^{-2})(1-q^{-2}zw)z
	+
	B(z)B^+(q^{-2}z)(q^2 -q^{-2})(w-q^{-2}z)z\ .
	\end{split} \label{eq:cc5}
\end{align}

\subsection{Terwilliger's correspondence}\label{sec:Tercor}
The precise relationship between the alternating generators of $^{\{a\}}O_q$ for the {\it standard} case (see Definition~\ref{def:stalt}) 
 and the root vectors
 of $^{\{r\}}O_q$ is given in~\cite[Secs.\,11\,\&\,12]{Ter21b}. In this section, we recall the main points that lead to this correspondence, and extend it to the {\it non-standard} case as well; See Lemmas~\ref{lem:GB},~\ref{lem:WB}, and~\ref{lem:BW}, with Definition~\ref{def:Bbar}. This is done by applying the map $\gamma$ in~\eqref{mapDeldel} to the results in~\cite[Sec.\,9]{Ter21b} adapted to our conventions. Meanwhile, we introduce an `auxiliary' generating function ${\bar B}(z)$ that will play a key role in further sections.

Recall the generating functions $B^\pm(z),B(z)$ for the root vectors given by~\eqref{rvgen1}. Firstly, in \cite[Sec.\,5]{Ter21b} the three relations (\ref{eq:qOnsCom1})-(\ref{eq:qOnsCom3}) are transformed into the set of `half-currents'  relations~\eqref{eq:cc1}-\eqref{eq:cc5}. Under the injective homomorphism $\imath: O_q \rightarrow {\cal A}_q$ given by~\cite[Lem.\,6.8]{Ter21b}:
\beqa
\imath(W_0) = \tW_0 \qquad \mbox{and} \qquad \imath(W_1) = \tW_1\ ,
\eeqa
it was shown that $\imath(B^\pm(z))$ and $\imath(B(z))$ can be expressed in terms of the  ${\cal A}_q$ generating functions ${\tW}_\pm(u)$, ${\tG}_\pm(u)$ introduced in Theorem~\ref{def:Aq0}, see~\cite[Props.\,9.6,~9.7\,\&\,Cor.\,10.5]{Ter21b}.
In particular, a {\it factorization formula}~\cite[Thm.\,10.1 \& 
Cor.\,10.5]{Ter21b} relates
the image via  $\imath$ of the imaginary root vectors of $O_q$ to the alternating generators of $\cal A_q$. We recall the quantum determinant $\Gamma(u)$ for $\cal A_q$ given by  (\ref{gamma}) together with the conventions (\ref{gamZvee}), and that  ${\tG}_-(u) + \rho/(q-q^{-1})$ is invertible~\cite[Lem.\,4.1]{Ter21d} because of the commutativity relations~\eqref{eq4-Aq}.
For the choice of $\rho$ given by (\ref{rhofix}), we then have the  factorization formula
\beqa
\imath(B(u^{-2})) = \xi^{-1}  \left({\tG}_-(u) + \frac{\rho}{q-q^{-1}}\right)^{-1} \left({\tG}_-(uq^{-1}) + \frac{\rho}{q-q^{-1}}\right)^{-1} 	\frac{\Gamma(uq^{-1})}{(q+q^{-1})^2(u^2-u^{-2}) } \ \label{BinGGZ}
\eeqa
with  
\beqa
\xi = -q(q-q^{-1})\rho^{-2}\ .\label{xi}
\eeqa

Then, we study the image via $\gamma$ of~\eqref{BinGGZ} in $O_q$. Using~\eqref{2delust}, let us define:
\beqa
d(u)^{(\epsp,\epsm)}_{\rho}&:=&\gamma(\Gamma(u))\ \label{eq:durho}\\
&=& 	\frac{(q+q^{-1})^2}{(q^2u^2-q^{-2}u^{-2})}\left( \epsp^2 + \epsm^2 + \epsp\epsm(q^2u^2+q^{-2}u^{-2}) - \frac{\rho(q^2u^2-q^{-2}u^{-2})^2}{(q^2-q^{-2})^2}\right)\ .\nonumber
\eeqa 
Applying $\gamma$ from the left to~\eqref{BinGGZ} and using  $\gamma \circ \imath=\id_{O_q}$, we get a new relation in $O_q[[u^{-2}]]$ that suggests to introduce a new generating function, defined as follows.
\begin{defn}\label{def:Bbar}
The `auxiliary' generating function $\bar B(z)$ is defined by
\beqa
{\bar B}(u^{-2}) :=  {G}_-(u) + \frac{\rho}{q-q^{-1}} \ \quad \mbox{with}\qquad {\bar B}(z) = \sum_{n=0}^{\infty} {\bar B}_{n}z^{n}\ , \qquad {\bar B}_{0}:=\rho/(q-q^{-1})\ ,\label{Bbargen}
\eeqa
where the coefficients $\{\bar B_n\}_{n\in {\mathbb N}}$ are such that the following relation holds in $O_q$, recall $\xi$ in~\eqref{xi}: 
\beqa
B(u^{-2}) = \xi^{-1}  {\bar B}(u^{-2})^{-1} {\bar B}(q^2u^{-2})^{-1} 	\frac{d(uq^{-1})^{(\epsp,\epsm)}_{\rho}}{(q+q^{-1})^2(u^2-u^{-2}) } \ .\label{BinGGred}
\eeqa
\end{defn}
Note that invertibility of $\bar B(u^{-2})$ follows from the fact that $G_-(u) + \rho/(q-q^{-1})$ is invertible, see the comment above~\eqref{BinGGZ}.
\begin{lem}\label{lem:invertBbar}
	The inverse generating function  $\bar B(u^{-2})^{-1}$ is given by:
	\beqa
	\bar B(u^{-2})^{-1} = \sum_{k=0}^\infty u^{-2k} {\bar B}^{inv}_{k} \ ,\qquad {\bar B}^{inv}_{0} = (q-q^{-1})/q \ ,\label{expBinv}
	\eeqa
	where the coefficients $\{ {\bar B}^{inv}_{k}\}_{k\in{\mathbb N}^*} \in O_q$ are determined by the recurrence relations:
	\beqa
	\sum_{k=0}^{p} {\bar B}_{p-k} {\bar B}^{inv}_{k} =0  \ ,\qquad p \in {\mathbb N}^*\ .\label{eq:Bbarinv}
	\eeqa
\end{lem}
\begin{proof}
The relation $\bar B(u^{-2})\bar B(u^{-2})^{-1}=1$ is equivalent to   $ {\bar B}^{inv}_{0} = (q-q^{-1})/\rho$ and the recurrence relations~\eqref{eq:Bbarinv}.
\end{proof}

\begin{rem}\label{rem:PTmap}
The map $\gamma$ given in~\eqref{mapDeldel} with the definition of $\delta_k$'s in \eqref{delpair}-\eqref{delimp}, is a generalization of the map $\gamma$ introduced in~\cite[Lem.\,11.2]{Ter21b}. Indeed, the latter maps all central elements $\Delta_k$ to $0$ -- compare the conventions in~\cite{Ter21b} versus our conventions, as given in Remark~\ref{rem:convention}. Thus, it is associated with the standard case $\varepsilon_\pm=0$ according to \eqref{delpair},~\eqref{delimp}. This is also seen by comparing~\eqref{BinGGred} to~\cite[Eq.\,(131)]{Ter21b}.
\end{rem}

We now explain how the coefficients $\{\bar B_{n}\}_{n\in {\mathbb N}}$ are determined in terms of $B_{n\delta}$ explicitly. Applying $\gamma$ to $[\tG_-(u),\tG_-(v)]=0$ in~(\ref{eq4-Aq}) and using~\eqref{BinGGred} it follows respectively:
\beqa
[{\bar B}_{m},{\bar B}_n]=0\qquad \mbox{and}\qquad [B_{m\delta},{\bar B}_n]=0  \ .\label{eq:comBbarBdel}
\eeqa
 In Appendix~\ref{ap:Bbarsol}, the set of recursive relations that determine the coefficients  $\bar B_{n}$ in terms of $B_{n\delta}$ is derived from~\eqref{BinGGred}-\eqref{eq:comBbarBdel}, see Proposition~\ref{prop:Bbar}.
\begin{example}\label{ex:Bbarexplicit} For any $\varepsilon_\pm\in{\mathbb C}^*$,
	\beqa
	&&	 {\bar B}_1 =-\left(q +q^{-1}\right) B_\delta +\frac{ q^{-1}\rho \epsp\epsm}{\left(q^{2}-q^{-2}\right)}\ ,\nonumber\\
	&& {\bar B}_2 =    -\frac{(q+q^{-1})^2}{q(q^2+q^{-2})}B_{2\delta} 
	-\frac{(q^3-q^{-3})}{(q^2+q^{-2})(q-q^{-1})^2}B_{\delta}^2   \nonumber\\ &&\qquad \quad - q^{-1} \epsp\epsm  B_{\delta}  -\frac{q^{-2}(q-q^{-1})}{(q^2+q^{-2})}\left(  (q+q^{-1})^2(\epsp^2+\epsm^2)  - \epsp^2\epsm^2 \right)\ .\nonumber
\eeqa
\end{example}
\begin{rem} \label{rem:Oq-gradingbar}  
With respect to the ${\mathbb Z}_2$-grading of $O_q$ introduced in Remark~\ref{rem:Oq-grading},
it follows:
\beqa
 \deg^{O_q}(\bar B_{m})=0 \ ,\quad m\in{\mathbb N}^*\ .\label{eq:degBbar}
\eeqa
\end{rem}

We can now give a complete correspondence between the generating functions for $^{\{a\}}O_q$ from~\eqref{c1oq}, \eqref{c2oq} and those for $^{\{r\}}O_q$ from \eqref{rvgen1} and \eqref{Bbargen}. In the proofs of the next three lemmas we use a dictionary between our notations and those of~\cite{Ter21b} explained in Remark~\ref{rem:convention}.

\begin{lem}\label{lem:GB}  The following relations hold in $O_q$:
	\beqa
	G_+(u) + \frac{\rho}{q-q^{-1}} &=& {\bar B}(u^{-2})\left( \frac{B(u^{-2})}{q^{-2}-1} - \frac{(q-q^{-1})(B^{+}(u^{-2}))^2 + q^{-1}(q^2u^2-q^{-2}u^{-2})B^{-}(q^{-2}u^{-2})B^{+}(u^{-2}) }{u^2(qu^2-q^{-1}u^{-2})} \right),\nonumber\\
	G_+(u) + \frac{\rho}{q-q^{-1}} &=&  \left( \frac{B(q^{-2}u^{-2})}{q^{-2}-1} + \frac{(q-q^{-1})(B^{+}(q^{-2}u^{-2}))^2 - q(u^2-u^{-2})B^{+}(q^{-2}u^{-2})B^{-}(u^{-2})}{q^{2}u^{2}(qu^2-q^{-1}u^{-2})} \right) {\bar B}(u^{-2}).\nonumber 
	\eeqa
\end{lem}
\begin{proof}
Apply the algebra map $\gamma$ given by~\eqref{mapDeldel} to the adapted (according to Remark~\ref{rem:convention}) version of the third and second equations in~\cite[Prop.\,9.11]{Ter21b}, respectively.  Then, we use the first equality in~\eqref{Bbargen}.
\end{proof}
 A comparison of the power series expansion on both sides of the equations in~\eqref{Bbargen} and Lemma~\ref{lem:GB} gives the following results.
\begin{example} For the standard case $\varepsilon_\pm=0$, the correspondence for the first generators is:
	\beqa
&&	{\tilde G}_1= -qB_{\delta}\ ,\qquad 
{\tilde G}_2=   -\frac{q}{(q^2+q^{-2})} B_{2\delta} - \frac{q^2(q^3-q^{-3})}{(q^2+q^{-2})(q^2-q^{-2})^2}B_{\delta}^2 \ ,\nonumber\\
&&	{G}_1= q^3B_{\delta} + (q^3 - q^{-1})B_{\alpha_0}B_{\alpha_1}\ ,\nonumber\\
&& 	{G}_2= \frac{q^5}{q^2+q^{-2}}B_{2\delta} + \frac{q^6(q-q^{-1})}{(q^2-q^{-2})^2(q^2+q^{-2})}B_{\delta}^2 + (q-q^{-1})B_{\alpha_0}B_{\delta+\alpha_1} + (q^3-q)B_{\delta+\alpha_0}B_{\alpha_1} \nonumber\\ && \qquad + \frac{q^2}{q+q^{-1}}B_{\alpha_0}B_{\alpha_1}B_{\delta} + q(q-q^{-1})^2  B_{\alpha_0}^2 \   .\nonumber
	\eeqa
\end{example}
\begin{lem}\label{lem:WB} 
The following relations hold in $O_q$:
	\beqa
	W_+(u)&=& -\frac{q^{-2}u^{-2}B^{+}(q^{-2}u^{-2}) + B^{-}(u^{-2})}{(q^2-q^{-2})(qu^2-q^{-1}u^{-2})}{\bar B}(u^{-2})\ ,\label{WB1}\\
	W_-(u)&=& -\frac{ B^{+}(q^{-2}u^{-2}) + u^{-2}B^{-}(u^{-2})}{(q^2-q^{-2})(qu^2-q^{-1}u^{-2})}{\bar B}(u^{-2})\ ,\label{WB2}
	\eeqa
	and
	\beqa
	W_+(u)&=& -{\bar B}(u^{-2})\frac{u^{-2}B^{+}(u^{-2}) + B^{-}(q^{-2}u^{-2})}{(q^2-q^{-2})(qu^2-q^{-1}u^{-2})}\ ,\label{WB3}\\
	W_-(u)&=& - {\bar B}(u^{-2})\frac{ B^{+}(u^{-2}) + q^{-2}u^{-2}B^{-}(q^{-2}u^{-2})}{(q^2-q^{-2})(qu^2-q^{-1}u^{-2})}\ ,\label{WB4}
	\eeqa
\end{lem}
\begin{proof}
Apply $\gamma$ to the set of equalities in~\cite[Thms~9.8,~9.9]{Ter21b}, adapted to our conventions.
\end{proof}
\begin{example} For the standard case $\varepsilon_\pm=0$, the correspondence for the first generators is:
	\beqa
	&& W_0=	B_{\alpha_0}\ ,\quad W_1 = B_{\alpha_1}\ ,\nonumber\\
	 	 && W_{-1}= \frac{1}{(q+q^{-1})(q^2-q^{-2})}\left( qB_{\alpha_0}B_\delta + (q^2-q^{-2})(qB_{\delta+\alpha_0} + q^{-1}B_{\alpha_1})  \right)	 \ ,\nonumber\\
	 	 && W_{2}= \frac{1}{(q+q^{-1})(q^2-q^{-2})}\left( qB_{\alpha_1}B_\delta + (q^2-q^{-2})(q^{-1}B_{\delta+\alpha_1} + qB_{\alpha_0})  \right) \ .\nonumber
	\eeqa
\end{example}
\begin{lem}\label{lem:BW} 
	 The following relations hold in $O_q$:
	\beqa
	B^+(u^{-2})&=& (q^2-q^{-2})\Big( G_-(u) + \frac{\rho}{q-q^{-1}} \Big)^{-1}\left( q^{-1}W_+(u) - qu^2W_-(u) \right) \ ,\label{BWrel1}\\
	B^-(u^{-2})&=& (q^2-q^{-2})\left( q^{-1}W_-(u) - qu^2W_+(u) \right)\Big( G_-(u) + \frac{\rho}{q-q^{-1}}\Big)^{-1}\  .\label{BWrel2}
	\eeqa
\end{lem}
\begin{proof} Apply $\gamma$ to  the adapted version of the second relation in~\cite[Prop.\,9.7]{Ter21b} and first relation in~\cite[Prop.\,9.6]{Ter21b}, respectively.
\end{proof}
\begin{rem}
	For the Onsager algebra, the generators $\{A_n,G_m|n\in {\mathbb Z},m \in {\mathbb N}^*\}$ \cite{Ons44} are the classical analogs of the root vectors $B_{(n-1)\delta+ \alpha_1}$ and $B_{m\delta}$, respectively~\cite{BK17}.  The precise relationship between the Onsager's generators $\{A_n,G_m|n\in {\mathbb Z},m \in {\mathbb N}^*\}$ and the classical alternating generators  $\{w_{-k},w_{k+1},{\tilde g}_{k+1}|k\in {\mathbb N}\}$, analogous to Lemmas~\ref{lem:GB} and~\ref{lem:BW}, is given by \cite[Eq.\,(4.8)]{BC18}. 
\end{rem}

\subsection{The Drinfeld type presentation $^{Dr}O_q$}\label{sec:Drinpres}
In~\cite{LW20}, the authors proposed a Drinfeld type presentation for $O_q$ that we denote  $^{Dr}O_q$.
While several equivalent formulations of $^{Dr}O_q$ were given, we  consider the most compact one \cite[Thm.\,5.1]{LW20} in terms of two generating functions only, denoted respectively ${\bf B}(z)$ for {\it real} root vectors and ${\bf \Theta}(z)$ for {\it imaginary} root vectors. 

In the present paper, we essentially follow the notations in~\cite{BK17,Ter21b} that slightly differ from the ones in~\cite{LW20,P24}. The correspondence between both notations will be needed in Section~\ref{sec:DrOq}, displayed in the following Table.

	\bigskip
{	\begin{center}
					\begin{tabular}[t]{cc}
			{\rm Notations in \cite{Ter21b}} & {\rm Notations in \cite{LW20,P24}}
			\\
			\hline
			$q$ & $v$ 
			\\
			$W_0$ & $v^{1/2} (v-v^{-1})B_0$
			\\
			$W_1$ & $v^{1/2}(v-v^{-1})B_1$
			\\
			 $B_{n\delta+\alpha_0} $ &  \qquad $ v^{1/2}(v-v^{-1}) B_{1,-(n+1)}$ 
			\\
			$B_{n\delta+\alpha_1} $ & $ v^{1/2}(v-v^{-1}) B_{1,n}$
			\\
			$B_{n\delta} $ & $ -v^{-1}(v-v^{-1})^2 \acute{\Theta}_n$
			\\
			$1$, $1$, $1$ & $\mathbb K_{n\delta+\alpha_1}$, $\mathbb K_\delta$, $C$
			\\
			$\lbrack X, Y\rbrack_q$ & $v \lbrack X,Y\rbrack_{v^{-2}}$
	\end{tabular}
\captionof{table}{Notations}
\label{tab:Nota}
\end{center}}
\bigskip
	
Using  this correspondence, following \cite{LW20}  let us introduce new  generating functions. Recall the real and imaginary `half-currents' $B^\pm(z),B(z)$ given by (\ref{rvgen1}). We define:
\beqa
{\bf B}(z) &=&  B^+(z) + z^{-1}B^{-}(z^{-1}) \ \qquad \mbox{(corresponding to $\sum_{n\in {\mathbb Z}} B_{1,n} z^n$ in  \cite{LW20})}\ , \label{currentB}\\
{\bf \Theta}(z) &=& \frac{(1-z^2)}{(1-q^{-2}z^2)} B(z) \ \qquad \qquad \mbox{(corresponding to \cite[eq.\,(2.31)]{LW20})}\ ,\label{currentTheta}\\
D(z)& =& \sum_{k\in {\mathbb Z}}  z^k\ .\nonumber 
\eeqa

By definition~\eqref{currentTheta}, with \eqref{rvgen1}, ${\bf \Theta}(z)$ is a power series in the generators $\{B_{n\delta}\}_{n\in{\mathbb N}^*}$. In~\cite{LW20}, for convenience the authors further introduce new elements $\{{\acute H}_m\}_{m\in{\mathbb N}^*}$ that are related with the $B_{n\delta}$'s using the following power series identity.
\begin{defn}\label{def:Hacute} The elements ${\acute H}_m$  are defined as: 
\beqa
{\bf \Theta} (z)= -q^{-1}(q-q^{-1}) \frac{(1-z^2)}{(1-q^{-2}z^2)}  \exp\left((q-q^{-1}) \sum_{m=1}^{\infty} {\acute H}_m z^m\right)\ .\label{eq:ThetaHmbis}
\eeqa
\end{defn}
\begin{rem}
The relation~\eqref{eq:ThetaHmbis} is mimicking the relation in Drinfeld's second realization $^{Dr}\Loop$ between the generating functions $\Psig(z)$, $\Phig(z)$ and exponentials of the power series in the Drinfeld generators $\{\tho_m\}$, $\{\tho_{-m}\}$ for $m\in{\mathbb N}^*$, respectively, see in further section~\eqref{def:psi},~\eqref{def:phi}.
\end{rem}
\begin{prop}[\cite{LW20}]\label{prop:DrOq}  $^{Dr}O_q$ is generated by $\{B_{n\delta + \alpha_1},{\acute H}_m|n \in {\mathbb Z},m\in {\mathbb N}^*\}$,  with the defining relations:
\begin{align}
	{\bf \Theta} (z)  {\bf \Theta}(w)
	&= {\bf \Theta} (w)  {\bf \Theta}(z)\ ,
	\label{gf1} \\
	{\bf \Theta} (z)  {\bf B}(w)
	& =  \frac{(1 -q^{-2}zw^{-1}) (1 -q^{2} zw)}{(1 -q^{2}zw^{-1})(1 -q^{-2}zw)}
	{\bf B}(w) {\bf \Theta} (z)\ ,
	\\
	(q^2z-w) {\bf B}(z) {\bf B}(w) & +(q^{2}w-z) {\bf B}(w) {\bf B}(z)
	\label{gf3} \\
	& = - D(zw) \big( (q^2z-w){\bf \Theta}(w) +(q^2w-z){\bf \Theta}(z) \big)\ .\notag
\end{align}
\end{prop}

Let us notice that the extra central elements ${\mathbb K}_1$ and $C$ in~\cite[Eqs.\,(1.1)-(1.3)]{LW20} can be included though the substitutions $\bold B(z)\rightarrow{\mathbb K}_1^{-1/2}\bold B(z)$  and then  $z\rightarrow C^{1/2}z$ and $w\rightarrow C^{1/2}w$ in above defining relations. For simplicity, it is thus sufficient to set ${\mathbb K}_1=C=1$.
The `half-currents' and `full-currents' formulations~\eqref{eq:cc0}-\eqref{eq:cc5} and~\eqref{gf1}-\eqref{gf3}, respectively,  are equivalent to the relations in Theorem \ref{thm:BK17}, see \cite[Section 5]{Ter21b} and \cite[Section 2.2]{LW20} for a proof, respectively.

To conclude this section, for further analysis we show how ${\bar B}(z)$ and its inverse relate to the Drinfeld type presentation $^{Dr}O_q$ in Proposition~\ref{prop:DrOq}.
Recall that the auxiliary generating function $\bar B(z)$ satisfies~\eqref{BinGGred}. In the recursive process determining $\bar B(z)$ in Appendix~\ref{ap:Bbarsol}, it is shown that $d(u)^{(\epsp,\epsm)}_{\rho}$ factorizes in terms of inverses of ($q$-shifted) power series $\eta(z)$, see~\eqref{etaetaAp}. Inserting~\eqref{BinGGred} into  ~\eqref{currentTheta} and using the factorized formula~\eqref{etaetaAp}, it immediately follows: 
	\beqa
	{\bf \Theta}(z) &=&  \frac{(1-z^2)}{(1-q^{-2}z^2)} \frac{\rho^2}{1-q^2}\eta(z)^{-1} \eta(q^2z)^{-1} {\bar B}(z)^{-1} {\bar B}(q^2z)^{-1}\label{eq:ThetaBBar}
	 \ .
	\eeqa
\begin{prop}\label{prop:BbarH} The auxiliary generating function ${\bar B}(z)$ is such that:
	\beqa
	{\bar B}(z) = \frac{\rho}{(q-q^{-1})} \eta(z)^{-1} \exp\left(-(q-q^{-1}) \sum_{m=1}^{\infty} \frac{{\acute H}_mz^m} {1 + q^{2m}}\right) \ ,\label{BbarHm} 
	\eeqa
	where $\eta(z)^{-1}$ is given in Lemma~\ref{lem:eta-recur} with $\eta_0=1$.
\end{prop}
\begin{proof} Recall the generating function ${\bf \Theta}(z)$ written as a power series in the generators ${\acute H}_m$, see~\eqref{eq:ThetaHmbis}. 
Knowing that the leading  term of the power series ${\bar B}(z)$ is $\rho/(q-q^{-1})$, and  the leading term of $\eta(z)^{-1}$ is $\eta^{inv}_0=1$, see Example~\ref{exC1}, we can make the following Ansatz:
\beqa
{\bar B}(z) = \frac{\rho}{(q-q^{-1})} \eta(z)^{-1} \exp b(z) \quad \mbox{with}\quad  b(z)= \sum_{m=1}^{\infty} b_m\acute{H}_m z^m\ . 
\eeqa
The comparison of the r.h.s.\ of~\eqref{eq:ThetaHmbis} to the r.h.s.\ of~\eqref{eq:ThetaBBar} shows that the $b_m$'s are fixed uniquely, as displayed in~\eqref{BbarHm}.
\end{proof}

\section{The counital comodule algebra $(O_q,\delta_{O_q},\epsilon_{O_q})$: FRT type presentation}\label{sec:FRTOq}
For the alternating central extension of the $q$-Onsager algebra $\cal A_q$, an R-matrix presentation is known \cite{BSh1,BasBel}, as well as a coaction and counit maps~\cite{LBG,BGL25}. In this section, we use the map $\gamma:\cal A_q \rightarrow O_q$ given in~\eqref{mapDeldel} and Terwilliger's correspondence discussed in Section~\ref{sec:Tercor} to get the FRT type presentation for $(O_q,\delta_{O_q},\epsilon_{O_q})$, see Theorem~\ref{thm:FRTOqDr}. As a preliminary, in Section~\ref{sec:Hopf} we recall the FRT presentation for~$\Loop$.
In Section~\ref{sec:RmatOq}, we give two different R-matrix presentations for $O_q$, in the alternating generators from Section~\ref{subalt} and in terms of root vectors  from Section~\ref{sec:rootpres}. In this framework, the coaction and counit maps are given in Section~\ref{sec:delcoun}, see Propositions~\ref{prop:deltaK},~\ref{prop:Kaq-Km}.

\subsection{The Hopf algebra $(\Loop,\Delta,S,\epsilon)$: FRT presentation}\label{sec:Hopf}
The FRT presentation of the Hopf algebra $(\Loop,\Delta,S,\epsilon)$ follows from the work~\cite{DF93}, whose results are adapted to our conventions. Let  $R(u)$ be the symmetric $R$-matrix with formal parameter $u$:
	\begin{align}
		R(u) =\left(
		\begin{array}{cccc} 
			uq -  u^{-1}q^{-1}    & 0 & 0 & 0 \\
			0  &  u -  u^{-1} & q-q^{-1} & 0 \\
			0  &  q-q^{-1} & u -  u^{-1} &  0 \\
			0 & 0 & 0 & uq -  u^{-1}q^{-1}
		\end{array} \right) \ .\label{R}
	\end{align}
	It satisfies the quantum Yang-Baxter equation\footnote{As usual, $R_{nm}$
are the R-matrices on corresponding
products of spaces: 
$$R_{12} = R\otimes I\!\!I, \quad R_{23}= I\!\!I \otimes R, \quad R_{13} = (P \otimes I\!\!I)R_{23}(P \otimes I\!\!I)\ ,$$ 
with $P(a \otimes b) = b \otimes a$ and $I\!\!I$ is the two-by-two identity matrix.}
	\begin{align}
		R_{12}(u/v)R_{13}(u)R_{23}(v)=R_{23}(v)R_{13}(u)R_{12}(u/v)\ .\label{YB}
	\end{align}
	Let us introduce the following generating functions in $\Loop$-root vectors introduced in Appendix~\ref{ap:rootvecdef} (with $\lambda=q-q^{-1}$):
\beqa
{\te}^+(u)&=& \lambda u^{-1}\sum_{k=0}^\infty   q^k(-u^{-2})^{k}\te_{k\delta+\alpha_1} \ , \qquad \quad \quad \quad {\te}^-(u)=-\lambda u^{-1}\sum_{k=0}^\infty q^{-k+1} (-u^{-2})^{k}\tf_{k\delta+\alpha_0}\tk_1\ ,\label{eq:curve}\\
{\tf}^+(u)&=&\lambda u^{-1}\sum_{k=0}^\infty q^{k-1} (-u^{-2})^{k}\tk_1^{-1}\te_{k\delta+\alpha_0} \ , \qquad \ {\tf}^-(u)=-\lambda u^{-1}\sum_{k=0}^\infty q^{-k}  (-u^{-2})^{k}\tf_{k\delta+\alpha_1}\ ,\label{eq:curvf}\\
\quad  {\tk}_{1}^+(u)&=&\tk_1^{\h}\exp\left(-\lambda\sum_{k=1}^\infty \frac{(-u^{-2})^k}{q^k+q^{-k}}\te_{k\delta} \right)\  \ ,\quad \quad \  {\tk}_{2}^+(u)=\tk_1^{-\h}\exp\left(\lambda\sum_{k=1}^\infty \frac{q^{2k}(-u^{-2})^k}{q^k+q^{-k}}\te_{k\delta} \right)\ ,\label{eq:curvkp}\\
\qquad \quad {\tk}_{1}^-(u)^{-1}&=&\tk_1^{\h}\exp\left(-\lambda\sum_{k=1}^\infty \frac{(-u^{-2})^k}{q^k+q^{-k}}\tf_{k\delta} \right)\  \ ,\quad   \  {\tk}_{2}^-(u)^{-1}=\tk_1^{-\h}\exp\left(\lambda\sum_{k=1}^\infty \frac{q^{-2k}(-u^{-2})^k}{q^k+q^{-k}}\tf_{k\delta} \right)\ .\label{eq:curvk}
\eeqa
\begin{thm}[\cite{DF93}]\label{thm:FRTUq} The Hopf algebra $(\Loop,\Delta,\epsilon)$ admits an FRT presentation in the form of a Yang-Baxter algebra. The defining relations are given by:\footnote{As usual, $L_{n}$
are L-operators  on corresponding
products of spaces: $L_{1} = L\otimes I\!\!I$ and $L_{2}= I\!\!I \otimes L$.}
	\begin{align}
		&R(u/v) \mathbf{L}^\pm_1(u) \mathbf{L}^\pm_2(v) =  \mathbf{L}^\pm_2(v) \mathbf{L}^\pm_1(u) 	R(u/v) \ , \label{RLpm1} \\ 
		&R(u/v) \mathbf{L}^\pm_1(u) \mathbf{L}^\mp_2(v) =  \mathbf{L}^\mp_2(v) \mathbf{L}^\pm_1(u) R(u/v) \ ,\label{RLpm2}
	\end{align}
	with R-matrix~\eqref{R}, and the L-operators have the unique  LDU decompositions:
\beqa
	{\bf L}^+(u) &=&	\begin{pmatrix} 1 & 0 \\ \te^+(u) & 1   \end{pmatrix}
	\begin{pmatrix} \tk_1^+(u) & 0 \\ 0 & \tk_2^+(u)   \end{pmatrix}
	\begin{pmatrix} 1 & \tf^+(u) \\ 0 & 1  \  \end{pmatrix}\ ,\label{expLp}\\
	{\bf L}^-(u) &=&	\begin{pmatrix} 1 & 0 \\ \te^-(u^{-1}) & 1   \end{pmatrix}
	\begin{pmatrix} \tk_1^-(u^{-1}) & 0 \\ 0 & \tk_2^-(u^{-1})   \end{pmatrix}
	\begin{pmatrix} 1 & \tf^-(u^{-1}) \\ 0 & 1  \  \end{pmatrix}\ .\label{expLm2}
	\eeqa
	The coproduct, antipode and counit of $\Loop$ are given, respectively, by:
	\begin{align} \label{coprodLpm}
		(\Delta \otimes \id)( {\bf L}^{\pm}(u)) &= \left ( {\bf L}^{\pm}(u) \right )_{[1]} \left( {\bf L}^{\pm}(u) \right )_{[2]} \ , \\ \label{antipLpm}
		(S \otimes \id)( {\bf L}^{\pm}(u)) &= ( {\bf L}^{\pm}(u))^{-1} \  , \\ \label{epsLpm}
		({\epsilon} \otimes \id)( {\bf L }^\pm (u)) &= 1 \ .
	\end{align}
\end{thm}
\begin{rem}\label{rem:LfromR} It is known that the FRT presentations of quantum affine algebras can be derived from the universal R-matrix~\cite{FR92}. 
In the case of $\Loop$, we consider 1-component evaluation of the universal R-matrix $\mathfrak{R}$~\cite{Tolstoy1991,Boos2012} on the formal evaluation  representations $\pi_u^{\h}\colon \Loop \rightarrow \End(\mathbb{C}^2)[u^{\pm1}]$, see~\cite[Sec.\,3.5]{LBG}. 
The L-operators~\eqref{expLp},~\eqref{expLm2} are identified with\footnote{See~\cite[App.\,A.3]{LBG} for details of calculations, where the evaluation of the root vectors in our conventions is given in~\eqref{app:piroot}.}
\begin{equation} \label{defLpm}
	\mathbf{L}^+(u)= ( \id \otimes \pi_{u^{-1}}^{\h})(\mathfrak{R}) \ , \qquad \mathbf{L}^-(u)= \lbrack (\id \otimes \pi_{u^{-1}}^{\h}) (\mathfrak{R}_{21}) \rbrack^{-1} \ ,
\end{equation}
while  the R-matrix~\eqref{R} is, up to a power series in $u/v$, given by:
\beqa
{R}(u/v) \propto(\pi^{\h}_u \otimes \pi^{\h}_v) ( \mathfrak{R}) \ .
\eeqa
A proper treatment of~\eqref{defLpm} requires the use of a bi-gradation on $\Loop$, we refer for more details to~\cite[Sec.\,4.1]{LBG}.
\end{rem}

Note that the L-operators ${\bf L}^{\pm}(u)$ are invertible~\cite[Eq.\,(4.9)]{DF93}. In particular, we have:
\beqa
\lbrack \mathbf{L}^-(u^{-1}) \rbrack^{-1}=	\begin{pmatrix} 1 & -\tf^-(u) \\ 0 & 1   \end{pmatrix}
	\begin{pmatrix} \tk_1^-(u)^{-1} & 0 \\ 0 & \tk_2^-(u)^{-1}   \end{pmatrix}
	\begin{pmatrix} 1 & 0 \\ -\te^-(u) & 1  \  \end{pmatrix}\ .\label{expLm}
\eeqa
From the definition of $\te^\pm(u)$, $\tf^\pm(u)$, $\tk_1^\pm(u)$, $\tk_2^\pm(u)$, in particular it follows:
\beqa 
\lbrack{\bf L}^\pm(u^{\pm 1})\rbrack^{\pm 1}\in \Loop[[u^{-1}]]  \otimes \End(\mathbb{C}^2) \ .
\eeqa
Also, let us mention that the coproduct (\ref{coprod}), antipode (\ref{antipode}) and counit (\ref{counit}) for $\Loop$ can be derived from the power series expansion of (\ref{coprodLpm}), (\ref{antipLpm}) and (\ref{epsLpm}), respectively. 

For the analysis in Section \ref{sec:DrOq}, we give the following lemma. The proof is straightforward using the relations~\eqref{comim} and~\eqref{app:Cartroot}, thus we skip it. 

\begin{lem}\label{lem:k1k2psiphi}
	\beqa
		\tk_1^+(u)^{-1}\tk_1^+(uq^{-1})^{-1} &=&\tk_1^{-1} \exp\left(\lambda\sum_{k=1}^\infty (-u^{-2}q)^k\te_{k\delta} \right)\ ,\label{eq:fullime}\\	
		\tk_2^-(u)\tk_2^-(uq^{-1}) &=& \tk_1 \exp\left(-\lambda\sum_{k=1}^\infty (-u^{2}q)^{-k}\tf_{k\delta}  \right)\ .\label{eq:fullimf}
	\eeqa
\end{lem}

\subsection{R-matrix presentations of $O_q$}\label{sec:RmatOq} Two R-matrix presentations associated respectively with $^{\{a\}}O_q$ or $^{\{r\}}O_q$ follow from the R-matrix presentation of ${\cal A}_q$ recalled in Appendix~\ref{sec:RpresAq}. They are given either in terms of a K-operator ${\bold K}^{\{a\}}(u)$ whose entries are power series in $u^{-1}$ with coefficients given by the alternating generators, or in terms of a K-operator ${\bold K}^{\{r\}}(u)$  whose entries are power series  in $u^{-1}$ with coefficients given by the root vectors of $O_q$. In the text below, by  Remark~\ref{rem:equiv-classes} we choose 
\beqa
k_-:= \frac{\rho}{k_+(q+q^{-1})^2}\quad \mbox{with $\rho$ fixed as \ \eqref{rhofix}}  .\nonumber
\eeqa

\subsubsection{R-matrix presentation of alternating type}
The alternating presentation $^{\{a\}}O_q$  from Definition~\ref{def:altOq} admits the following  R-matrix presentation.
	\begin{prop}\label{prop:FRTalt} $^{\{a\}}O_q$ admits a presentation in the form of a reflection algebra and a $q$-determinant condition.
 The defining relations are given by\footnote{As usual, $K_{n}(u)$
denote the K-operators such that $K_{1} = K\otimes I\!\!I, K_{2}= I\!\!I \otimes K$. Also, $\text{t}_{12}$ stands for the trace over $({\mathbb C}^2)_1\otimes ({\mathbb C}^2)_2$.}
	\beqa 
		R(u/v)\,{\bold K}^{\{a\}}_1(u)\,R(uv)\,{\bold K}^{\{a\}}_2(v)\
		&=& \ {\bold K}^{\{a\}}_2(v)\,R(uv)\,{\bold K}^{\{a\}}_1(u)\,R(u/v)\ 
		\label{REOq} \ ,\\
		\normalfont{\text{tr}}_{12}\big({\cal P}^{-}_{12} {\bold K}^{\{a\}}_1(u)  R(qu^2)  {\bold K}^{\{a\}}_2(u q) \big) &=&  d(u)^{(\epsp,\epsm)}_{\rho} \label{qdetOq} \ ,
	\eeqa
with R-matrix (\ref{R}) and $d(u)^{(\epsp,\epsm)}_{\rho}$ given by~\eqref{eq:durho}, and ${\bold K}^{\{a\}}(u)$ 
is given by 
	\begin{equation}
		{\bold K}^{\{a\}}(u)=\begin{pmatrix} 
			uq W_+(u)-u^{-1}q^{-1}W_-(u) &\frac{1}{k_-(q+q^{-1})}G_+(u)+\frac{k_+(q+q^{-1})}{q-q^{-1}} \\ 
			\frac{1}{k_+(q+q^{-1})}G_-(u) +\frac{k_-(q+q^{-1})}{q-q^{-1}}& uq W_-(u) -u^{-1}q^{-1}W_+(u) 
		\end{pmatrix}  \label{K-Oqa} \ ,
	\end{equation}
with the generating functions defined in~\eqref{c1oq} and~\eqref{c2oq}.
\end{prop}
\begin{proof} Recall $\gamma$ defined by~\eqref{mapDeldel}. Applying $\gamma\otimes \id$ to (\ref{K-Aq}) we  get~\eqref{K-Oqa}, and then~\eqref{REOq} follows from~\eqref{RE}. Applying $\gamma$ to (\ref{gammaform}) with (\ref{gamma}) and (\ref{2delust}), we  get~\eqref{qdetOq}.
	\end{proof}
We call ${\bold K}^{\{a\}}(u)$ from Proposition~\ref{prop:FRTalt}  the {\it K-operator} of alternating type for $O_q$.
 Applying $\gamma \otimes \id$ to~\eqref{K-Aqexp}, we note that:
	\beqa\label{Ka-form}
	{\bold K}^{\{a\}}(u) \in O_q[[u^{-1}]]\otimes \End(\mathbb{C}^2)\ .
	\eeqa 
The K-operator (\ref{K-Oqa}) is invertible. Indeed, applying $\gamma$ to (\ref{KinvAq}) and using~\eqref{eq:durho}, it yields to:
\begin{equation}
	\left [	{\bold K}^{\{a\}}(u) \right ]^{-1} = \frac{u^{-2}-u^{2}}{d(uq^{-1})^{(\epsp,\epsm)}_{\rho}}{\bold K}^{\{a\}}(u^{-1}) 
	\ .\label{KinvOqa}
\end{equation}
As $U|_{u\rightarrow u^{-1}}= U|_{u\rightarrow uq^{-1}}$, note that	$\left [	{\bold K}^{\{a\}}(u) \right ]^{-1}\in O_q[[u^{-1}]]\otimes \End(\mathbb{C}^2)$.

\subsubsection{R-matrix presentation of root vectors type}
Recall the root vectors presentation $^{\{r\}}O_q$ in Section~\ref{sec:rootpres}.
The R-matrix presentation of $^{\{r\}}O_q$ is obtained by applying Terwilliger's correspondence explained in  Section~\ref{sec:Tercor} to the R-matrix presentation in Proposition~\ref{prop:FRTalt}.
Recall that $\rho$  is fixed by (\ref{rhofix}).  
	\begin{prop}\label{prop:FRTOqrv} $^{\{r\}}O_q$ admits a presentation in the form of a reflection algebra.
	The defining relations are given by:
	\beqa 
	R(u/v)\,{\bold K}^{\{r\}}_1(u)\,R(uv)\,{\bold K}^{\{r\}}_2(v)\
	&=& \ {\bold K}^{\{r\}}_2(v)\,R(uv)\,{\bold K}^{\{r\}}_1(u)\,R(u/v)\ 
	\label{REOqr} \ ,
	\eeqa
with R-matrix (\ref{R}). ${\bold K}^{\{r\}}(u)$ is called the {\it K-operator} of root vectors type for $O_q$ given by
	\begin{equation}
		{\bold K}^{\{r\}}(u)=\begin{pmatrix} 
			-\frac{u^{-1}}{q^2-q^{-2}}B^-(u^{-2}){\bar B}(u^{-2})  &  \qquad A^{\{r\}}_{12}(u)  \\ 
			\frac{1}{k_+(q+q^{-1})}{\bar B}(u^{-2}) & 	-\frac{u^{-1}}{q^2-q^{-2}}{\bar B}(u^{-2})B^+(u^{-2})
		\end{pmatrix}  \label{Kr-Oq} \ 
	\end{equation}
	with
	\beqa
	\qquad {A}^{\{r\}}_{12}(u)=\frac{1}{k_-(q+q^{-1})}\left(  \frac{\rho^2\eta(u^{-2})^{-1}\eta(q^{2}u^{-2})^{-1}}{(q-q^{-1})^2}  {\bar B}(q^2u^{-2})^{-1}  - u^{-2} B^{-}(u^{-2}){\bar B}(u^{-2})  B^{+}(u^{-2})\right)
	\ , \label{eq:expordK12r}
	\eeqa 
	\normalsize
	where the generating functions  $B^\pm(z)$, ${\bar B}(z)$ are defined in (\ref{rvgen1}), (\ref{Bbargen}) and $\eta(z)^{-1}$ is given in Lemma~\ref{lem:eta-recur} and $\eta_0=1$.
\end{prop}
 We notice that in this presentation there is no   $q$-determinant condition~\eqref{qdetOq}, in fact it is  implicitly encoded in Definition~\ref{def:Bbar} of $\bar B(z)$.

\begin{proof} It is sufficient to apply  Terwilliger's correspondence of Section \ref{sec:Tercor} to Proposition~\ref{prop:FRTalt}. We convert the entries of~\eqref{K-Oqa} in terms of the generating functions for the root vectors using Lemmas~\ref{lem:GB},~\ref{lem:WB}, and simplify the resulting expressions. It produces (\ref{Kr-Oq}) with:
\small
\beqa
A^{\{r\}}_{12}(u)= \frac{{\bar B}(u^{-2})}{k_-(q+q^{-1})} \left( \frac{B(u^{-2})}{q^{-2}-1} - \frac{u^{-2}\left((q-q^{-1})(B^{+}(u^{-2}))^2 + q^{-1}(q^2u^2-q^{-2}u^{-2})B^{-}(q^{-2}u^{-2})B^{+}(u^{-2})\right) }{qu^2-q^{-1}u^{-2}} \right)\ , \nonumber
\eeqa 
\normalsize
and by construction it satisfies the reflection equation~\eqref{REOqr}, which is just rewriting of~\eqref{REOq}.
Also, by construction of $O_q$ as the central reduction of $\cal A_q$ the defining relation~\eqref{BinGGred} for $\bar B(z)$ is equivalent to~\eqref{qdetOq}. 
From the definitions of the generating functions $B^\pm(z),{\bar B}(z)$ given by~\eqref{rvgen1},~\eqref{Bbargen} with~\eqref{eq:comBbarBdel}, it is clear that the entries (11), (21), (22) of the K-operator (\ref{Kr-Oq}) are ordered according to
\begin{align}\label{eq:order-Kop}
 n\delta+\alpha_0 < k\delta < m\delta+\alpha_1\quad \mbox{for $n,m\in {\mathbb N}, k\in\N^*$}.
\end{align} 
Our next task is to reorder the r.h.s.\,of the above expression $A^{\{r\}}_{12}(u)$ in the same order. This is done using the relation
\beqa
\frac{(q-q^{-1}){\bar B}(u^{-2})B^{+}(u^{-2})  + q^{-1}(q^2u^2-q^{-2}u^{-2}){\bar B}(u^{-2})B^{-}(q^{-2}u^{-2}) }{qu^2-q^{-1}u^{-2}} = B^{-}(u^{-2}){\bar B}(u^{-2})\ ,\label{order1}
\eeqa
that is equivalent to the entry (21) of the reflection equation~\eqref{REOqr}, which in fact uses only the entries (11), (21), (22) of~\eqref{Kr-Oq},  taken at the specialization $v=qu$ and multiplied from the right by ${\bar B}(q^{-2}u^{-2})^{-1}$. Moreover, the product ${\bar B}(u^{-2})B(u^{-2})$ is rewritten as
\beqa\frac{{\bar B}(u^{-2})B(u^{-2}) }{q^{-2}-1} = \underbrace{\xi^{-1}(q^{-2}-1)^{-1}}_{=\rho^2(q-q^{-1})^{-2}}\eta(u^{-2})^{-1}\eta(q^{2}u^{-2})^{-1}{\bar B}(q^2u^{-2})^{-1} \ \label{BbarBinBbarinv}
\eeqa 
using~(\ref{BinGGred}) and~\eqref{etaetaAp}.
All together, we obtain (\ref{Kr-Oq}) with (\ref{eq:expordK12r}).	
\end{proof}

The inverse of the K-operator can be also computed and takes a rather compact form with respect to a different ordering of the root vectors. 
\begin{prop}\label{prop:Kr-inv}
${\bold K}^{\{r\}}(u)$ from (\ref{Kr-Oq}) is invertible, and its inverse is given by
	\beqa
		\big[{\bold K}^{\{r\}}(u)\big]^{-1}=-\frac{(u^{2}-u^{-2})}{d(u)^{(\epsp,\epsm)}_{\rho}} 
		\begin{pmatrix} 
			\frac{u^{-1}}{q^2-q^{-2}}B^+(u^{-2}){\bar B}(q^2u^{-2})  &  \qquad 	\bar{A}^{\{r\}}_{12}(u)  \\ 
			\frac{1}{k_+(q+q^{-1})}{\bar B}(q^2u^{-2}) & 	\frac{u^{-1}}{q^2-q^{-2}}{\bar B}(q^2u^{-2})B^-(u^{-2})
		\end{pmatrix}  \label{Kr-inv-Oq} \ 
		\eeqa
		with 
		\beqa
		\qquad \bar{A}^{\{r\}}_{12}(u)&=& \frac{1}{k_-(q+q^{-1})}\left(  \frac{\rho^2\eta(u^{-2})^{-1}\eta(q^{2}u^{-2})^{-1}}{(q-q^{-1})^2} {\bar B}(u^{-2})^{-1}  - u^{-2} B^{+}(u^{-2}){\bar B}(q^2u^{-2})  B^{-}(u^{-2})\right)
		\ . \label{eq:expordK12rinv}
		\eeqa
\end{prop}
\begin{proof} Consider (\ref{KinvOqa}). The entries $({\bold K}^{\{a\}}(u^{-1}))_{ij}$, $i,j=1,2$, are expressed in terms of $W_\pm(u^{-1}),G_\pm(u^{-1})$. Observe that  $U|_{u\rightarrow u^{-1}}= U|_{u\rightarrow uq^{-1}}$, that is used  to rewrite $({\bold K}^{\{a\}}(u^{-1}))_{ij}$ in terms of $W_\pm(uq^{-1}),G_\pm(uq^{-1})$. Then, 
we apply the correspondence of Lemmas \ref{lem:GB}, \ref{lem:WB}, and get (\ref{Kr-inv-Oq}) where ${A}^{\{r\}}_{12}(uq^{-1})$ equals
	\beqa
	 \frac{1}{k_-(q+q^{-1})} \left( \frac{B(u^{-2})}{q^{-2}-1} + \frac{u^{-2}\left((q-q^{-1})(B^{+}(u^{-2}))^2 - q(q^{-2}u^2-q^2u^{-2})B^{+}(u^{-2})B^{-}(q^2u^{-2})\right) }{q^{-1}u^2-qu^{-2}} \right) {\bar B}(q^2u^{-2})\ . \nonumber
	\eeqa  
	It is clear that the entries (11), (21), (22) of the K-operator (\ref{Kr-inv-Oq}) are ordered according to
\begin{align}\label{eq:order-Kopinv}
 n\delta+\alpha_1 < k\delta < m\delta+\alpha_0\quad \mbox{for $n,m\in {\mathbb N}, k\in\N^*$}.
\end{align} 
	To  reorder the entry (12), we use
	\beqa
	\frac{(q-q^{-1})B^-(u^{-2}){\bar B}(u^{-2})  + q^{-1}(q^2u^2-q^{-2}u^{-2}) B^+(q^{-2}u^{-2}){\bar B}(u^{-2}) }{qu^2-q^{-1}u^{-2}} = {\bar B}(u^{-2})B^{+}(u^{-2})\ ,\label{eq:43}
	\eeqa
	that is equivalent to the entry (43) of the reflection equation (\ref{REOqr}) at the specialization $v=qu$, multiplied from the left by ${\bar B}(q^{-2}u^{-2})^{-1}$, and this equation involves only the entries (11), (21), (22) of~\eqref{Kr-Oq}. 
	Using~\eqref{eq:43}  together with~\eqref{order1}, we reorder the r.h.s.\,of the above expression $\bar{A}^{\{r\}}_{12}(u)$ according to~\eqref{eq:order-Kopinv}. Together with~\eqref{BinGGred}, we finally obtain (\ref{Kr-inv-Oq}) with (\ref{eq:expordK12rinv}).	
\end{proof}	
\begin{rem}\label{rem:qtwist} 
A quotient of an example of twisted $q$-Yangian~\cite{MRS02}, denoted $SY^{\textsf{tw}}_q(\mathfrak{o}_2)$, is a quantum symmetric pair coideal subalgebra of standard type~\cite[Sec.\,11.2]{Ko12}: its generators take the form~\eqref{eq:iot1}, \eqref{eq:iot2} for fixed values of $k_\pm$ and $\varepsilon_\pm=0$, see (up to an automorphism)~\cite[Eqs.\,(11.11),~(11.12)]{Ko12}. Thus, as an algebra, $SY^{\textsf{tw}}_q(\mathfrak{o}_2) \cong O_q$. In~\cite[Sec.\,5.1]{Lu23}, an R-matrix presentation for $SY^{\textsf{tw}}_q(\mathfrak{o}_2)$ is given.
It takes the form of a twisted version of our K-operator (denoted $S(u)$) satisfying a twisted reflection equation  \cite[Def.\,2.4]{Lu23} and a $q$-determinant condition without free parameters~\cite[Lem.\,4.6]{Lu23}. Our R-matrix presentation in Proposition~\ref{prop:FRTOqrv}  for $O_q$ thus gives a generalization to the non-standard case and an independent confirmation of the R-matrix presentation for $SY^{\textsf{tw}}_q(\mathfrak{o}_2)$ given in~\cite{Lu23}.
\end{rem}

\subsection{Coaction $\delta_{O_q}$, counit $\epsilon_{O_q}$ and injection $\iota_{O_q}$}\label{sec:delcoun}
Recall $\delta_{O_q}$,  $\epsilon_{O_q}$ and  $\iota_{O_q}$ for $O_q$, see Proposition~\ref{prop:coact}, and equations~\eqref{eq:coun1} and~\eqref{eq:iot1}-\eqref{eq:iot2}, respectively. Here, we derive the analog versions in the framework of the R-matrix presentations previously introduced.

\subsubsection{The coaction map $\delta_{O_q}\colon O_q \rightarrow O_q \otimes \Loop$}
In the R-matrix presentation of $\Loop$, recall that the coproduct takes the form~\eqref{coprodLpm}.
Following~\cite[Sec.\,6.2]{LBG}, a coaction for $O_q$ can be obtained as follows, recall the notation~\eqref{eq:TTT} and that both K-operators ${\bold K}^{\{a\}}(u)$, ${\bold K}^{\{r\}}(u)$ take the tensor form~\eqref{Ka-form}.

\begin{prop} \label{prop:deltaK}
In the R-matrix presentations from Propositions~\ref{prop:FRTalt} and~\ref{prop:FRTOqrv}, the coaction $\delta_{O_q}: O_q \rightarrow O_q \otimes \Loop$ is such that 
\begin{equation} \label{copK}
		\delta_{O_q} \otimes \id \ :\quad  {\bold K}(u) \quad \mapsto \quad  \lbrack  {\bf L}^{-}(u^{-1}) \rbrack^{-1}_{[2]}  
		{\bold K}(u)_{[1]}   {\bf L}^+(u) _{[2]}\ 
\end{equation}
for ${\bold K}(u) \in\{ {\bold K}^{\{a\}}(u),{\bold K}^{\{r\}}(u)\}$. 
\end{prop}
\begin{proof}
Apply $\gamma \otimes \id$ to the coaction for $\cal A_q$ given in~\cite[Eq.\,(6.23)]{LBG}.
The fact that this coaction is an algebra map and satisfies the axioms of  Definition~\ref{def:right-counit-com} was proven in~\cite[Lem\,6.10]{LBG} and below it.
\end{proof}

\subsubsection{The counit map $\epsilon_{O_q}: O_q \rightarrow {\mathbb C}$} In the R-matrix presentation of $\Loop$, the counit has the form~\eqref{epsLpm}. A counit map for $O_q$ is derived as follows.
 All one-dimensional representations of $\mathcal{A}_q$ have been classified in~\cite[Prop.\,2.7]{BGL25}. For our purpose, it is sufficient to consider  the algebra map $\cal A_q \rightarrow {\mathbb C}$ given in~\cite[Ex.\,2.8]{BGL25}:
\begin{lem}\label{mapQSP} We have the algebra map   $\epsilon_{O_q} \colon O_q \to \mathbb{C}$ given by:
		\beqa
			\epsilon_{O_q}( G_{k+1}) &=& \epsilon_{O_q}(\tilde{G}_{k+1})= 0 \ , \notag\\  
			\epsilon_{O_q}(W_{k+1})&=&     \overline{k} \epsp + \overline{k+1} \epsm\ , \label{eq:eps-0}\\
			\epsilon_{O_q}(W_{-k}) &=&     \overline{k+1} \epsp + \overline{k} \epsm\ , \qquad k \in \mathbb{N} \ ,
		\eeqa  
where we denote $\overline{k} = k \,\mathrm{mod}\; 2$.
\end{lem}
In particular, setting $k=0$ in the last two equalities we recover~\eqref{eq:coun1}. 
From~\eqref{BWrel1},~\eqref{BWrel2} and~\eqref{Bbargen}, we can deduce a similar result for the root vectors.
\begin{cor}\label{mapQSPB}
	For any $n\in {\mathbb Z}$, $m,\ell\in {\mathbb N}^*$ we have:
	\beqa
&&{\epsilon}_{O_q}(B_{n\delta + \alpha_1})=  \overline{n+1}\epsm + \overline{n}\epsp\ ,\\
	&&{\epsilon}_{O_q}(B_{m\delta})=-q^{-1}(q+q^{-1})^{-2}q^{2m}c_m\delta_m \quad \left(\mbox{and}\quad {\epsilon}_{O_q}({\bar B}_{\ell})=0\right)
	\ , \label{counBnnc} 
	\eeqa
    where we denote $\overline{n} = n \,\mathrm{mod}\; 2$, and the $\delta_m$'s and $c_m$'s are given by~\eqref{delpair}, \eqref{delimp} and~\eqref{eq:ckp1}, respectively.
\end{cor}

Then, we compute the image via $\epsilon_{O_q} \otimes \id$ of the alternating K-operator ${\bold K}^{\{a\}}(u)$ given by~\eqref{K-Oqa}.  The proof of the following proposition is done along the same line as~\cite[Prop.\,2.9]{BGL25}, thus we skip it.
\begin{prop} \label{prop:Kaq-Km}
In the R-matrix presentations, the counit $\epsilon_{O_q}:O_q \rightarrow {\mathbb C}$ is such that:
	\begin{equation}
		 \epsilon_{O_q} \otimes \id \ :\quad {\bold K}^{\{a\}}(u)\ \ \mbox{or}\ \	{\bold K}^{\{r\}}(u) \quad  \mapsto \quad  \left(\frac{q+q^{-1}}{u^2-u^{-2}}\right)K_{scal}(u) 
	\ \label{counKanc}
	\end{equation}
with
	\begin{equation}\label{KM-spin1/2}
		K_{scal}(u)=\begin{pmatrix}
			u \epsp + u^{-1} \epsm & \kp (u^2-u^{-2})/(q-q^{-1}) \\ 
			\km (u^2-u^{-2})/(q-q^{-1}) & u\epsm + u^{-1} \epsp
		\end{pmatrix}\ .
	\end{equation}
\end{prop}
Let us point out that the two-by-two matrix with scalar entries $K_{scal}(u)$ is known to be  the most general solution of the reflection equation (\ref{REOq}) with R-matrix (\ref{R}) and substitution ${\bold K}^{\{a\}}(u)\rightarrow K_{scal}(u) $~\cite{GZ94,dVG}. Applying the counit~\eqref{counKanc} to the l.h.s.\,of ~\eqref{qdetOq} we indeed get
\beqa
	\frac{(q+q^{-1})^2}{(u^2-u^{-2})(q^2u^2-q^{-2}u^{-2})}\normalfont{\text{tr}}_{12}\big({\cal P}^{-}_{12} K_{scal\ 1}(u)  R(qu^2)  K_{scal\ 2}(uq) \big) \  = d(u)^{(\epsp,\epsm)}_{\rho}\ .\nonumber
\eeqa

\subsubsection{The injective map $\iota_{O_q}\colon O_q \rightarrow \Loop$}
From Propositions~\ref{prop:deltaK} and~\ref{prop:Kaq-Km}, it immediately follows: 
\begin{cor}\label{cor:iotK}
In the R-matrix presentations, the injective map $\iota_{O_q}\colon O_q \rightarrow \Loop$  
defined in~\eqref{iota-B}, recall also Remark~\ref{rem:iotadef},
is such that
\begin{equation} \label{iotK}
		\iota_{O_q} \otimes \id \ :\quad  {\bold K}(u) \quad \mapsto \quad \left(\frac{q+q^{-1}}{u^2-u^{-2}}\right)  \lbrack  {\bf L}^{-}(u^{-1}) \rbrack^{-1}   
	K_{scal}(u) {\bf L}^+(u)\ 
\end{equation}
for ${\bold K}(u) \in\{ {\bold K}^{\{a\}}(u),{\bold K}^{\{r\}}(u)\}$.
\end{cor}

In the literature, elements of the form as displayed on the r.h.s.\,of~\eqref{iotK} are usually called {\it dressed K-operators}. The dressed K-operator's entries will be studied in details in Section~\ref{sec:dressK}.

\subsection{The FRT type presentation of $(O_q,\delta_{O_q},\epsilon_{O_q})$}
We now collect the results of previous sections. The following theorem can be viewed as the analog for $O_q$ of the FRT presentation of $\Loop$~\cite[Def.\,3.2 \& ~eqs.\,(4.8),~(4.9)]{DF93} previously recalled in Theorem~\ref{thm:FRTUq}.
Recall first Definition~\ref{def:right-counit-com} of a counital comodule algebra  and that $(O_q,\delta_{O_q},\epsilon_{O_q})$ is such an algebra by Proposition~\ref{prop:coact}.

\begin{thm} \label{thm:FRTOqDr} The counital  comodule algebra $(O_q,\delta_{O_q},\epsilon_{O_q})$ admits an FRT type presentation in the form of a reflection algebra. The defining relations are given by:
	\beqa 
	R(u/v)\,{\bold K}^{\{r\}}_1(u)\,R(uv)\,{\bold K}^{\{r\}}_2(v)
	&=& {\bold K}^{\{r\}}_2(v)\,R(uv)\,{\bold K}^{\{r\}}_1(u)\,R(u/v)\ ,
	\nonumber
	\eeqa
	where	${\bold K}^{\{r\}}(u)$ has the unique `$\sigma_1$-twisted' LDU-decomposition: 
	\beqa
	{\bold K}^{\{r\}}(u)= \sigma_1  {\cal L}(u) {\cal D}(u) {\cal U}(u) \label{LDU-Kop}
	\eeqa
	with
	\begin{equation}
		\sigma_1= \begin{pmatrix} 0 & 1 \\ 1 & 0  \end{pmatrix}\ ,\qquad  {\cal D}(u)  =    \begin{pmatrix} \frac{1}{k_+(q+q^{-1})}{\bar B}(u^{-2})   & 0 \\  0 &  \frac  {\rho^2\eta(u^{-2})^{-1}\eta(q^{2}u^{-2})^{-1}}{k_-(q+q^{-1})(q-q^{-1})^2}{\bar B}(q^2u^{-2})^{-1} \end{pmatrix}  ,\nonumber
	\end{equation}
	\begin{equation}
		{\cal L}(u)=	\begin{pmatrix} 1 & 0\\  -\frac{k_+u^{-1}}{(q-q^{-1})}B^-(u^{-2}) \qquad & 1  \end{pmatrix}\ ,\qquad  {\cal U}(u)= \begin{pmatrix} 1 & \qquad -\frac{k_+u^{-1}}{(q-q^{-1})}B^+(u^{-2})\\  0 & 1 \end{pmatrix}\ \nonumber
	\end{equation}	
	\normalsize
	and generating functions (\ref{rvgen1}), (\ref{Bbargen}).  The coaction is such that:
	\beqa
	(\delta_{O_q} \otimes \id)({\bold K}^{\{r\}}(u))= \left ( \lbrack  {\bf L}^{-}(u^{-1}) \rbrack^{-1} \right)_{[2]}  
		{\bold K}^{\{r\}}(u)_{[1]} \left(  {\bf L}^+(u) \right)_{[2]}\ ,\nonumber
	\eeqa
	where the L-operators satisfy~\eqref{RLpm1},~\eqref{RLpm2}.
	The counit is  
	\beqa
	(\epsilon_{O_q} \otimes \id)({\bold K}^{\{r\}}(u))=\left(\frac{q+q^{-1}}{u^2-u^{-2}}\right)K_{scal}(u)\  \nonumber
	\eeqa
	with~\eqref{KM-spin1/2}.
\end{thm}
\begin{proof}
It only remains to show that the decomposition (\ref{LDU-Kop}) exists, and that it is unique.  Assume ${\bold K}^{\{r\}}(u)$ admits a twisted  LDU-decomposition of the form \eqref{LDU-Kop} where ${\cal L}(u)$ and ${\cal U}(u)$ are lower  and upper triangular matrices, respectively, with $1$ on the diagonal entries, whereas  ${\cal D}(u)$ is assumed to be a diagonal matrix. A comparison with the explicit expression (\ref{Kr-Oq}) shows that such a decomposition exists, and that the entries are uniquely determined as displayed. 
	\end{proof}

\section{Actions of  $\delta_{O_q}$, $\iota_{O_q}$ on the generators}\label{sec:actgen}
By Lemmas~\ref{mapQSP} and~\ref{mapQSPB}, the image via the counit $\epsilon_{O_q}$ of the alternating generators or root vectors, respectively, is known. The goal of this section is to compute the corresponding images of the coaction $\delta_{O_q}$ and the injective map $\iota_{O_q}$ defined in~\eqref{iota-B}. The results are given in terms of generating functions, see Corollaries~\ref{cor:DelWG},~\ref{cor:DeltBinKij},~\ref{cor:WKdress},~\ref{cor:BinKij}.

\subsection{Images of K-operators}\label{sec:dressK} 
Firstly, we consider  the image via $\delta_{O_q}\otimes \id$ of the K-operator ${\bold K}(u)$  in~\eqref{copK}. Let us denote its entries by
\beqa
\{\delta_{O_q}\left({\bold K}(u)_{ij}\right)\}_{i,j=1,2}\ .
\eeqa
\begin{lem}\label{copKdressij} For ${\bold K}(u) \in\{ {\bold K}^{\{a\}}(u),{\bold K}^{\{r\}}(u)\}$, we have:
\beqa
\delta_{O_q}\left({\bold K}(u)_{11}\right) &=&  {\bold K}(u)_{12} \otimes \left( \tk_1^-(u)^{-1} \te^+(u)\tk_1^+(u) + \tf^-(u) \tk_2^-(u)^{-1} \te^-(u)\te^+(u)\tk_1^+(u)  \right)\nonumber \\
&& - \  {\bold K}(u)_{21} \otimes  \tf^-(u)\tk_2^-(u)^{-1}\tk_1^+(u) \nonumber \\
&& + \  {\bold K}(u)_{11} \otimes \left( \tk_1^-(u)^{-1}\tk_1^+(u) +  \tf^-(u)\tk_2^-(u)^{-1}\te^-(u)\tk_1^+(u) \right)\nonumber \\
&& -\   {\bold K}(u)_{22} \otimes   \tf^-(u)\tk_2^-(u)^{-1}\te^+(u) \tk_1^+(u) \nonumber ,\\
\delta_{O_q}\left({\bold K}(u)_{22}\right) &=& 
- \ {\bold K}(u)_{12} \otimes \left( \tk_2^-(u)^{-1}\te^-(u)\te^+(u)\tk_1^+(u)\tf^+(u) +\tk_2^-(u)^{-1}\te^-(u)\tk_2^+(u) \right)\nonumber \\
&&+ \  {\bold K}(u)_{21} \otimes  \tk_2^-(u)^{-1} \tk_1^+(u)\tf^+(u) \nonumber \\
&&- \ {\bold K}(u)_{11} \otimes  \tk_2^-(u)^{-1} \te^-(u)\tk_1^+(u)\tf^+(u)  \nonumber \\
&&+ \  {\bold K}(u)_{22} \otimes \left( \tk_2^-(u)^{-1}\tk_2^+(u) + \tk_2^-(u)^{-1}\te^+(u)\tk_1^+(u)\tf^+(u)  \right)\ ,\nonumber\\
\delta_{O_q}\left({\bold K}(u)_{12}\right) &=& 
{\bold K}(u)_{12} \otimes \left(  \tk_1^-(u)^{-1}\tk_2^+(u) +    \tk_1^-(u)^{-1}\te^+(u) \tk_1^+(u)\tf^+(u) \right. \nonumber\\
&& \qquad \qquad\quad  \left. + \tf^-(u)\tk_2^-(u)^{-1}\te^-(u) \tk_2^+(u)+ \tf^-(u)\tk_2^-(u)^{-1} \te^-(u)\te^+(u)\tk_1^+(u)\tf^+(u)  \right)\nonumber \\
&& -\  {\bold K}(u)_{21} \otimes  \tf^-(u)\tk_2^-(u)^{-1} \tk_1^+(u)\tf^+(u) \nonumber \\
&& + \   {\bold K}(u)_{11} \otimes \left( \tk_1^-(u)^{-1}\tk_1^+(u)\tf^+(u) + \tf^-(u)\tk_2^-(u)^{-1}\te^-(u)\tk_1^+(u)\tf^+(u) \right)\nonumber \\
&&- \  {\bold K}(u)_{22} \otimes \left(\tf^-(u)\tk_2^-(u)^{-1}\tk_2^+(u)  + \tf^-(u)\tk_2^-(u)^{-1}\te^+(u)\tk_1^+(u)\tf^+(u)  \right)\ ,\nonumber\\
\delta_{O_q}\left({\bold K}(u)_{21}\right) &=& 
-\ {\bold K}(u)_{12} \otimes \tk_2^-(u)^{-1}\te^-(u)\te^+(u)\tk_1^+(u) \nonumber \\
&&+ \  {\bold K}(u)_{21} \otimes  \tk_2^-(u)^{-1} \tk_1^+(u)\nonumber \\
&&- \ {\bold K}(u)_{11} \otimes    \tk_2^-(u)^{-1} \te^-(u)\tk_1^+(u)  \nonumber \\
&&+ \  {\bold K}(u)_{22} \otimes  \tk_2^-(u)^{-1}\te^+(u)\tk_1^+(u)  \nonumber\ .
\eeqa
\end{lem}        
\begin{proof} Consider the r.h.s.\,of \eqref{copK}.
 Expanding the matrix product with the explicit expressions for the L-operators (\ref{expLp}), (\ref{expLm}), the result follows.
\end{proof}

Secondly, let us denote the r.h.s.\,of~\eqref{iotK} by 
	\beqa 
{\bold K}_{dress}(u) =	(\iota_{O_q} \otimes \id){\bold K}^{\{r\}}(u) \ ,
\label{mapKaltKncdress}
\eeqa
which also equals $(\iota_{O_q} \otimes \id){\bold K}^{\{a\}}(u)$, and introduce the following notation
\beq\label{eq:beta}
\beta= \frac{q+q^{-1}}{q-q^{-1}}\ .
\eeq

\begin{lem}\label{entreeKdr} The entries of the dressed K-operator~\eqref{mapKaltKncdress} are given by:
	\beqa
{\bf K}_{dress}(u)_{ij} =   {\bf K}^0_{dress}(u)_{ij} + \Big(\frac{q+q^{-1}}{u^2-u^{-2}}\Big)\Box(u)_{ij}\ ,\quad i,j=1,2 \ ,\label{Kdressnc}
\eeqa
\smallskip
	where the entries $\{{\bf K}^0_{dress}(u)_{ij}\}_{i,j=1,2}$ associated with the standard case $\varepsilon_\pm=0$ are given by
	\beqa
	{\bold K}^0_{dress}(u)_{11}&=& \beta\Big(k_+\left(\tk_1^-(u)^{-1}\te^+(u) \tk_1^+(u) + \tf^-(u)\tk_2^-(u)^{-1}\te^-(u)\te^+(u)\tk_1^+(u)\right) 
	- k_-\tf^-(u) \tk_2^-(u)^{-1}\tk_1^+(u)   \Big)\ ,\nonumber\\
	{\bold K}^0_{dress}(u)_{22}&=& \beta\Big(-k_+\left(\tk_2^-(u)^{-1}\te^-(u) \tk_2^+(u) + \tk_2^-(u)^{-1}\te^-(u)\te^+(u)\tk_1^+(u)\tf^+(u) \right)
	+k_- \tk_2^-(u)^{-1}\tk_1^+(u)\tf^+(u)   \Big)\ ,\nonumber\\
	{\bold K}^0_{dress}(u)_{12}&=& \beta 
	\Big(-k_-\tf^-(u)\tk_2^-(u)^{-1} \tk_1^+(u)\tf^+(u) 
	\ \nonumber\\
	&&\qquad \qquad\qquad + 
	k_+\left(\tk_1^-(u)^{-1} \tk_2^+(u) + \tk_1^-(u)^{-1} \te^+(u) \tk_1^+(u)\tf^+(u) + \tf^-(u)\tk_2^-(u)^{-1} \te^-(u)\tk_2^+(u)
	\right.\nonumber\\
	&&\qquad \qquad\qquad \qquad \qquad\qquad \qquad \qquad\qquad + \left.
	\tf^-(u)\tk_2^-(u)^{-1}\te^-(u)\te^+(u)\tk_1^+(u)\tf^+(u) \right)  
	\Big)\ ,\nonumber\\
	{\bold K}^0_{dress}(u)_{21}&=& \beta
	\Big(k_-\tk_2^-(u)^{-1} \tk_1^+(u) - k_+\tk_2^-(u)^{-1} \te^-(u) \te^+(u) \tk_1^+(u)
	\Big)\ ,\nonumber
	\eeqa
	and the contributions $\{\Box(u)_{ij}\}_{i,j=1,2}$ for the non-standard case are such that: \smallskip
	\beqa
	\Box_{11}(u)&=&  (u\epsp + u^{-1}\epsm)\left(  \tk_1^-(u)^{-1}\tk_1^+(u) + \tf^-(u)\tk_2^-(u)^{-1}\te^-(u)\tk_1^+(u)\right)\nonumber\\
	&& - (u\epsm + u^{-1}\epsp)\left(   \tf^-(u)\tk_2^-(u)^{-1}\te^+(u)\tk_1^+(u)  \right)  \ ,\nonumber\\
	\Box_{22}(u)&=&  -(u\epsp + u^{-1}\epsm)  \tk_2^-(u)^{-1}\te^-(u)\tk_1^+(u)\tf^+(u) \nonumber\\
	&& + (u\epsm + u^{-1}\epsp)\left( \tk_2^-(u)^{-1}\tk_2^+(u) +   \tk_2^-(u)^{-1}\te^+(u)\tk_1^+(u)\tf^+(u)   \right) \ ,\nonumber\\
	\Box_{12}(u)&=&  (u\epsp + u^{-1}\epsm)\left( \tf^-(u)\tk_2^-(u)^{-1}\te^-(u)\tk_1^+(u)\tf^+(u) +   \tk_1^-(u)^{-1}\tk_1^+(u)\tf^+(u)   \right)    \ \nonumber\\
	&& -(u\epsm + u^{-1}\epsp)\left( \tf^-(u)\tk_2^-(u)^{-1}\te^+(u)\tk_1^+(u)\tf^+(u) +   \tf^-(u)\tk_2^-(u)^{-1}\tk_2^+(u)   \right) \ ,\nonumber\\
	\Box_{21}(u)&=&  -(u\epsp + u^{-1}\epsm)\left( \tk_2^-(u)^{-1}\te^-(u)\tk_1^+(u)    \right)    \ \nonumber\\
	&&+(u\epsm + u^{-1}\epsp)\left( \tk_2^-(u)^{-1}\te^+(u)\tk_1^+(u)   \right) \ .\nonumber
	\eeqa
	\end{lem}
\begin{proof}
 Inserting (\ref{expLp}), (\ref{expLm}) into the r.h.s.\,of~\eqref{mapKaltKncdress}  with~\eqref{iotK}  and expanding the matrix product, the result follows. Note that we may proceed alternatively by applying  $\epsilon_{O_q}\otimes \id$ from the left to the set of equalities in Lemma~\ref{copKdressij} for ${\bold K}(u) \in\{ {\bold K}^{\{a\}}(u),{\bold K}^{\{r\}}(u)\}$ and using~\eqref{counKanc}.
\end{proof}

 As will be shown, ${\bold K}_{dress}(u)_{21}$ is,  up to an overall factor, the image via~$\iota_{O_q}$ of the auxiliary generating function $\bar{B}(u)$ introduced in Definition~\ref{def:Bbar}, see Corollary~\ref{cor:BinKij} for more details. Some properties of the entry ${\bold K}_{dress}(u)_{21}$ are now studied.
 By Corollary~\ref{cor:iotK}, ${\bold K}_{dress}(u)$ satisfies the reflection equation~\eqref{REOq} from which we  extract in particular the commutation relation:
\begin{equation}\label{eq:K12-comm}
[{\bold K}_{dress}(u)_{21} ,{\bold K}_{dress}(v)_{21}]=0\; .
\end{equation}
 Moreover, according to (\ref{eq:curve})-(\ref{eq:curvk}), the entries ${\bold K}_{dress}(u)_{ij}$ displayed in Lemma~\ref{entreeKdr} are power series in $u^{-2}$. In particular, the explicit expression for ${\bold K}_{dress}(u)_{21}$ can be formally written as:
\beqa
{\bold K}_{dress}(u)_{21} = \sum_{k=0}^\infty u^{-2k} \kappa_{21,k} \ ,\qquad {\kappa}_{21,0} = \beta k_- \ ,
\eeqa
where the coefficients $\{{\kappa}_{21,k}\}_{k\in{\mathbb N}^*} \in \Loop$, and $\beta$ is from~\eqref{eq:beta}. From~\eqref{eq:K12-comm}, one has $[{\kappa}_{21,k}, {\kappa}_{21,\ell}]=0$. Furthermore, as we have the constant term ${\kappa}_{21,0}$ is non-zero,  ${\bold K}_{dress}(u)_{21}$ is invertible by \cite[Lem.\,4.1]{Ter21d}. The following lemma is proven analogously to Lemma~\ref{lem:invertBbar}:
\begin{lem}\label{lem:invertK21}
	The generating function  ${\bold K}_{dress}(u)_{21}$ is invertible. Its inverse is given by:
	\beqa
	{\bold K}_{dress}(u)_{21}^{-1} = \sum_{k=0}^\infty u^{-2k} \kappa^{inv}_{21,k} \ ,\qquad {\kappa}^{inv}_{21,0} = (\beta k_-)^{-1} \ ,\label{expKinv}
	\eeqa
	where the coefficients $\{{\kappa}^{inv}_{21,k}\}_{k\in{\mathbb N}^*} \in \Loop$ are determined by the recurrence relations:
	\beqa
	\sum_{k=0}^{p} {\kappa}_{21,p-k}  {\kappa}^{inv}_{21,k} =0  \ ,\qquad p \in {\mathbb N}^*\ .\label{eq:recK21}
	\eeqa
\end{lem}

\begin{lem}\label{factBbarnc}
The generating function 	${\bold K}_{dress}(u)_{21}$ can be written as a factorized expression, in two different ways:
			\beqa
			{\bold K}_{dress}(u)_{21} &=& \beta k_- \tk_2^-(u)^{-1} \tk_1^+(u) {\cal F}_1(u) \label{facK211}\\
			&=&
			\beta k_-  {\cal F}_2(u) \tk_2^-(u)^{-1} \tk_1^+(u)\ ,\label{facK212}
			\eeqa
			where 
			\beqa
			\qquad  {\cal F}_1(u)&=& 1 -  \frac{k_+}{k_-} \tk_1^+(u)^{-1} \te^-(u) \te^+(u) \tk_1^+(u) \label{eq:F1}\\
		&&
		 +\left(\frac{q+q^{-1}}{u^2-u^{-2}}\right) \frac{1}{\beta k_-}\tk_1^+(u)^{-1}\Big( -(u\epsp + u^{-1}\epsm)\te^-(u) + (u\epsm +u^{-1}\epsp) \te^+(u)   \Big)\tk_1^+(u) \nonumber
		\eeqa
		and
		\beqa
			 \qquad \ {\cal F}_2(u) &=& 
			 1 -  \frac{k_+}{k_-} \tk_2^-(u)^{-1} \te^-(u) \te^+(u) \tk_2^-(u) \label{eq:F2}\\
		&&+\left(\frac{q+q^{-1}}{u^2-u^{-2}}\right) \frac{1}{\beta k_-}\tk_2^-(u)^{-1}\Big( -(u\epsp + u^{-1}\epsm)\te^-(u) + (u\epsm +u^{-1}\epsp) \te^+(u)   \Big)\tk_2^-(u)\ . \nonumber
			\eeqa
		\end{lem}
		\begin{proof}
		We notice by a straightforward calculation that~\eqref{Kdressnc} for $i=2,j=1$ factorizes in two different ways, as displayed.
		\end{proof}
\begin{lem}\label{lem:invertF12}
	The generating functions  ${\cal F}_i(u)$, $i=1,2$, are invertible.
\end{lem}
\begin{proof} By Lemma \ref{lem:invertK21}, ${\bf K}_{dress}(u)_{21}$ is invertible. Also, $\tk^+_1(u)$ is invertible, see~\eqref{eq:curvkp}. From  (\ref{facK211}), it follows  ${\cal F}_1(u) = (\beta k_-)^{-1} \tk_1^+(u)^{-1} \tk_2^-(u)  {\bold K}_{dress}(u)_{21}$ is invertible. We use similar arguments to show ${\cal F}_2(u)$ is invertible too.
\end{proof}

\subsection{Images of the generators via $\delta_{O_q}$} The coaction of the generators associated with each presentations $^{\{a\}}O_q$ or  $^{\{r\}}O_q$ can be derived from Proposition~\ref{prop:deltaK} using generating functions. Recall the expressions of $\delta_{O_q}\big({\bold K}(u)_{ij}\big)$, $i,j=1,2$, in Lemma~\ref{copKdressij}.
\begin{cor}\label{cor:DelWG} The coaction $\delta_{O_q}$  of the alternating generating functions (\ref{c1oq}), (\ref{c2oq}) is given by:
	\beqa
	\delta_{O_q}(W_{+}(u))&=& (q^2u^2-q^{-2}u^{-2})^{-1}\Big( uq	\delta_{O_q}\left({\bold K}^{\{a\}}(u)_{11}\right) + u^{-1}q^{-1}\delta_{O_q}\left({\bold K}^{\{a\}}(u)_{22}\right)\Big)\ ,\label{DeltaW0}\\
	\delta_{O_q}(W_{-}(u))&=&(q^2u^2-q^{-2}u^{-2})^{-1}\Big( uq	\delta_{O_q}\left({\bold K}^{\{a\}}(u)_{22}\right) + u^{-1}q^{-1}\delta_{O_q}\left({\bold K}^{\{a\}}(u)_{11}\right)\Big) \ ,\label{DeltaW1}\\
	\delta_{O_q}(G_{+}(u))&=& k_-(q+q^{-1}) \delta_{O_q}\left({\bold K}^{\{a\}}(u)_{12}\right) -\frac{\rho}{q-q^{-1}} 1 \otimes 1 \ ,\\
	\delta_{O_q}(G_{-}(u))&=&  k_+(q+q^{-1}) \delta_{O_q}\left({\bold K}^{\{a\}}(u)_{21}\right) -\frac{\rho}{q-q^{-1}} 1 \otimes 1 \  ,\label{coprodWG} 
	\eeqa
\end{cor}
Using expressions in Lemma~\ref{copKdressij}, one can expand the r.h.s.\,of each of the four equations above, and then
extract an explicit form of the coaction on the alternating generators. 
\begin{example} For the first alternating generators, the coaction acts as: 
\begin{align*}
	\delta_{O_q}\left(W_0\right)& = W_0 \otimes \tk_1 + 1 \otimes {\bf w}_0 \ ,\\
	\delta_{O_q}\left(W_1\right)& = W_1 \otimes \tk_0 + 1 \otimes {\bf w}_1 \ ,\\
	\delta_{O_q}(\tilde{G}_1) &= \tilde{G}_1 \otimes 1 + 1 \otimes \tilde{{\bf g}}_1 + W_0 \otimes \lbrack \tk_1,{\bf w}_1 \rbrack_q + W_1 \otimes \lbrack {\bf w}_0, \tk_0 \rbrack_q
\end{align*}
with 
\beqa
{\bf w}_0 &=&
 k_+q\te_1 \tk_1 + k_-\tf_1 \  ,\nonumber\\
{\bf w}_1
 &=& k_-q\te_0 \tk_0 + k_+\tf_0 \  ,\nonumber\\
 \tilde{\bf g}_1 &=&   k_+^2(q^2-q^{-2})\te_1\tk_1\tf_0 + k_+k_-(q\te_\delta - q^{-1}\tf_\delta) \ . \nonumber
\eeqa
\end{example}
At the leading order of (\ref{DeltaW0}), (\ref{DeltaW1}), one recovers (\ref{eq:del1}), (\ref{eq:del2}). As a consistency check, the leading order of (\ref{coprodWG}) produces $\delta_{O_q}(\tilde{G}_1)=[\delta_{O_q}(W_0),\delta_{O_q}(W_1)]_q$ in agreement with Example \ref{ex1}. \vspace{1mm}  

The coaction of the generating functions of the root vectors is obtained  from Corollary~\ref{cor:DelWG} using expressions in Lemma~\ref{lem:BW} and Definition~\ref{def:Bbar}, and ${\bold K}^{\{r\}}(u)_{ij}$ instead of ${\bold K}^{\{a\}}(u)_{ij}$.

\begin{cor}\label{cor:DeltBinKij} The coaction $\delta_{O_q}$ of the generating functions in the root vectors (\ref{rvgen1}) and (\ref{Bbargen}) is 
	\beqa
	\delta_{O_q}\left(B^+(u^{-2})\right) &=& -\frac{(q-q^{-1})}{k_+}u \lbrack  \delta_{O_q} ( {\bold K}^{\{r\}}(u)_{21}) \rbrack^{-1}  \delta_{O_q}( {\bold K}^{\{r\}}(u)_{22})     \ ,\nonumber\\
\delta_{O_q}\left(B^-(u^{-2})\right)	&=& -\frac{(q-q^{-1})}{k_+}u  \   \delta_{O_q} ( {\bold K}^{\{r\}}(u)_{11}   )   \lbrack   \delta_{O_q} ( {\bold K}^{\{r\}}(u)_{21}) \rbrack^{-1}  \ ,\nonumber\\
\delta_{O_q}\left(B(u^{-2})\right) &=& \frac{\xi^{-1}}{k_+^2(q+q^{-1})^4} \frac{d(uq^{-1})^{(\epsp,\epsm)}_{\rho}}{(u^2-u^{-2})} \lbrack   \delta_{O_q} ( {\bold K}^{\{r\}}(u)_{21}) \rbrack^{-1} \lbrack   \delta_{O_q} ( {\bold K}^{\{r\}}(uq^{-1})_{21}) \rbrack^{-1} \ ,  \nonumber\\
\delta_{O_q}\left({\bar B}(u^{-2})\right) &=& k_+(q+q^{-1})  \delta_{O_q}( {\bold K}^{\{r\}}(u)_{21})  \ ,\nonumber
	\eeqa
	with (\ref{xi}) and the expressions in Lemma \ref{copKdressij} are used.
\end{cor}

\begin{rem} Recall the Drinfeld type presentation of $O_q$ in Proposition~\ref{prop:DrOq}. We note that the coaction on Lu-Wang's generating function ${\bold B}(z)$ for real roots is  the sum \eqref{currentB} of coactions of $B^\pm(z^{\pm 1})$ that are the first two lines in Corollary~\ref{cor:DeltBinKij};
The coaction on Lu-Wang's generating function ${\bold \Theta}(z)$ for imaginary roots follows from~\eqref{currentTheta} and the third line of Corollary~\ref{cor:DeltBinKij}.
\end{rem}

\subsection{Images of the generators via $\iota_{O_q}$}\label{sec:images}
Recall the expressions of ${\bold K}_{dress}(u)_{ij}$, $i,j=1,2$, in Lemma~\ref{entreeKdr}.
\begin{cor}\label{cor:WKdress} The image via $\iota_{O_q}$ of the alternating generating functions (\ref{c1oq}), (\ref{c2oq})  is such that:
	\beqa
	W_+(u) &\mapsto& (q^2u^2-q^{-2}u^{-2})^{-1}\left( uq	{\bold K}_{dress}(u)_{11} + u^{-1}q^{-1}	{\bold K}_{dress}(u)_{22}\right)\ ,\label{mapWp}\\
	 W_-(u) &\mapsto&  (q^2u^2-q^{-2}u^{-2})^{-1}\left( uq	{\bold K}_{dress}(u)_{22} + u^{-1}q^{-1}	{\bold K}_{dress}(u)_{11}\right)\ ,\label{mapWm}\\
      G_+(u) &\mapsto& k_-(q+q^{-1}) {\bold K}_{dress}(u)_{12} -\frac{\rho}{q-q^{-1}} \ ,\\
	 G_-(u) &\mapsto& k_+(q+q^{-1}) {\bold K}_{dress}(u)_{21} -\frac{\rho}{q-q^{-1}} \ , \label{expGm}
	\eeqa
	where  ${\bold K}_{dress}(u)_{ij}$ are given in terms of $\te^\pm(u)$, $\tf^\pm(u)$, $\tk_1^\pm(u)$, $\tk_2^\pm(u)$ in Lemma \ref{entreeKdr}.
	\end{cor}

Recall that $\tk_1^{-1}=\tk_0$. Comparing the power series in $u^{-1}$ on both sides of the above maps, the coefficients of each power series can be compared. 
\begin{example} For the standard case $\varepsilon_\pm=0$, the image via $\iota_{O_q}$ of the first alternating generators is
	\beqa
	&& W_0 \mapsto k_+ q \te_1\tk_1 + k_-\tf_1\ ,\qquad  	W_1 \mapsto k_- q \te_0\tk_0 + k_+\tf_0\ ,\nonumber\\
	&&	\tilde{G}_1 \mapsto  k_+^2(q^2-q^{-2})q^2\te_1\tf_0\tk_1 + k_+k_-(q\te_\delta - q^{-1}\tf_\delta) \ ,\nonumber\\
	&&W_{-1} \mapsto \frac{1}{q+q^{-1}} \Big ( k_- (q^2 \te_0 \tk_0 -\tf_{\delta+\alpha_1}) + k_+\left(q \tf_0 -q^3 \te_{\delta+\alpha_1}\tk_1 +(q^2-1)^2 \te_1\tf_1\tf_0\tk_1\right) \Big) \nonumber\\
&&\qquad\quad	+ \frac{q-q^{-1}}{(q+q^{-1})^2} \Big ( k_+ q^2 \te_1 (\te_\delta+\tf_\delta)k_1 + k_-(q \te_\delta \tf_1 - q^{-1} \tf_1 \tf_\delta) \Big)\ , \nonumber\\
        && W_{2} \mapsto  
        \frac{1}{q+q^{-1}} \Big ( k_- q (\tf_1-q^2 \te_{\delta+\alpha_0} \tk_1^{-1} ) +k_+ (q^2 \te_1 \tk_1 - \tf_{\delta+\alpha_0})  \Big ) \nonumber\\
        &&\qquad \quad+ \frac{q-q^{-1}}{q(q+q^{-1})^2} \Big ( k_+ \big ( (q^4-1) \te_1\te_0\tf_0 - (q^4 \te_\delta + \tf_\delta) \tf_0 \big )+ k_- q(q^2 \te_\delta \te_0- \te_0\tf_\delta)k_1^{-1} \Big )\ .\nonumber
     \eeqa
	\end{example}
Alternatively, let us consider the root vectors presentation of $O_q$ and its image in $\Loop$. 
Recall that 
${\bold K}_{dress}(u)_{21}$ is invertible by Lemma \ref{lem:invertK21}.  Using \eqref{BinGGred}, (\ref{BWrel1}), (\ref{BWrel2}) and  Corollary~\ref{cor:WKdress}, it then follows:
\begin{cor}\label{cor:BinKij} The image via $\iota_{O_q}$ of the generating functions in the root vectors (\ref{rvgen1}) and (\ref{Bbargen}) is 
	\beqa
	B^+(u^{-2}) &\mapsto& -\frac{(q-q^{-1})}{k_+}u {\bold K}_{dress}(u)_{21}^{-1}{\bold K}_{dress}(u)_{22}\ ,\label{eq:Bpdressp}\\
	B^-(u^{-2}) &\mapsto& -\frac{(q-q^{-1})}{k_+}u {\bold K}_{dress}(u)_{11}{\bold K}_{dress}(u)_{21}^{-1} \ ,\label{eq:Bpdressm}\\
	B(u^{-2}) &\mapsto& \frac{\xi^{-1}}{k_+^2(q+q^{-1})^4}  \frac{d(uq^{-1})^{(\epsp,\epsm)}_{\rho}}{(u^2-u^{-2})}{\bold K}_{dress}(u)_{21}^{-1}{\bold K}_{dress}(uq^{-1})_{21}^{-1} \qquad \mbox{with (\ref{xi})} \ ,\label{eq:Bpdressim}\\
	{\bar B}(u^{-2}) &\mapsto& k_+(q+q^{-1}){\bold K}_{dress}(u)_{21} \ ,\nonumber
	\eeqa
	where  ${\bold K}_{dress}(u)_{ij}$ are given in terms of $\te^\pm(u)$, $\tf^\pm(u)$, $\tk_1^\pm(u)$, $\tk_2^\pm(u)$ in Lemma \ref{entreeKdr}.
\end{cor}
 Comparing the power series on both sides of~\eqref{eq:Bpdressp}-\eqref{eq:Bpdressim}  using~\eqref{Kdressnc} and~\eqref{expKinv}, 
  we get  explicit expression for the image of all root vectors of $O_q$ in the loop algebra $\Loop$,
and for the leading coefficients,  as displayed in the next example.

\begin{example} For the standard case $\varepsilon_\pm=0$, the image via $\iota_{O_q}$ of the first root vectors of $O_q$ is
	\begin{align*}
	&& B_{\alpha_0} &\mapsto k_+ q \te_1\tk_1 + k_-\tf_1\ ,\\
	&&  	B_{\alpha_1} &\mapsto k_- q \te_0\tk_0 + k_+\tf_0\ ,\\
	&&	B_{\delta} &\mapsto -k_+^2(q^2-q^{-2})q\te_1\tf_0\tk_1 - k_+k_-(\te_\delta - q^{-2}\tf_\delta) \ .
	\end{align*}
For the non-standard case, it is easily checked that one recovers (\ref{eq:iot1}), (\ref{eq:iot2}). \end{example}

\section{Applications to the Drinfeld type presentation of $O_q$}\label{sec:DrOq}
In this section, firstly we recall the defining relations of  Drinfeld's second realization ${}^{Dr}\!\Loop$ as well as Beck's correspondence \cite{Beck,BCP} in order to relate the root vectors of Section \ref{sec:Hopf} to Drinfeld generators. To this end, the anti-automorphism $\Omega:\Loop \rightarrow \ {}^{Dr}\!\Loop$ is introduced. The explicit image via $\Omega\circ \iota_{O_q}$ of Lu-Wang's Drinfeld type currents ${\bf B}(z)$, ${\bf \Theta}(z)$ of $^{Dr}O_q$ in ${}^{Dr}\!\Loop$ is obtained, see Theorem~\ref{thm:imOqUqnc}. In particular, the image in  ${}^{Dr}\!\Loop$ of ${\bf \Theta}(z)$ -- that encodes the imaginary root vectors of $^{Dr}O_q$ -- factorizes, see Lemma~\ref{lem:factform1nc}. Finally, the explicit image via the twisted coproduct $\Delta^{\Omega}$ of ${}^{Dr}\!\Loop$ is also given in Proposition~\ref{prop:imDelOqUqnc}. Both results enhance and independently confirm  \cite[Thm.\,1.1,Thm.\,1.2]{P24}, the latter being corollaries of our construction.

\subsection{Drinfeld second realization ${}^{Dr}\!\Loop$}
The so-called Drinfeld's second realization of $U_q(\widehat{sl}_2)$ has been introduced and studied in \cite{Dr0,Dr88}. For our purpose, we consider its quotient denoted ${}^{Dr}\!\Loop$.
\begin{defn}[\cite{Dr88}]\label{def:UqDrloop}
	Drinfeld's second realization ${}^{Dr}\!\Loop$ has generators $\{{\tx}_k^{\pm}, 
	 \tho_\ell, 
	\tK^{\pm 1} |k\in {\mathbb Z},\ell\in {\mathbb Z}^* \}$ subject to the following relations:  
	\beqa
	&&   \tK\tK^{-1}=\tK^{-1}\tK=1 \ ,\label{dl1}\\
    &&\big[ \tho_k,\tho_\ell   \big]  = 0 \ , \quad \big[ \tho_k,\tK   \big]  = 0\ ,\label{dl2}\\
    && \big[ \tho_k, \tx^\pm_\ell \big]  = \pm \frac{1}{k}\big[ 2k\big]_q \tx^\pm_{k+\ell}\ ,\label{dl3}\\
	&& \tK \tx^\pm_k \tK^{-1}    =  q^{\pm 2} \tx^\pm_{k} \ ,\label{dl31}\\
	&& \tx^\pm_{k+1} \tx^\pm_\ell  -q^{\pm 2} \tx^\pm_{\ell} \tx^\pm_{k+1} = q^{\pm 2} \tx^\pm_{k} \tx^\pm_{\ell+1} - \tx^\pm_{\ell+1} \tx^\pm_k \ ,\label{dl4}\\     
	&& \big[ \tx^+_k, \tx^-_\ell   \big]  = \frac{ \psi_{k+\ell} - \phi_{k+\ell}}{q-q^{-1}}\ ,\label{dl5}
	\eeqa
	where the $\psi_{k}$ and $\phi_{k}$  are defined by the following equalities of formal power series in the indeterminate $z$:
	\beqa
	\Psig(z)&=&\sum_{k=0}^\infty \psi_{k} z^{k} = \tK \exp\left( (q-q^{-1}) \sum_{k=1}^{\infty} \tho_k  z^{k} \right) \ ,\label{def:psi}\\
	\Phig(z)&=& \sum_{k=0}^\infty \phi_{-k} z^{-k} = \tK^{-1} \exp\left(- (q-q^{-1}) \sum_{k=1}^{\infty} \tho_{-k}  z^{-k} \right)\ ,\label{def:phi}
	\eeqa
    and $\psi_{-k},\phi_k = 0$ for $k>0$.
\end{defn}
Let us mention that ${}^{Dr}\!\Loop$ can also be regarded as a graded algebra~\cite{Da12}. Following \cite[Sec.\,8.1]{P24}, to each Drinfeld generator is assigned the so-called {\it Drinfeld degree}  $\deg^{Dr}\colon {}^{Dr}\!\Loop \rightarrow {\mathbb Z}$ such that
\beqa
\deg^{Dr}(\tx^\pm_k) = \pm1, \quad \deg^{Dr} (\tho_\ell) = \deg^{Dr}(\tK^{\pm 1}) = 0\ ,\quad \ k\in {\mathbb Z}, \quad \ell \in {\mathbb Z^*}\ . \label{degxh}
\eeqa
The subalgebra:
\beqa
U_{+} =  \langle \ y\in  {}^{Dr}\!\Loop \, |\, \deg^{Dr}(y) >0  \ \rangle \ \quad \left( \mbox{resp.} \quad U_{-} =  \langle \ y\in  {}^{Dr}\!\Loop \, | \, \deg^{Dr}(y) <0  \ \rangle\ \right) \label{def:Upm}
\eeqa
is called the `Drinfeld positive half' (resp. `Drinfeld negative half'). We also introduce the subalgebra associated with degree zero:
\beqa
U_{0} =  \langle \ y\in  {}^{Dr}\!\Loop \, |\,  \deg^{Dr}(y) =0 \ \rangle \ .\label{def:U0}
\eeqa

\subsection{The anti-automorphism 
$\Omega$ of $\Loop$}
We first establish a precise relationship between the Drinfeld generators from Definition~\ref{def:UqDrloop} and our system of root vectors of 
$\Loop$ given in~\eqref{app:rootvect}. Then, we give a set of relations between Drinfeld's full currents $\tx^\pm(z)$ (generating functions of $\tx^\pm_k$) and the images via $\Omega$ of the half-currents $\te^\pm(u)$ and $\tf^\pm(u)$ previously introduced  in terms of the root vectors of $\Loop$.

\begin{defn}\label{def:Om}
 Let $\Omega$ be a degree preserving involutive algebra anti-automorphism  of $\Loop$ such that: 
\beqa
\Omega(\te_i)=\te_i\ , \quad \Omega(\tf_i)=\tf_i\ , \quad \Omega(\tk_i)=\tk_i^{-1}\ .\label{mapOmega0} 
\eeqa
\end{defn}
For the image of the root vectors~\eqref{app:rootvect} under $\Omega$ we use the capital letter notations:
\beqa
E_{k\delta + \alpha_i}\ ,\quad E_{\ell\delta}\ ,\quad  F_{k\delta+\alpha_i}\ ,\quad F_{\ell\delta}\ ,\quad K_i^{\pm 1}\ ,\quad k\in{\mathbb Z}\ ,\ \ell \in {\mathbb N}^*, \
i=0,1\ ,\eeqa
i.e.\ these are defined by setting
\beqa
E_i \equiv \Omega(\te_i)\ , \quad F_i \equiv\Omega(\tf_i)\ , \quad K_i^{-1} \equiv\Omega(\tk_i)\ ,\label{mapOmega} 
\eeqa
and by the recursions as in \eqref{app:rootvect} but with the opposite order of terms. This is exactly the recursion relations for the root system of~\cite[Prop.\,1.2]{BCP}.
With this root system and using the relation to the  Drinfeld generators~\cite[Lem.\,1.5]{BCP},  the correspondence is summarized in Table~\ref{tab:Omx}.\vspace{1mm}
\medskip

\begin{center}
	\begin{tabular}[H]{ccc}
		{\rm $x$ \cite{Tolstoy1991}} \quad \quad & \quad \quad {\rm $\Omega(x)$ \cite{BCP}} \quad\quad  &\quad \quad  $\Omega(x)$  \cite{Dr88}
		\\
		\hline
		$\te_{k\delta+\alpha_1}$ & $E_{k\delta+\alpha_1}$ &  $\tx_{k}^+$
		\\
		$\tf_{k\delta+\alpha_1}$ & $F_{k\delta+\alpha_1}$ &  $\tx_{-k}^-$
		\\
		$\te_{k\delta+\alpha_0}$ & $E_{k\delta+\alpha_0}$ &  $-\tK^{-1}\tx_{k+1}^-$
		\\
		$\tf_{k\delta+\alpha_0}$ & $F_{k\delta+\alpha_0}$ &  $-\tx_{-k-1}^+\tK$
			\\
		$\te_{\ell\delta}$ & $E_{\ell\delta}$ &  $\tho_{\ell}$
			\\
		$\tf_{\ell\delta}$ & $F_{\ell\delta}$ &  $\tho_{-\ell}$
			\\
		$\tk_1$ & $K_1^{-1}$ &  \ \ $\tK^{-1}$
			\\
		$\tk_0$ & $K_0^{-1}$ &  $\tK$
\end{tabular}
\captionof{table}{ Action
	of $\Omega: \Loop \rightarrow{}^{Dr}\!\Loop$, with $k\in{\mathbb N}$, $\ell\in{\mathbb N}^*$.}
\label{tab:Omx}
\end{center}
\bigskip

In what follows, we  also need the principal gradation homomorphism denoted $\phi_{v}$, $v\in {\mathbb C}^*$, see~\eqref{gradp}:
\beqa
&&\phi_{v}(E_{k\delta + \alpha_i})= v_{}^{-2k-1}E_{k\delta + \alpha_i}\ ,\quad \phi_{v}(F_{k\delta + \alpha_i})= v_{}^{2k+1}F_{k\delta + \alpha_i}\ ,\label{grad1} \\
&&\phi_{v}(E_{k\delta})= v_{}^{-2k}E_{k\delta}\ ,\qquad \ \ \qquad \phi_{v}(F_{k\delta})= v_{}^{2k}F_{k\delta}\ ,\label{grad2} \\
&&\phi_{v}(K_i^{\pm 1})= K_i^{\pm 1}\ .\label{grad3} 
\eeqa

According to Definition \ref{def:UqDrloop}, introduce the generating functions (Drinfeld `full currents'):
\beqa
\tx^\pm(z) = \sum_{k\in {{\mathbb Z}}}  \tx^\pm_k z^{k} \  .\label{xcu}
\eeqa
They can be actually expressed in terms of the image via $\Omega$ of the half-currents $\te^\pm(u),\tf^\pm(u),\tk_1^\pm(u),\tk_2^\pm(u)$. Indeed, from~(\ref{eq:curve}),(\ref{eq:curvf}) and using Table~\ref{tab:Omx}, one finds: 
\beqa
\tx^+(-u^{-2}q^{-1}) &=&  (q-q^{-1})^{-1}uq\left( \Omega(\te^+(uq)) -  \Omega(\te^-(u^{-1}q^{-1})) \right) \ , \label{eq:xpe} \\
\tx^-(-u^{-2}q^{-1}) &=& (q-q^{-1})^{-1}u^{-1}q^{-1}\left( \Omega(\tf^+(uq)) -  \Omega(\tf^-(u^{-1}q^{-1})) \right) \ .  \label{eq:xpf}
\eeqa
For the imaginary root vectors, by Lemma~\ref{lem:k1k2psiphi} and using Table~\ref{tab:Omx} we have
\beqa
\Psig(-u^{-2}q) &=& \tK \exp\left(\lambda\sum_{k=1}^\infty (-u^{-2}q)^k\tho_k \right) =     \Omega\left(\tk_1^+(u)^{-1}\tk_1^+(uq^{-1})^{-1}\right) \ ,\label{eq:Omegapsi}\\	
\Phig(-u^{2}q) &=& \tK^{-1} \exp\left(-\lambda\sum_{k=1}^\infty (-u^{2}q)^{-k}\tho_{-k} \right)=   \Omega\left(\tk_2^-(u)\tk_2^-(uq^{-1})\right)\ ,\label{eq:Omegaphi}
\eeqa
where the definitions \eqref{def:psi}, \eqref{def:phi}  are used. Moreover, we have:
\beqa
&& \qquad\Omega(\tk_1^{+}(u)^{-1}) =  \tK^\h\exp\left(\lambda\sum_{k=1}^\infty \frac{(-u^{-2})^k}{q^k+q^{-k}}\tho_{k} \right)\ , \quad  \Omega(  {\tk}_{2}^-(uq^{-1}))=\tK^{-\h}\exp\left(-\lambda\sum_{k=1}^\infty \frac{(-u^{2})^{-k}}{q^k+q^{-k}}\tho_{-k} \right)\ .\label{kqsquareroot}
\eeqa

\subsection{Images of Lu-Wang currents in ${}^{Dr}\!\Loop$}\label{sec:imagesLW}
We now consider the action of the injective anti-ho\-mo\-mor\-phism, recall $\iota_{O_q}$ defined by~\eqref{iota-B} together with Corollary~\ref{cor:iotK},
\beqa
\Omega \circ \iota_{O_q}\colon \quad {}^{Dr}O_q \ \rightarrow \ {}^{Dr}\!\Loop
\eeqa
on the `full current' in the real root vectors ${\bf B}(z)$ as defined in (\ref{currentB}), as well as its action on the `half-current' in the  imaginary root vectors  ${\bf \Theta}(z)$ as defined in (\ref{currentTheta}).

\begin{thm}\label{thm:imOqUqnc}
The images via \ $\Omega \circ\iota_{O_q}$  of the generating functions ${\bf B}(z)$, ${\bf \Theta}(z)$  are expressed in terms of the half-currents generating functions $\Omega(\te^{\pm}(u))$, $\Omega(\tf^{\pm}(u))$, $\Omega(\tk_{1}^{\pm}(u))$,  $\Omega(\tk_{2}^{\pm}(u))$ in the Drinfeld generators $\{\tx_k^{\pm}, \tho_{\ell}|k\in{\mathbb Z},\ell\in{\mathbb Z}^*\}$ of  ${}^{Dr}\!\Loop $ as follows:
\beqa
 \qquad \quad \Omega \circ\iota_{O_q}\Big({\bf B}(u^{-2})\Big) &=&  -\frac{(q-q^{-1})}{k_+}u \Big( {\bold K}_{dress}^{\Omega}(u)_{22}\lbrack{\bold K}_{dress}^{\Omega}(u)_{21}\rbrack^{-1} +  \lbrack{\bold K}_{dress}^{\Omega}(u^{-1})_{21}\rbrack^{-1}{\bold K}_{dress}^{\Omega}(u^{-1})_{11} \Big),\label{BDrfinnc}
\eeqa
\beqa
 \Omega \circ\iota_{O_q}\Big({\bf \Theta}(u^{-2})\Big) =    \frac{q\xi^{-1}}{k_+^2(q+q^{-1})^4} \frac{d(uq^{-1})^{(\epsp,\epsm)}_{\rho}}{(qu^2-q^{-1}u^{-2})} 
\bigl\lbrack{\bold K}_{dress}^{\Omega}(uq^{-1})_{21}{\bold K}_{dress}^{\Omega}(u)_{21}\bigr\rbrack^{-1} \ \label{ThetaDrfinnc} \ ,
\eeqa
with (\ref{xi}) and (\ref{eq:durho}),  and we set,  recall~\eqref{Kdressnc}:
\beqa
\left({\bold K}_{dress}^\Omega(u)\right)_{ij}\equiv \Omega \left({\bold K}_{dress}(u)\right)_{ij} \ ,
\label{eq:OmKdress}
\eeqa
where we used~\eqref{kqsquareroot} and 3rd column of Table~\ref{tab:Omx} for $\Omega$ images of modes of $\te^{\pm}(u)$ and $\tf^{\pm}(u)$.
For the inverse, we have
\beqa\label{eq:K21-Om-inv}
\bigl[\left({\bold K}^{\Omega}_{dress}(u)\right)_{21}\bigr]^{-1} = \Omega \big(\bigl[\left({\bold K}_{dress}(u)\right)_{21}\bigr]^{-1} \big)
\eeqa
and use~\eqref{expKinv}. 
\end{thm}
\begin{proof} We compute $\iota_{O_q}({\bf B}(z))$ and $\iota_{O_q}({\bf \Theta}(z))$ with~\eqref{currentB},~\eqref{currentTheta} using the expressions given in Corollary~\ref{cor:BinKij} to get the results. Finally, for the inverses of~\eqref{eq:OmKdress} we used the fact that $\Omega$ is an anti-homomorphism of algebras, and thus it preserves inverses.
\end{proof}

Following \cite[Eq.\,(3.1)]{P24}, define the new generating function of imaginary root vectors:
\beqa
\Thgsr(z) = (q-q^{-1})\left(\frac{1-q^{-2}z^2}{1-z^{2}} \right){\bf \Theta}(z)\ .\label{newTheta}
\eeqa 
\begin{lem}\label{lem:factform1nc}   The image via $\Omega \circ\iota_{O_q}$ of the new generating function $\Thgsr(z)$ in the imaginary root vectors of $^{Dr}O_q$ takes a factorized form:
	\beqa
	\Omega \circ\iota_{O_q}\left(\Thgsr(u^{-2})\right) &=& -q^{-1}(q-q^{-1})^{2} \left(1+ \frac{ \epsp^2 + \epsm^2 + \epsp\epsm(u^2+u^{-2})}{(u^2-u^{-2})^2} \right) \label{thetBrnc}\\
	&& \qquad\qquad \quad \times \quad  {\cal F}^{\Omega}_2(uq^{-1})^{-1} \Phig(-u^2q) \Psig(-u^{-2}q){\cal F}^{\Omega}_1(u)^{-1}\ ,\nonumber
	\eeqa
	where  we denote, recalling \eqref{eq:F1} and \eqref{eq:F2},
    $$
    {\cal F}_i^\Omega(u)\equiv\Omega \left({\cal F}_i(u)\right)\ ,
    $$
	and with the inverse ${\cal F}_i^\Omega(u)^{-1} = \Omega \left({\cal F}_i(u)^{-1}\right)$,  with Lemma~\ref{lem:invertF12} used.
\end{lem}
\begin{proof} We apply $\Omega\circ\iota_{O_q}$ to~\eqref{newTheta}. Then, we set $z=u^{-2}$ and replace $\Omega\circ\iota_{O_q}({\bf \Theta}(u^{-2}))$ by the expression given in Theorem~\ref{thm:imOqUqnc}, where 
using second equality of~\eqref{etaetaAp}. Finally, we  use~\eqref{facK211}, \eqref{facK212} and~\eqref{eq:Omegapsi},~\eqref{eq:Omegaphi}.
\end{proof}

 Now, recall (\ref{def:psi}), (\ref{def:phi}) and (\ref{eq:curve})-(\ref{eq:curvk}). Using (\ref{degxh}), it follows:
\begin{lem} \label{lem:degcur}
	\beqa
	&&\qquad \deg^{Dr}\left(\Psig(-u^{-2}q)\right) =0 \ , \quad \ \deg^{Dr}\left(\Phig(-u^{-2}q)\right)=0\ , \label{degpsiphi}\\ 
	&&\qquad\deg^{Dr}\left(\Omega(\tk_1^{\pm}(u))\right) =0 \ , \quad\ \  \deg^{Dr}\left(\Omega(\tk_2^\pm(u))\right)=0\ , \\
	&&\qquad \deg^{Dr}\left( \Omega(\te^{\pm}(u))\right)=+1 \ , \quad \deg^{Dr}\left( \Omega(\tf^{\pm}(u))\right)=-1 \ .\label{degefk}
	\eeqa
	\end{lem}

According to the explicit expressions for the ${\cal F}_i(u)$'s in  Lemma \ref{factBbarnc},  and using Lemmas~\ref{lem:invertF12} and~\ref{lem:degcur},  we find that  
\beqa
{\cal F}^{\Omega}_i(u)^{-1} \in U_+[[u^{-2}]] \ ,\qquad i=1,2\ .
\eeqa
It immediately follows from~\eqref{thetBrnc}:
\begin{cor} The following ``factorization formula'' holds:\label{cor:finnc}
	\beqa
	\Omega \circ\iota_{O_q}\left(\Thgsr(u^{-2})\right) &=&   -q^{-1}(q-q^{-1})^{2}   \left(1+ \frac{ \epsp^2 + \epsm^2 + \epsp\epsm(u^2+u^{-2})}{(u^2-u^{-2})^2} \right)  \Phig(-u^2q) \Psig(-u^{-2}q) \label{eq:corPr}\\
	&& \qquad \qquad \qquad \qquad\qquad \qquad \qquad\qquad \qquad\qquad \mod  U_+[[u^{-2}]]\ .\nonumber
	\eeqa  
\end{cor}

\begin{rem}\label{rem:compTom}
A variation of (\ref{eq:corPr}) has been independently proven from a different approach in~\cite[Thm.\,1.1]{P24} for the standard case $\varepsilon_\pm=0$ and~\cite[Thm.\,7.5]{P24} for the non-standard case. We now show that our result is consistent with~\cite{P24}.
The first step is to relate  the generators for $O_q$ with the convention~\eqref{rhofix} to the generators  $B_0^{Pr}, B_1^{Pr}$ from~\cite[Eq.\,(2.12)]{P24} with  $q$-Dolan-Grady relations~\cite[Eq.\,(2.7)]{P24}. They read:
\beqa
B_0^{Pr} = e_0^-  - c_0 e_0^+K_0^{-1} + s_0K_0^{-1} \ ,\qquad  B_1^{Pr} = e_1^-  - c_1 e_1^+K_1^{-1} +s_1K_1^{-1} \ ,\nonumber
\eeqa
where $\{e_i^\pm,f_i^\pm,K_i^{\pm 1}|i=0,1\}$ denote the Chevalley generators of $\Loop$ in \cite{P24}. In our case, the central element ${\mathbb K}_i \rightarrow 1$ in ~\cite[Remark 2.2]{P24}, i.e. $c_0=c_1=q^{-2}$. 
Moreover, the relation between the Drinfeld generators  and  $e_i^\pm,f_i^\pm,K_i^{\pm 1}$ \cite[Section~2.2]{P24} gives:
\beqa
e_0^+ \rightarrow E_0\ ,\quad e_0^- \rightarrow F_0\ ,\quad  e_1^+ \rightarrow E_1\ ,\quad e_1^- \rightarrow F_1\ ,\quad K_i \rightarrow K_i\ ,\ i=0,1\ , \label{eq:idP24}
\eeqa
in our conventions, see~Table~\ref{tab:Omx}. 
Furthermore, recall the principal gradation homomorphism $\phi_{v_0}$ for $v_0\in{\mathbb C}^*$ given by~\eqref{gradp}. From~\eqref{eq:iot1},~\eqref{eq:iot2} and using Definition~\ref{def:Om} we get  
\beqa
\phi_{v_0}\circ\Omega\circ \iota_{O_q}(W_0)&=& k_+q^{-1}v_0^{-1}E_1K_1^{-1} + k_-v_0F_1 + \varepsilon_+K_1^{-1} \ ,\nonumber \\ \phi_{v_0}\circ\Omega\circ \iota_{O_q}(W_1)&=& k_-q^{-1}v_0^{-1}E_0K_0^{-1} + k_+v_0F_0 + \varepsilon_-K_0^{-1} \ .\nonumber
\eeqa
Then, recall $\rho=k_+ k_-(q+q^{-1})^2$ with (\ref{rhofix}). A comparison of the $q$-Dolan-Grady relations~\eqref{qDG1},~\eqref{qDG2} with the ones~\cite[Eqs.\,(2.7)]{P24} implies $\phi_{v_0}\circ\Omega\circ\iota_{O_q}(W_0)=q^{1/2}(q-q^{-1})B_1^{Br}$ and $\phi_{v_0}\circ\Omega\circ\iota_{O_q}(W_1)= q^{1/2}(q-q^{-1})B_0^{Br}$. This determines (up to overall signs) the scalars $k_+,k_-,v_0$ and $\varepsilon_\pm$:
\beqa
k_\pm = -\mathsf{i}(q-q^{-1})\ ,\quad \varepsilon_+=q^{1/2}(q-q^{-1})s_1\ ,\quad \varepsilon_-=q^{1/2}(q-q^{-1})s_0\ , \quad v_0=\mathsf{i} q^{\h}\ .\label{val}
\eeqa
We now return to  (\ref{eq:corPr}). According to the correspondence displayed in Table~\ref{tab:Nota}, define
\beqa
 \Thgsr(z)|_{\mbox{\cite{P24}}} =\phi_{v_0}\circ\Omega \circ \iota_{O_q} \Big(-q(q-q^{-1})^{-2}\Thgsr(z)\Big) \label{def:ThetBr23}
 \eeqa
with $v_0$ in~\eqref{val}. Using~\eqref{eq:Omegapsi},~\eqref{eq:Omegaphi} and ~\eqref{grad1},~\eqref{grad2} one has:
\beqa
 \phi_{v_0}(\Psig(z)) =   \Psig(zv_0^{-2})  \quad   \mbox{and}\quad  \phi_{v_0}(\Phig(z)) =   \Phig(zv_0^{-2})	\ . \nonumber
\eeqa
From~\eqref{eq:corPr} and using~\eqref{val}, we get:
\beqa
\qquad \Thgsr(u^{-2})|_{\mbox{\cite{P24}}} =  \left(1+ \frac{ q(q-q^{-1})^2\big(s_0^2+s_1^2 + s_0s_1(u^2+u^{-2}) \big)}{(u^2-u^{-2})^2} \right)  \Phig(u^2) \Psig(u^{-2})   \quad \mbox{mod} \  U_+[[u^{-2}]]\ \label{resTh1P23}
\eeqa
which matches \cite[Thm.\,7.5 (1) \& Lem.\,6.1]{P24} for the identification $C=1$, $ z=u^{-2}$.
\end{rem}
\medskip

The image via $\Omega \circ \iota_{O_q}$ of Lu-Wang's generators $\{{\acute H}_m\}_{m\in{\mathbb N}^*}$ of $^{Dr}O_q$ can be also computed as follows. Recall the expression of $\bar B(z)$ in $\{{\acute H}_m\}_{m\in{\mathbb N}^*}$ given by~\eqref{BbarHm}.
\begin{prop}\label{prop:HmOqUq}
	The image via \ $\Omega \circ\iota_{O_q}$ of the generating function ${\bar  B}(u^{-2})$ is expressed in terms of the half-currents generating functions $\Omega(\te^{\pm}(u))$,  $\Omega(\tk_{1}^{+}(u))$,  $\Omega(\tk_{2}^{-}(u))$ in the Drinfeld generators $\{\tx_k^{+}, \tho_{\ell}|k\in{\mathbb Z},\ell\in{\mathbb Z}^*\}$ of  ${}^{Dr}\!\Loop $ as follows:
	\beqa
	\Omega \circ \iota_{O_q}({\bar B}(u^{-2})) &=& k_+(q+q^{-1})({\bold K}^\Omega_{dress}(u))_{21} \ 
	 ,\label{HDrfin} 
	\eeqa
	where we use~\eqref{facK211}  together with the notation~\eqref{eq:OmKdress}, 
    $\Omega$ images in~\eqref{kqsquareroot} and 3rd column of Table~\ref{tab:Omx} for the images of modes of $\te^{\pm}(u)$.
\end{prop}
\begin{proof}
We apply $\Omega$  to the last equality in  Corollary~\ref{cor:BinKij}.
\end{proof}
Inserting (\ref{BbarHm}) with $z\rightarrow u^{-2}$ into the l.h.s.\ of  (\ref{HDrfin}), the image in ${}^{Dr}\!\Loop$ of Lu-Wang's generators $\acute{H}_m$ can be extracted from the comparison between the power series in $u^{-2}$ on both sides of (\ref{HDrfin}). For instance, for the standard case $\varepsilon_\pm=0$ we find:
\beqa
\Omega \circ \iota_{O_q}(\acute{H}_1) =  q^{-1}\tho_{-1} - q\tho_1 + \frac{k_+}{k_-}q^{-2}(q^2-q^{-2})\tx_0^+\tx_{-1}^+\ .\label{eq:Hcute1}
\eeqa 
\begin{rem} For the identification~\eqref{val} and using~\eqref{dl4},~\eqref{grad1},~\eqref{grad2}, one finds that $\phi_{v_0}\circ \Omega \circ (\acute{H}_1)$  matches with~\cite[Lem.\,4.4 for $C=1$, $c_1=q^{-2}$]{P24}.
\end{rem}

\subsection{Coproduct of images of Lu-Wang's currents}
 In~\cite{P24}, a partial information on a coproduct formula of the image in ${}^{Dr}\!\Loop$ of Lu-Wang's current associated with the imaginary root vertors has been obtained. Here, we give its complete image in ${}^{Dr}\!\Loop \otimes {}^{Dr}\!\Loop $, i.e.\ the coproduct of ${}^{Dr}\!\Loop$ on $\Omega \circ \iota_{O_q}({\bf \Theta}(z))$, together with the coproduct on real root vectors from $\Omega \circ \iota_{O_q}({\bf B}(z))$. 

Recall our choice of coproduct
$\Delta$ for $\Loop$ associated with the generators $\{\te_i,\tf_i,\tk_i^{\pm 1}\}$ given by (\ref{coprod}). 
The $\Omega$-twisted coproduct 
$$\Delta^{\Omega}\ : \ {}^{Dr}\!\Loop \rightarrow {}^{Dr}\!\Loop \ \otimes \ {}^{Dr}\!\Loop $$ 
is given by 
\beqa
\Delta^{\Omega}=  (\Omega \otimes \Omega) \circ \Delta \circ \Omega^{-1}\ .\label{eq:DeltaOm}
\eeqa
From (\ref{coprod}) and using (\ref{mapOmega}) we get:
\beqa\label{eq:Delta-Om-gen}
\Delta^{\Omega}(E_i) = E_i \otimes 1 + K_i\otimes E_i\ ,\quad \Delta^{\Omega}(F_i) = F_i \otimes K_i^{-1} + 1 \otimes F_i\ ,\quad  \Delta^{\Omega}(K_i^{\pm 1}) = K_i^{\pm 1}\otimes K_i^{\pm 1}\ ,
\eeqa
which coincides with the coproduct given in~\cite[Prop.\,1.1]{BCP}. We notice that while $\Omega$ is the algebra anti-isomorphism from $\Loop$ to ${}^{Dr}\!\Loop$, it is a coalgebra isomorphism between $\Loop$ with $\Delta$ and ${}^{Dr}\!\Loop$ with the coproduct~\eqref{eq:Delta-Om-gen} and the same counit~$\epsilon$, noticing that $\Omega$ preserves $\epsilon$.

\begin{prop}\label{prop:imDelOqUqnc}   The  coproduct $\Delta^\Omega$ on the generating functions $\Omega \circ\iota_{O_q}\left({\bf B}(z)\right)$ and $\Omega \circ \iota_{O_q}({\bf \Theta}(z))$ is expressed in terms of the tensor product of half-currents generating functions  $\Omega(\te^{\pm}(u))$, $\Omega(\tf^{\pm}(u))$, $\Omega(\tk_{1}^{\pm}(u))$,  $\Omega(\tk_{2}^{\pm}(u))$ in the Drinfeld generators $\{\tx_k^{\pm}, \tho_{\ell}|k\in{\mathbb Z},\ell\in{\mathbb Z}^*\}$ of  ${}^{Dr}\!\Loop $ as follows:
    \beqa
\Delta^\Omega\left(\Omega \circ\iota_{O_q}\left({\bf B}(u^{-2})\right)\right)
=  -\frac{(q-q^{-1})}{k_+}u \Delta^\Omega\Big( {\bold K}_{dress}^{\Omega}(u)_{22}\lbrack{\bold K}_{dress}^{\Omega}(u)_{21}\rbrack^{-1} +  \lbrack{\bold K}_{dress}^{\Omega}(u^{-1})_{21}\rbrack^{-1}{\bold K}_{dress}^{\Omega}(u^{-1})_{11} \Big) \ ,\label{copBDrfinnc}\nonumber
\eeqa
	\beqa
&&\Delta^\Omega\left(\Omega \circ\iota_{O_q}\left({\bf \Theta}(u^{-2})\right)\right) =  \frac{q\xi^{-1}}{k_+^2(q+q^{-1})^4} \frac{d(uq^{-1})^{(\epsp,\epsm)}_{\rho}}{(qu^2-q^{-1}u^{-2})}
\Bigl\lbrack\Delta^\Omega\bigl({\bold K}_{dress}^{\Omega}(uq^{-1})_{21}{\bold K}_{dress}^{\Omega}(u)_{21}\bigr)\Bigr\rbrack^{-1}  
	\label{copThetafinnc} 
	\eeqa
	with~\eqref{xi}, \eqref{eq:durho}, and
	\beqa
\Delta^\Omega\bigl({\bold K}_{dress}^{\Omega}(u)_{ij}\bigr) = (\Omega \otimes \Omega)\circ (\iota_{O_q}\otimes \id)\circ \delta_{O_q}\bigl({\bold K}^{\{r\}}(u)_{ij}\bigr)\ ,\label{eq:DelOmKOm}
	\eeqa
	 with the explicit expressions for $\delta_{O_q}\bigl({\bold K}^{\{r\}}(u)_{ij}\bigr)$ given in Lemma~\ref{copKdressij},  and using Theorem~\ref{thm:imOqUqnc} for the image under $\Omega\circ \iota_{O_q}$.
\end{prop}
\begin{proof} 
We apply $\Delta^{\Omega}$ to~\eqref{BDrfinnc},~\eqref{ThetaDrfinnc} to get the first two equations.
In particular, \eqref{copThetafinnc}  is determined by expression of $\Delta^\Omega\bigl({\bold K}_{dress}^{\Omega}(u)_{21}\bigr)$ due to the algebra anti-isomorphism property of $\Delta^\Omega$.
Then, to prove~\eqref{eq:DelOmKOm} we recall~\eqref{eq:OmKdress} and~\eqref{mapKaltKncdress}, and rewrite its l.h.s.\,as 
$$
\Delta^\Omega\bigl({\bold K}_{dress}^{\Omega}(u)_{ij}\bigr) = 
 (\Omega \otimes \Omega) \circ \Delta \circ \iota_{O_q} \bigl({\bold K}^{\{r\}}(u)_{ij}\bigr) \stackrel{\eqref{iota-B}}{=} (\Omega \otimes \Omega) \circ  (\epsilon_{O_q} \otimes \id ) \circ (\id_{O_q}\otimes \Delta) \circ \delta_{O_q} \bigl({\bold K}^{\{r\}}(u)_{ij}\bigr)\ .
 $$
Finally, using the coassociativity~\eqref{coassoc} of $\delta_{O_q}$ and one more time the definition~\eqref{iota-B}, we get~\eqref{eq:DelOmKOm}.
\end{proof}

\begin{rem}
In practice, we calculate the r.h.s.\,of~\eqref{copThetafinnc} as follows. First,  we have
${{\bold K}^{\{r\}}(u)}^{-1}_{21}\propto \bar B(u^{-2})^{-1}$ from~\eqref{Kr-Oq}, whose coefficients are  recursively expressed in terms of $\bar B_{n\delta}$'s, see Lemma~\ref{lem:invertBbar}. Then, using Lemma~\ref{copKdressij} and Corollary~\ref{cor:DeltBinKij}, we obtain a power series for $\delta_{O_q}\bigl({{\bold K}^{\{r\}}(u)}^{-1}_{21}\bigr)$. Finally, we take its image  by $(\Omega \otimes \Omega)\circ (\iota_{O_q}\otimes \id)$, using Corollary~\ref{cor:BinKij} together with  Table~\ref{tab:Omx}.
By~\eqref{eq:DelOmKOm}, this gives a power series for 
  $\Delta^\Omega\bigl({{\bold K}_{dress}^{\Omega}(u)}^{-1}_{21}\bigr)=\big[\Delta^\Omega\bigl({\bold K}_{dress}^{\Omega}(u)_{21}\bigr)\big]^{-1}$.

\end{rem}

 Now, consider the entries of the dressed K-operator in Lemma~\ref{entreeKdr}. By Lemma~\ref{lem:degcur} we get:
\begin{lem}\label{lem:degOmK}
	\beqa
\Omega(u{\bf K}_{dress}(u)_{11}) &\in& (U_+ + U_-)[[u^{-2}]]  \ ,\qquad \qquad \Omega({\bf K}_{dress}(u)_{12}) \in (U_- + U_0) [[u^{-2}]]
 \ ,\nonumber\\
 \Omega({\bf K}_{dress}(u)_{21}) &\in&  (U_+ + U_0) [[u^{-2}]] \ ,\qquad\  \qquad
 \Omega(u{\bf K}_{dress}(u)_{22}) \in (U_+ + U_-)[[u^{-2}]]\ .\nonumber
\eeqa
\end{lem}

Finally, we consider $\delta_{O_q}\left({\bold K}(u)_{21}\right)$  from Lemma \ref{copKdressij}. It contains four different terms, and with respect to each of them, the Drinfeld degree of its image via $\Omega \otimes \Omega$ is such that: 
\begin{lem} \label{lem:degOmDelK}
	\beqa
			 \Omega({\bf K}_{dress}(u)_{12})  \otimes   \Omega\left( -(\tk_2^-(u))^{-1}\te^-(u)\te^+(u)\tk_1^+(u) \right) \quad &\in&  \quad (U_- + U_0) \otimes U_+ [[u^{-2}]]\ ,\nonumber\\
		   \Omega({\bf K}_{dress}(u)_{21})  \otimes   \Omega\left( (\tk_2^-(u))^{-1} \tk_1^+(u)\right) \quad &\in&  \quad (U_+ + U_0 )\otimes U_0 [[u^{-2}]]\ ,\nonumber \\
	     \Omega({\bf K}_{dress}(u)_{11})  \otimes \Omega\left(
	     - (\tk_2^-(u))^{-1} \te^-(u)\tk_1^+(u)   \right) \quad &\in&  \quad (U_+ + U_-)  \otimes U_+ [[u^{-2}]]\ ,\nonumber \\
	     \Omega({\bf K}_{dress}(u)_{22})  \otimes \Omega\left(  (\tk_2^-(u))^{-1}\te^+(u)\tk_1^+(u)  \right)\quad &\in&  \quad(U_+ + U_-)  \otimes U_+ [[u^{-2}]]\ .\nonumber
	\eeqa
	\end{lem}

Using~\eqref{newTheta},~\eqref{etaetaAp}, and Lemma~\ref{lem:degOmDelK} it immediately follows:
\begin{cor}\label{cor:DelTh} The generating function $\Omega \circ \iota_{O_q}(\Thgsr(z))$ is ``approximately group-like'' in the following sense: 
	\beqa
	\Delta^\Omega	\left(\Omega \circ\iota_{O_q}\left(\Thgsr(u^{-2})\right)\right) &=& -q^{-1}(q-q^{-1})^{2}   \left(1+ \frac{ \epsp^2 + \epsm^2 + \epsp\epsm(u^2+u^{-2})}{(u^2-u^{-2})^2} \right)   \nonumber\\
	&& \!\!\!\!\!\!\!\!  \times \ \Phig(-u^2q) \Psig(-u^{-2}q)  \otimes \Phig(-u^2q) \Psig(-u^{-2}q) \quad \mod  \left( U_+\otimes U_+ + U_-\otimes U_+\right)[[u^{-2}]]\ .\nonumber\label{Deltmodnc} 
	\eeqa
\end{cor}
\begin{rem} As before, we can compare Corollary~\ref{cor:DelTh}  with the results in~\cite{P24}. Using (\ref{def:ThetBr23}) with~\eqref{val}, we recover
 \cite[Thm.\,7.5 (2)]{P24} for $\Loop$.
\end{rem}

\section{Refined $q$-characters for $\Loop$} \label{sec:refqchar}
Consider the monoidal  category $\mathrm{Rep}\Loop$, whose objects are finite-dimensional $\Loop$-modules  of type~$1$, and morphisms are homomorphisms of $\Loop$-modules. Type 1 means that the action of $\tk_1$ is diagonalizable with eigenvalues of the form $q^{\lambda}$, for a weight $\lambda\in \mathbb{Z}$. 
 Let $[\mathrm{Rep}\Loop]$ be the Grothendieck ring of the category $\mathrm{Rep}\Loop$: it has a $\mathbb{Z}$-basis  $\{[V]\}_{V\in\mathrm{Irr}}$  indexed by irreducible finite-dimensional $\Loop$-modules, and with the relation $[W] = [V] + [U]$ for any short exact sequence \ $0 \rightarrow V \rightarrow W \rightarrow U \rightarrow 0$, and with multiplication $[V]\cdot[W]:=[V\otimes W]$, whose unit is the trivial one-dimensional module.
 In~\cite{FR99,FMu99}, $q$-characters of $\Loop$-modules are introduced and studied
  in order to describe relations  in this basis of $[\mathrm{Rep}\Loop]$.
  
In Section~\ref{sec:refined-q-char}, we introduce a new type of $q$-characters for $\Loop$ -- called \textit{refined} $q$-characters -- using the FRT presentation of $O_q$  given by Theorem~\ref{thm:FRTOqDr}, see~Definition~\ref{def:chicute} and Theorem~\ref{thm:homhcute}.
The corresponding $q$-character map is a ring homomorphism:
 \beqa
\acute{\chi}_q\ : \ [\mathrm{Rep}\Loop] \ \rightarrow \ O_q\tilde{H}[[z]]\ ,\label{eq:typechi0}
\eeqa
where  $O_q\widetilde{H}\subset O_q$ is the commutative subalgebra
 \begin{equation}\label{eq:OqHt}
         O_q\tilde{H}= \langle \acute{H}_m|m\in{\mathbb N}^*\rangle 
         \end{equation}
generated by  Lu-Wang's generators $\acute{H}_m$ \cite{LW20} from Definition~\ref{def:Hacute} and which are polynomials in the imaginary root vectors $B_{n\delta}$ via~\eqref{currentTheta} and~\eqref{rvgen1}. 
For small spin-$j$ evaluation representations, the refined $q$-characters are computed explicitly in Proposition~\ref{prop:acutechi-half}, Example~\ref{valchi1} and Appendix~\ref{app:refchar}. 

Furthermore, in Section~\ref{sec:defchi} we discuss how  the refined $q$-characters are related to a conjectured universal K-matrix for $O_q$, culminating in Conjecture~\ref{conj:chiunivK}. Finally, below this conjecture we give a representation-theoretic meaning of the refined $q$-characters, we describe an algorithm of extracting relations in $[\mathrm{Rep}\Loop]$. 

\medskip

Similarly to the {\it standard} $q$-characters of $\Loop$~\cite{FR99}  reviewed below, the refined $q$-character map sends an arbitrary finite-dimensional $\Loop$-module  to a polynomial in some generating functions $\acute{Y}_a^{\pm 1}$, $a\in{\mathbb C}^*$, introduced in~\eqref{eq:Ya}. 
To prepare the logic we follow for the construction of the refined $q$-characters, as a preliminary we briefly recall the construction of the standard $q$-characters using the FRT framework for $\Loop$~\cite{FR99}.

\subsection{Standard $q$-characters for $\Loop$ from L-operators}  Given the universal R–matrix $\mathfrak{R}$ of a quantum affine algebra, Frenkel-Reshetikhin constructed an injective homomorphism of commutative rings~\cite{FR99}, see also~\cite[Thm.\,2.2\,\&\,eq.\,(2.3)]{FMu99}:
		\beqa
		&& \chi_q\colon \ [\mathrm{Rep}\Loop] \ \rightarrow \ {\cal Y} = {\mathbb Z}[Y_{a}^{\pm 1}]_{a\in{\mathbb C}^*} \label{eq:char}
		\eeqa
		 with 
		 \beqa
		  Y_a= \tK^{-\h}\exp\Big(-(q-q^{-1})\sum_{k=1}^{\infty} \frac{z^ka^k}{q^k+q^{-k}}\tho_{-k}\Big)\ \ \in \  U_q\tilde{h}[[z]]\ ,\label{eq:chiY}
		\eeqa
		 where 
         \begin{equation}\label{eq:Uqth-def}
         U_q\tilde{h}= \langle \tK^{\pm \h},\tho_{-\ell}|\ell\in{\mathbb N}^*\rangle \ ,
         \end{equation}
         recall the Drinfeld generators in Definition~\ref{def:UqDrloop}.  The  image  $\chi_q(V)$ is called the $q$–character of the finite-dimensional $\Loop$-module $V$ of $\Loop$.  
         
         In the construction of the standard $q$-characters, there are two main ingredients. The first one is a ring  homomorphism called \textit{the universal transfer matrix map}:
    \begin{equation}\label{eq:t-defR}
         t : \; [\mathrm{Rep}\Loop] \to U_qb_-[[z]]\ ,
         \end{equation}
          where 
          $$
          U_qb_-=\langle \tx_{-k}^-,\tx_{-k-1}^+,\tho_{-\ell},\tK^{\pm \h}|k\in{\mathbb N},\ell\in{\mathbb N}^*\rangle 
          $$
          is the (negative) Borel subalgebra of $\Loop$, 
         see Table~\ref{tab:Omx} where the correspondence with root vectors is given. The second ingredient is an analog of the Harish-Chandra  map
$h_q: U_qb_- \to U_q \tilde{h}$, see~\cite[Def.\,2]{FR99}. Below, we review these ingredients in detail.

         The image of $t$ on $V$ is denoted by $t_V(z)$, and adapting~\cite[Eq.\,(2.1)]{FMu99} to our conventions,\footnote{The difference with~\cite{FMu99} is that the evaluation map  is in the principal (not homogeneous) gradation, and we introduce the $\Omega$ term because  the Khoroshkin-Tolstoy universal R-matrix~\cite[Eq.\,(4.23)]{Boos2012} that we use is given in terms of root vectors, and finally, we take the inverse due to the fact that our $\tk_i$ is the inverse of $k_i$ of~\cite{FR99}.} it is given by
         \beqa\label{eq:tV-def}
         t_V(z) =  \normalfont{\text{tr}}_V\Big( \big[\big(  \Omega \otimes (\pi^V \circ \phi_{z^{-1}}) \big){\mathfrak R}_{21}\big]^{-1}\Big)
         \ ,
         \eeqa
         where $(V,\pi^V)$ is the finite-dimensional representation of $\Loop$, and $\Omega$ is given in Table~\ref{tab:Omx}, and the evaluation map $\phi_z:\Loop\rightarrow \Loop \otimes {\mathbb C}[z,z^{-1}]$ is in the principal gradation~\eqref{gradp}. 
         
         The image of $t$ is in fact a commutative subalgebra of $U_qb_-$, see details in~\cite[Prop.\,2]{FR99}. The map $t$ has the following properties~\cite[Lem.\,2]{FR99}:   
 \begin{enumerate}
\item commutativity: $[t_V(z),t_W(w)]=0$\ ,
\item additivity: given a short exact sequence \ $0 \rightarrow V \rightarrow W \rightarrow U \rightarrow 0$, $t_{W}(z) = t_V(z)+t_U(z)$\ ,
\item multiplicativity: $t_{V\otimes W}(z) = t_V(z)t_W(z)$\ .
 \end{enumerate}
 The proof that $t$ is a ring homomorphism uses the Yang-Baxter equation -- see~\eqref{RLpm1} with~\eqref{defLpm} -- which implies~(1), whereas (2)-(3) follows from  the properties of the trace and the coproduct formula $(\Delta \otimes \id)(\mathfrak R)={\mathfrak R}_{12}{\mathfrak R}_{13}$ for the universal R-matrix, respectively. In addition, for the trivial representation $\mathbb C$ of $\Loop$ we have $t_{{\mathbb C}}(z)=1 \in U_qb_-$.
 
The action of the map $t$ on $V$ is now discussed in more details. 
Any finite-dimensional $\Loop$-module $V$ 
admits a finite Jordan-Hölder (composition) series whose factors are irreducible $\Loop$-modules $V_p$, $p\in{\mathbb N}^*$, with multiplicities $m_p$. By the additive property (2) 
it thus follows 
\beqa
t_V(z) = \sum_p m_p t_{V_p}(z)\ .\label{eq:tVgen}
\eeqa
Among the examples of irreducible modules, the so-called evaluation modules $V_v^{(j)}$, $j\in\h{\mathbb N}^*$, $v\in{\mathbb C}^*$,  are of special importance \cite[Sec.\,4]{CP91}, see Appendix~\ref{sec:formevrep} for their definition.
As shown in \cite[Thm.\,4.11]{CP91},
$V_p$ 
is isomorphic to a tensor product of
 evaluation modules $V_{v_i}^{(j_i)}$ of dimension $2j_i+1$, $j_i\in \h{\mathbb N}^*$, for some 
$v_i\in{\mathbb C}^*$ and integer $N$:
\beqa
V_p \cong  V_{v_1}^{(j_1)}\otimes V_{v_2}^{(j_2)} \otimes \cdots \otimes V_{v_N}^{(j_N)}\ .\label{eq:Vptens}
\eeqa
By the multiplicative property (3), we thus have $t_{V_p}(z)= \prod_{i=1}^Nt_{V_{v_i}^{(j_i)}}(z)$. Furthermore, each $t_{V_{v}^{(j)}}(z)$ can be 
uniquely determined recursively from $t_{V_{v}^{(\h)}}(z)$ using the property (2) specialized to the short exact sequence~\cite[Prop.\,4.9\,(b)]{CP91}:\footnote{We note that our coproduct and the evaluation map  ${\mathsf{ev}}_{v}\colon \Loop \rightarrow U_qsl_2$ with~\eqref{app:evek}, associated with the principal gradation,  differ from those in~\cite{CP91}. This is why we recalculate the associated Drinfeld polynomials, and the distinction in the coproduct requires changing $q$ to $q^{-1}$ in~\cite[Prop.\,4.9]{CP91}. We also refer to~\cite[Lem.\,3.1]{LBG} for explicit construction of the embedding map in this exact sequence.}
\beqa
0 \longrightarrow  V_{v}^{(j)} \longrightarrow V_{vq^{-j+\h}}^{(\h)}  \otimes V_{vq^{\h}}^{(j-\h)} \longrightarrow   V_{vq}^{(j-1)} \longrightarrow 0\ .\label{eq:short} 
\eeqa 
Indeed, according to the Drinfeld polynomial~\eqref{app:Drpoly} computed in the principal gradation, see~\eqref{app:DrPol1}, the evaluation representation $V_{v}^{(j)}$ is characterized by the $q$-string $\{xq^{2j-1},xq^{2j-3},\cdots,xq^{-2j+1}\}$ with $x\equiv -q^{-1}v^{-2}$, 
recall~\cite[Cor.\,4.2 \& Def.\,4.6]{CP91}. 
Therefore, from~\eqref{eq:short} and by property (2)  we get the so-called TT-relations:\footnote{Noticing from~\eqref{eq:tV-def} that $t^{(j)}_v(u)=t^{(j)}_1(vu)$.}
\beqa
		  \label{TTnorm-relj}
		  	t_{v}^{(j)}(u) =   t_{vq^{-j+\h}}^{(\h)}(u) t_{vq^\h}^{(j-\h)}(u)- t_{vq}^{(j-1)}(u) \quad \mbox{with}\quad t_{v}^{(0)}(u)=1\quad \mbox{ for all $j \in  \h{\mathbb N}^*$}\ ,
		  \eeqa
where we introduced the shorthand notation 
\beqa
t_v^{(j)}(u)  &\equiv& t_{V_v^{(j)}}(u^{-1})\ .
\eeqa
Altogether, this shows that  $t_{v}^{(\h)}(u)$
determines the map $t$ uniquely.

Thus, let us determine $t_{v}^{(\h)}(u)$ explicitly. Using~\eqref{defLpm}, the spin-$\h$ universal transfer matrix
takes the form:
 \beqa
		t_v^{(\h)}(u)  &=& \normalfont{\text{tr}}_{{\mathbb C}^2}\big((\Omega \otimes \id)\mathbf{L}^{-}(u^{-1}v^{-1})\big)\ ,
		 \eeqa
         where $\mathbf{L}^{-}(u)$ is defined in~\eqref{expLm2} and we used $V_v^{(\h)} = \mathbb{C}^2$ as a vector space. We finally get that 
		  \beqa
	 t_v^{(\h)}(u) &=&  \Omega\Big(\tk_1^-(uv) + \tk_2^-(uv) + \te^-(uv)\tk_1^-(uv)\tf^-(uv)\Big)\ \in U_qb_-[[u^{-2}]]\label{eq:thalfu}\label{eq:tvhalf}
		 \ .
		 \eeqa
In particular, the Drinfeld degree can be computed. Using Table~\ref{tab:Omx} and~\eqref{degxh}, we get $\deg^{Dr}(t_v^{(\h)}(u))=0$. By the TT-relations~\eqref{TTnorm-relj} and~\eqref{eq:tVgen}, it follows:
\beqa\label{eq:deg-t}
\deg^{Dr}(t_V(z))=0 \ .
\eeqa

The  Harish-Chandra map $h_q$ is a ``projection''  in the following sense\footnote{Note that the argument presented here slightly differs from~\cite[Sec.\,3.2]{FR99}. It avoids the introduction of $U_q\tilde g=\langle \tx_{-k}^\pm,\tho_{-\ell},\tK^{\pm \h}|k\in{\mathbb N},\ell\in{\mathbb N}^*\rangle $, and doesn't require to introduce augmentation ideals of $U_q\tilde n_\pm=\langle \tx_{-k}^\pm|k\in{\mathbb N}\rangle $. In the next section, our argument extends naturally to the case of $O_q$, which explains our choice of presentation.}. According to the PBW basis of $U_qb_-$ associated with the ordering $\tx_{-k}^- < \tho_{-\ell}  < \tK^{\pm \h} < \tx_{-m-1}^+$, see \cite[Sec.\,3]{Da14},
we can write:
 \beqa\label{eq:hq-def}
\quad U_q b_- &=& U_q\tilde h \oplus \big(  U_-^{(1)}\cdot U_q b_- + U_q b_-  \cdot U_+^{(1)} \big) \quad\mbox{where} \quad U_-^{(1)} = \bigoplus_{k\in {\mathbb N}} {\mathbb C} \cdot\tx_{-k}^-\ ,\quad U_+^{(1)} = \bigoplus_{k\in {\mathbb N}} {\mathbb C} \cdot\tx_{-k-1}^+.
\eeqa
The map $h_q: U_q b_- \rightarrow U_q\tilde h$ is defined as the projection along the last two summands. 
This map alone is not a homomorphism but it becomes so after restricting to the image of $t$ which is of degree 0 according to~\eqref{eq:deg-t}. We then notice by~\eqref{degxh} that any element of degree 0 should have a balanced number of $\tx^{\pm}$, and thus it can be written  as an element of either $U_q\tilde h \oplus  U_-^{(1)}\cdot U_q b_-$ or $U_q\tilde h \oplus U_q b_-  \cdot U_+^{(1)}$. The rest of the proof that $h_q\circ t$ is a homomorphism follows the end of the proof of~\cite[Lem.\,3]{FR99}. 

Finally, the $q$-character map $\chi_q$ from~\eqref{eq:char} is defined by~\cite[Def.\,2 \& Lem.\,3]{FR99}
 \beqa
 \chi_q(V)=h_q \circ t_V(z)\ .\label{def:charV}
 \eeqa
From~\eqref{eq:tVgen}, it follows that $\chi_q(V) = \sum_p m_p\chi_q(V_p)$,
with $V_p$ from~\eqref{eq:Vptens} and $\chi_q(V_p)= \prod_{i=1}^N\chi_q(V_{v_i}^{(j_i)})$, and  $\chi_q(V_{v}^{(j)})$ is
uniquely determined recursively through the image via $h_q$ of the TT-relation~\eqref{TTnorm-relj}
\beqa
\chi_q\bigl(V_v^{(j)}\bigr)= \chi_q\Bigl(V_{vq^{-j+\h}}^{(\h)}\Bigr)\chi_q\Bigl(V_{vq^\h}^{(j-\h)}\Bigr) - \chi_q\bigl(V_{vq}^{(j-1)}\bigr)\ \quad \mbox{with}\quad \chi_q(V_v^{(0)})=1\ ,\quad v\in{\mathbb C}^* \ .\label{eq:chichirec}
\eeqa
It is thus sufficient to compute the fundamental $q$-character $\chi_q\bigl(V_{v}^{(\h)}\bigr)$.
  We get: 
		  \beqa
	 \chi_q(V_v^{(\h)}) = h_q\Big( t_v^{(\h)}(u)\Big) 
     &=&  \Omega\Big(\tk_1^-(uv) + \tk_2^-(uv)\Big)\ \in U_q\tilde{h}[[u^{-2}]]\label{eq:chihalf}\\
		 &=& \Big( Y_{q^{-2}v^{-2}}+Y_{v^{-2}}^{-1} \Big) \Bigg|_{z\rightarrow -u^{-2}}\ ,\nonumber
		 \eeqa
using~\eqref{eq:tvhalf}, and with the generating function introduced in~\eqref{eq:chiY}. 
The above result matches with~\cite[Eq.\,(4.3) for $r=1$, $a=q^{-2}v^{-2}$]{FR99}. 
		Solving~\eqref{eq:chichirec} recursively width the  initial condition~\eqref{eq:chihalf}, we obtain:
		 \beqa\label{eq:chij}
		 \chi_q(V_v^{(j)}) = \prod_{k=1}^{2j} Y_{v^{-2}q^{2j-2k-1}}
		 \Bigg(\sum_{\ell =0}^{2j}\prod_{k=1}^\ell A^{-1}_{v^{-2}q^{2j-2k}}\Bigg) \Bigg|_{z\rightarrow -u^{-2}}
    \eeqa
    where
    \beq
         A_a = Y_{aq} Y_{aq^{-1}}
        = \Phig(z^{-1}a^{-1})\ ,\label{eq:defAphi}
		 \eeq
		 which agrees with~\cite[Eqs.\,(4.3) for $r=2j$, $a=q^{-2}v^{-2}$]{FR99}.
    \begin{rem}\label{rem:Aphi} 
    Let $V=\bigoplus_{\ell}V_\ell$ be the (generalized) eigenspace decomposition with respect to the action on $V$ of the coefficients $\psi_k$ and $\phi_{-k}$ in the power series $\Psig(z)$ and $\Phig(z)$ from~\eqref{def:psi},~\eqref{def:phi}.  The standard $q$-character of $V$ takes the form~\cite{FR99,FMu99}:
    \beqa
    \chi_q(V) = \sum_{\ell} \dim(V_\ell) M_\ell
    \eeqa
where  $M_{\ell}$ are monomials in $Y^{\pm 1}_{a}$ for some $a\in{\mathbb C}^*$. By~\cite[Prop.\,2.4]{FMu99}, there exists a one-to-one correspondence between  the $M_{\ell}'s$ and the common (generalized)  eigenvalues of $\Psig(z)$ and $\Phig(z)$.
    In particular, with respect to our basis $\{{\mathsf v}_0,{\mathsf v}_1,...,{\mathsf v}_{2j}\}$ in~\eqref{app:actKEF} and the decomposition $V=\bigoplus_{\ell=0}^{2j}{\mathbb C}{\mathsf v}_\ell$, the eigenvalues of $\Psig(-u^{-2})$ and $\Phig(-u^{-2})$ have multiplicity $1$, and are given  by the same expression~\eqref{app:psivn} but considered as power series in $u^{-2}$ and $u^2$, respectively. The correspondence associates each product (within the brackets) in
\beqa
M_\ell = \Bigr(\prod_{k=1}^{2j}Y_{aq^{2j-2k+1}}\Bigl)\Bigr( \prod_{k=1}^{\ell}Y_{aq^{2j-2k+3}}^{-1}Y_{aq^{2j-2k+1}}^{-1}\Bigl)\Bigg|_{z\rightarrow -u^{-2}}\ \quad \mbox{with}\quad a=q^{-2}v^{-2}\ \label{eq:Ml}
\eeqa
to the polynomials $R(z)$ and $Q_\ell(z)$ in~\eqref{app:R},~\eqref{app:Q}, respectively,  where the arguments of $Y$'s are the inverse roots of these polynomials.
For spin-$\h$, we see from~\eqref{eq:chihalf} two such monomials corresponding to the two distinct  common eigenvalues (of multiplicity $1$) associated with $\ell=0,1$.
    \end{rem}
\begin{example}\label{ex:q-char-spin-1}
 The image of ${\chi}_q$ on the spin-$1$ evaluation representation is given by:
 \beqa
 \chi_q(V_v^{(1)}) = Y_{aq}Y_{aq^{-1}} + Y_{aq^{-1}}Y_{aq^{3}}^{-1} + Y_{aq}^{-1}Y_{aq^{3}}^{-1} \Big|_{z\rightarrow -u^{-2}}\  
 \eeqa
 with $a=q^{-2}v^{-2}$.  Each of the the three monomials occuring in the above expression correspond, respectively, to~\eqref{eq:Ml} for $\ell=0,1,2$, weighted by  $\dim(V_\ell)=1$. 
\end{example}

		     \subsection{Refined $q$-characters of $\Loop$ from K-operators} \label{sec:refined-q-char}
             
             Recall the subalgebras  $U_q\tilde{h}\subset {}^{Dr}\!\Loop$ from~\eqref{eq:Uqth-def} and   $O_q\widetilde{H}\subset O_q$ from~\eqref{eq:OqHt} generated by imaginary root vectors of $\Loop$ and $O_q$, respectively. There are two types of  $q$-characters that may be considered:
		    \beqa
		    \chi_q^\iota \colon \ [\mathrm{Rep}O_q] \ \rightarrow \ U_q\tilde{h}[[z]]\ \quad \mbox{or} \quad\quad \acute{\chi}_q\colon \ [\mathrm{Rep}\Loop] \ \rightarrow \ O_q\tilde{H}[[z]]\ .\label{eq:typechi}
		   \eeqa
		   The first type is the so-called `boundary' $q$-character, which is just a linear map\footnote{Note that $[\mathrm{Rep}\, O_q]$ has no natural  multiplication because $\mathrm{Rep}\, O_q$ is not a monoidal category, rather a module category over $\mathrm{Rep}\,\Loop$.},
		   defined in~\cite[Def.\,3.1]{LPr25} specialized to $O_q$. Here, we study the second type and construct the ring homomorphism
		  \beqa
		  && \acute{\chi}_q\ : \ [\mathrm{Rep}\Loop] \ \rightarrow \ \acute{\cal Y} = {\mathbb Z}[\acute{Y}_{a}^{\pm 1}]_{a\in{\mathbb C}^*}\  \in \ O_q\widetilde{H}[[z]] \label{def:refchar}
		  \eeqa
		  with
          \beqa
		  \acute{Y}_a = \exp\left(-(q-q^{-1}) \sum_{k=1}^{\infty} \frac{z^ka^k} {q^{k} + q^{-k}}{\acute H}_k\right)  \ \in \ O_q\widetilde{H}[[z]]\ .\label{eq:Ya}
		  \eeqa
The  image  $\acute{\chi}_q(V)$ is called the {\it refined} $q$–character of the finite-dimensional representation $V$ of $\Loop$. 

The refined $q$-characters are defined and constructed within the framework of the FRT presentation of $O_q$ given in Theorem~\ref{thm:FRTOqDr}, using two main ingredients. The first one is a ring homomorphism
    \begin{equation}\label{eq:t-def}
         \acute{t} : \; [\mathrm{Rep}\Loop] \to O_q[[z]]
         \end{equation}
         called \textit{the universal boundary transfer matrix map}.  
         The image of $\acute{t}$ on a finite-dimensional representation $V$ of $\Loop$, denoted $\acute{t}_V(z)$, is a commutative subalgebra of $O_q$. We refer to~\cite[Chap.\,4]{L23} and~\cite[Sec.\,3.4]{BGL25} for the case of  simple evaluation modules over $\Loop$. The second ingredient is an analog of the Harish-Chandra map $h_q$
         discussed below~\eqref{eq:hq-def}:
\beqa
\acute{h}_q: O_q \to O_q \tilde{H}\ .
\eeqa
Below, we introduce and discuss both the ingredients in detail. 

\subsubsection{Definition of $\acute{t}$}\label{sec:def-acute}
To define $\acute{t}$, let us first introduce the (normalized) spin-$\h$ universal boundary transfer matrix 
in the root vectors of $O_q$:
 \beqa
		\acute{t}_v^{(\h)}(u)  &=& \acute{t}_{V_{v}^{(\h)}}(u) \label{eq:tcuteVhalf} \\ &:=& 
        \normalfont{\text{tr}}_{{\mathbb C}^2}\bigl(\bar K^{+}(uvq^{\h})\bar{\bold  K}^{\{r\}}(uvq^{\h})\bigr)  \quad \in \quad O_q((u^{-2}))
		  	\ ,\label{eq:acutet12} 
		 \eeqa
where the normalized fundamental K-operator for $O_q$ and normalized `dual' K-matrix are given by
          \beqa
         &&\qquad  \bar {\bold K}^{\{r\}}(u) = \bar\nu(u){\bold K}^{\{r\}}(u) \ ,\label{eq:Krbar12}\\
		  &&\qquad  \bar K^{+}(u)=\bar\nu_+(u) \frac{(q+q^{-1})}{(q^2u^2-q^{-2}u^{-2})}	\begin{pmatrix}
		  	uq \ov{\varepsilon}_+ + u^{-1}q^{-1} \ov{\varepsilon}_- & \ov{k}_+ (q^2u^2-q^{-2}u^{-2})/(q-q^{-1}) \\ 
		  	\ov{k}_- (q^2u^2-q^{-2}u^{-2})/(q-q^{-1}) & uq\ov{\varepsilon}_- + u^{-1}q^{-1} \ov{\varepsilon}_+
            \end{pmatrix}
		  \label{eq:Kpcute}\ ,
		  \eeqa
with~\eqref{Kr-Oq} and $\ov{\varepsilon}_\pm,\ov{k}_\pm \in{\mathbb C}$, and the normalization factors $\bar\nu(u)$ and $\bar\nu_+(u)$ are scalar power series in $u^{-1}$ defined by
\beqa
&&\pi^{\h}(\mu(u^2q)) \bar\nu(u) \bar\nu(uq) d_\rho^{(\varepsilon_+,\varepsilon_-)}(u) = 1 \ ,\label{eq:relmu1}\\
&&\pi^{\h}(\mu(u^2q^2)) \bar\nu_+(u) \bar\nu_+(uq) d_{\bar\rho}^{(\bar\varepsilon_+,\bar\varepsilon_-)}(u) = 1 \ ,\label{eq:relmu2}
\eeqa
where $d_c^{(a,b)}(u)$ is defined by~\eqref{eq:durho} and:
\beqa
\pi^{\h}(\mu(u)) = u^{-1}q^{-\h}  \exp \sum_{k=1}^\infty \frac{(q^{2k}+q^{-2k})}{1+q^{2k}}\frac{u^{-2k}}{k}\ .\label{eq:defmu}
\eeqa
The above functional relations have a  unique  (up to an overall sign) solution, see e.g.~\cite[Lem.\,6.1]{LBG}.
\begin{rem}  We notice that ${\bold K}^{\{r\}}(u)$ satisfies the same $q$-determinant relation~\eqref{qdetOq} as ${\bold K}^{\{a\}}(u)$. The `normalized K-operator' terminology is due to  the fact that $\bar{\bold K}^{\{r\}}(u)$ satisfies the normalized $q$-determinant relation, thanks to~\eqref{eq:relmu1}, 
\beqa
\normalfont{\text{tr}}_{12}\big({\cal P}^{-}_{12} \bar{\bold K}^{\{r\}}_1(u)  {\cal R}(u^2q)  \bar{\bold K}^{\{r\}}_2(u q) \big) &=& 1\ ,
\eeqa
where ${\cal R}(u_1/u_2)$ is the evaluation of the universal R-matrix, see~\cite[Sec.\,4.3]{Boos2012}:
$$(\pi^\h_{u_1} \otimes \pi^\h_{u_2}) (\mathfrak R) =\pi^{\h}(\mu(u_1/u_2)) R(u_1/u_2)\ .$$
\end{rem}

Now, we define the image $\acute{t}$ on the spin-$j$ evaluation representations:
\beqa
\acute{t}_{V_{v}^{(j)}}(u) := \acute{t}_v^{(j)}(u)\ \label{eq:tcutej}
\eeqa
defined recursively
\beqa
\label{TTacute-j}
		  	\acute{t}_{v}^{(j)}(u) =   \acute{t}_{vq^{-j+\h}}^{(\h)}(u) \acute{t}_{vq^\h}^{(j-\h)}(u)- \acute{t}_{vq}^{(j-1)}(u) \quad \mbox{with}\quad \acute{t}_{v}^{(0)}(u)=1\quad \mbox{ for all $j \in  \h{\mathbb N}^*$}\ .
\eeqa
Furthermore, 
\begin{itemize}
\item
for each $V_p$ of the form~\eqref{eq:Vptens}, 
we define $\acute{t}_{V_p}(z):= \prod_{i=1}^N\acute{t}_{V_{v_i}^{(j_i)}}(z)$;  

\item for each short exact sequence \ $0 \rightarrow V_p \rightarrow W \rightarrow V'_{p'} \rightarrow 0$, we define $t_{W}(z) := t_{V_p}(z)+t_{V'_{p'}}(z)$;

\item as every finite-dimensional $\Loop$-module $V$ 
admits a finite Jordan-Hölder series whose factors are of the form $V_p$ from~\eqref{eq:Vptens}, $p\in{\mathbb N}^*$, with multiplicities $m_p$, we thus define
$\acute{t}_V(z) := \sum_p m_p \acute{t}_{V_p}(z)$.
\end{itemize}
By construction, we have thus defined  a ring homomorphism~\eqref{eq:t-def}.

For our further purpose, we need to identify the subset of universal boundary transfer matrices of degree~$0$ with respect to the  grading map $\deg^{O_q}$ introduced in Remark~\ref{rem:Oq-grading}. 
\begin{lem}\label{lem:degtcute}
We have
\beqa
\deg^{O_q}(\acute{t}_v^{(\h)}(u)) = 0 \quad\mbox{iff} \quad \ov{\varepsilon}_\pm=0\ .
\eeqa
\end{lem}
\begin{proof} Recall Remarks~\ref{rem:Oq-grading},~\ref{rem:Oq-gradingbar} and the definition of ${\bold K}^{\{r\}}(u)$ in~\eqref{Kr-Oq}. Using~\eqref{eq:degOq} and \eqref{eq:degBbar}, we have $\deg^{O_q}\big({\bold K}^{\{r\}}(u)_{ii}\big)=1$ and $\deg^{O_q}\big({\bold K}^{\{r\}}(u)_{ij}\big)=0$ for $i\neq j$. Expanding the r.h.s.\,of~\eqref{eq:acutet12}, it follows that the terms ${\bold K}^{\{r\}}(u)_{ii}$ do not contribute iff $\ov{\varepsilon}_\pm=0$.
\end{proof}

\begin{cor}\label{cor:degtcute}
    For $\ov{\varepsilon}_\pm=0$, the image of the homomorphism $\acute{t}$ is of degree $\deg^{O_q}$ equal 0.
\end{cor}

\subsubsection{Definition of $\acute{h}_q$}\label{sec:def-acute-h}
 The second ingredient is the map
$\acute{h}_q: O_q \to O_q \tilde{H}$, which is  a projection  in the following sense. According to the PBW basis of $O_q$ with  ordering $B_{n\delta+\alpha_0}< B_{m\delta}  < B_{n\delta+\alpha_1}$~\cite[Thm.\,5.3]{Zh26b},
we can write:
 \beqa
\qquad O_q &=& O_q\tilde H \oplus \big(  O_-^{(1)}\cdot O_q + O_q\cdot O_+^{(1)} \big) \quad\mbox{where} \quad O_-^{(1)} = \sum_{k\in {\mathbb N}} {\mathbb C}\cdot  B_{n\delta+\alpha_0}\ ,\quad O_+^{(1)} = \sum_{k\in {\mathbb N}} {\mathbb C} \cdot B_{n\delta+\alpha_1} \ .\label{eq:decompO}
\eeqa
The map $\acute{h}_q$  is defined as the projection along the last two summands.
As we show below, this map becomes a  homomorphism after restricting its domain to the subspace of elements of $O_q$ of degree $0$. 

\begin{defn}\label{def:chicute} The map $\acute{\chi}_q\colon \ [\mathrm{Rep}\Loop] \ \rightarrow \ O_q\tilde{H}[[z]]$ is the composition of \ $\acute{t} : \; [\mathrm{Rep}\Loop] \to O_q[[z]]$  for $\ov{\varepsilon}_\pm=0$ and 
$\acute{h}_q: O_q \to O_q \tilde{H}$.
\end{defn}
\begin{thm}\label{thm:homhcute} The map $\acute{\chi}_q\colon \ [\mathrm{Rep}\Loop] \ \rightarrow \ O_q\tilde{H}[[z]]$ is a ring homomorphism.
\end{thm}
\begin{proof}  By construction, $\acute{t}$ is a ring homomorphism. The addition operation and the unit of $O_q$ are preserved by $\acute{h}_q$ due to the decomposition~\eqref{eq:decompO} and the fact that the unit belongs to $O_q\tilde{H}\subset O_q$. It thus remains to prove that $\acute{h}_q$ restricted on the image of $\acute{t}$ preserves the multiplication. Indeed, consider two elements $a,b\in O_q$ of degree $0$, where we 
recall Corollary~\ref{cor:degtcute}. Because any element of degree $0$ should have a balanced number of $B_{n\delta + \alpha_i}$, $i=0,1$, recall  Remark~\ref{rem:Oq-grading}, it can be written  as an element of either 
 $O_q\tilde H \oplus  O_-^{(1)}\cdot O_q$ 
or $O_q\tilde H \oplus O_q\cdot O_+^{(1)}$.  Using~\eqref{eq:decompO}, we write  $a=a_0+a_1$ with $a_0=\acute{h}_q(a)$, $a_1\in O_-^{(1)} \cdot O_q$ and $b=b_0+b_1$, with $b_0=\acute{h}_q(b)$, $b_1\in  O_q\cdot O_+^{(1)}$. We then have $\acute{h}_q(ab)=\acute{h}_q(a_0b_0+ a_1b_0 + a_0b_1 + a_1b_1)= a_0b_0$, because $a_1b_1$ belongs to both $O_-^{(1)} \cdot O_q$ and  $O_q\cdot O_+^{(1)}$, thus the multiplication is preserved.
\end{proof}

\subsubsection{Calculation of the refined $q$-character $\acute{\chi}_q$}
The refined $q$-character $\acute{\chi}_q(V_{v}^{(\h)})$ is now explicitly computed using the FRT presentation of Theorem~\ref{thm:FRTOqDr}. 
From~\eqref{eq:tcuteVhalf} and Definition~\ref{def:chicute}, we have:
		  \beqa
	 \acute{\chi}_q(V_v^{(\h)}) = \acute{h}_q\Big( \acute{t}_v^{(\h)}(u)\big|_{\bar\varepsilon_\pm=0}\Big)  \ \label{eq:acutechi12}\ .
		 \eeqa

\begin{prop}\label{prop:acutechi-half} Let $\ov{k}_\pm \in{\mathbb C}$ and $k_\pm\in\mathbb{C}^*$.  The refined $q$-character $\acute{\chi}_q(V_{v}^{(\h)})$ is given by
\beqa
\acute{\chi}_q(V_v^{(\h)})    
	 = \kappa(a)  \Bigg[ \frac{\overline{k}_+}{k_+}\acute{Y}_{a} + \frac{\overline{k}_-}{k_-}\acute{Y}^{-1}_{aq^2 }\Bigg]\Bigg|_{z\rightarrow u^{-2}} \label{eq:chirefhalf} \ ,
		  \eeqa
with $a=q^{-2}v^{-2}$ and $\acute{Y}_a$ from~\eqref{eq:Ya}, and we introduce the power series in $u^{-2}$:
          \beqa
          \kappa(a) = \bar\zeta(u^{-2}qa)
          \eeqa
defined uniquely  (up to an overall sign)  by the functional equation:
\beqa
\bar\zeta(u^{-2})  \bar\zeta(q^{2}u^{-2}) =
\frac{c(u^2q^{-1})c(u^2q)}{c(u^2)^2}\ \quad \mbox{with}\quad c(x)=x-x^{-1}\ .\label{eq:zetzet}
\eeqa
\end{prop}
\begin{proof}   
Recall from the proof of Lemma~\ref{lem:degtcute} that the terms ${\bold K}^{\{r\}}(u)_{ii}$ do not contribute in the r.h.s.\,of~\eqref{eq:acutet12} if $\ov{\varepsilon}_\pm=0$. Thus, it is sufficient to compute the images of the off-diagonal entries  $\bold{K}^{\{r\}}(u)_{ij}$, $i\neq j$. 
From~\eqref{Kr-Oq}, we have:
 \beqa
\acute{h}_q \circ \big(\bold K^{\{r\}}(u)\big)_{21} &=&	\frac{1}{k_+(q+q^{-1})}{\bar B}(u^{-2}) \ \label{eq:hqK21}
\eeqa
and
\beqa
\acute{h}_q \circ \big(\bold K^{\{r\}}(u)\big)_{12} &=&\acute{h}_q \circ \frac{1 }{k_-(q+q^{-1})}\Bigl(  \frac{\rho^2\eta(u^{-2})^{-1}\eta(q^{2}u^{-2})^{-1}}{(q-q^{-1})^2}  {\bar B}(q^2u^{-2})^{-1}  - u^{-2} B^{-}(u^{-2}){\bar B}(u^{-2})  B^{+}(u^{-2})\Bigr)
	 \nonumber \\
	 &=& \frac{\rho^2\eta(u^{-2})^{-1}\eta(q^{2}u^{-2})^{-1}}{k_-(q+q^{-1})(q-q^{-1})^2}  {\bar B}(q^2u^{-2})^{-1} \ , \label{eq:hqK12}
	\eeqa 
 with $\eta(u^{-2})^{-1}$ introduced in~\eqref{etagenAp} with~\eqref{etaetaAp}, and  in the  2nd equality, we  used that by definition of $\acute{h}_q$:
 $$\acute{h}_q\bigl( B^{-}(u^{-2}){\bar B}(u^{-2}) B^{+}(u^{-2})\bigr)=0\ ,
 $$
   after expanding the product in the argument according to~\eqref{eq:decompO}.

From~\eqref{eq:acutet12} with~\eqref{eq:Kpcute} and using~\eqref{eq:hqK21} with~\eqref{eq:hqK12} it follows:	
		  \beqa
		  \acute{\chi}_q(V_v^{(\h)}) 
          &=& 
		   \bar\zeta(u^{-2}v^{-2}q^{-1}) \nonumber\\
           && \times\Bigg[\frac{\overline{k}_+}{k_+}\underbrace{\Big(\frac{\rho}{q-q^{-1}}\Big)^{-1} \eta(u^{-2}){\bar B}(u^{-2})}_{}
		  + \frac{\overline{k}_-}{k_-} \underbrace{\Big(\frac{\rho}{q-q^{-1}}\Big)
		  	\eta(q^{2}u^{-2})^{-1}{\bar B}(q^2u^{-2})^{-1}}_{} \Bigg]\Bigg|_{u\rightarrow uvq^{\h}}\nonumber\ , 
\eeqa
where we denote
\beqa
\bar\zeta(u^{-2}) =  \frac{\rho \eta(u^{-2})^{-1}}{(q-q^{-1})^2}\bar\nu_+(u)\bar\nu(u)\ .
\label{eq:zet}
\eeqa
 Using \eqref{BbarHm} and \eqref{eq:Ya} for the two underlined terms in brackets, we get~\eqref{eq:chirefhalf}.
Then, using~\eqref{etaetaAp} with~\eqref{eq:relmu1} and~\eqref{eq:relmu2}, we see that~\eqref{eq:zetzet} indeed holds at the specialization $\bar\varepsilon_\pm=0$. 
Finally, the equation~\eqref{eq:zetzet} has a  unique  (up to a sign) solution, as in~\cite[Lem.\,6.1]{LBG}.
\end{proof}

\begin{rem}
Note that~\eqref{eq:zetzet} implies, up to an overlall sign factor, that the power series $\bar\zeta(u^{-2})=\sum_{n=0}^{\infty} \bar\zeta_{n} u^{-2n}$ has leading term $\bar\zeta_0=1$. The explicit expressions for the coefficients $\bar\zeta_k$ can be derived recursively, along the lines described in Lemma~\ref{lem:eta-recur}.
Moreover, we observe that $\bar\zeta(u^{-2})$  -- and thus $\kappa(a)$ -- doesn't depend on the scalars $k_\pm,\bar k_\pm,\varepsilon_\pm$.
\end{rem}

\begin{rem}
Similarly to~\eqref{eq:defAphi}, we have the following relation:
    \beqa
\Bigl(\kappa(aq^{-1})\acute Y_{aq^{-1}}\Bigr)^{-1} \Bigl(\kappa(aq^{-1})^{-1}\acute Y_{aq}\Bigr)^{-1} = -q(q-q^{-1})^{-2}\Thgsr(za)\ 
    \label{eq:YYcute}
    \eeqa
where the generating function $\Thgsr(z)$ was introduced in~\eqref{newTheta} with~\eqref{currentTheta}, and whose relation to Drinfeld's currents of $\Loop$ is given in Corollary~\ref{cor:finnc}.
\end{rem}

Now using~\eqref{TTacute-j}, the image of $\acute{\chi}_q$ on the spin-$j$ evaluation representations $\acute{\chi}_q(V_{v}^{(j)})$ is
uniquely determined recursively by
\beqa\label{eq:chichirecOq}
\acute{\chi}_q(V_v^{(j)})= \acute{\chi}_q\Bigl(V_{vq^{-j+\h}}^{(\h)}\Bigr)\acute{\chi}_q\Bigl(V_{vq^\h}^{(j-\h)}\Bigr) - \acute{\chi}_q\Bigl(V_{vq}^{(j-1)}\Bigr)\ \quad \mbox{with}\quad \acute{\chi}_q(V_v^{(0)})=1\ .
\eeqa
\begin{example}\label{valchi1} Similarly to Example~\ref{ex:q-char-spin-1},
 the image of $\acute{\chi}_q$ on the spin-$1$ evaluation representation is 
 \beqa
\acute\chi_q(V_v^{(1)}) &=&
\acute\chi_q\Bigl(V_{vq^{-\h}}^{(\h)}\Bigr) \acute\chi_q\Bigl(V_{vq^\h}^{(\h)}\Bigr)  - \acute\chi_q\bigl(V_{vq}^{(0)}\bigr) \label{eq:chi-spin11stline}\\
&=&\kappa(aq)\kappa(aq^{-1})\Bigl( \frac{\bar k_+^2}{k_+^2}\acute Y_{aq}\acute Y_{aq^{-1}} + \frac{\bar k_+\bar k_-}{k_+k_-}\acute Y_{aq^{-1}}\acute Y_{aq^{3}}^{-1} + \frac{\bar k_-^2}{k_-^2}\acute Y_{aq}^{-1}\acute Y_{aq^{3}}^{-1} + \frac{\bar k_+\bar k_-}{k_+k_-}  \Bigr)\Bigg|_{z\rightarrow u^{-2}}\! \!\!\!\!\! -1 \ ,\label{eq:acutchi-spin1}
 \eeqa
 with the identification $a=q^{-2}v^{-2}$, where we have used~\eqref{eq:chirefhalf}. We notice that the arguments of ${\acute Y}^{\pm1}$ are exactly same as those of $Y^{\pm1}$'s in Example~\ref{ex:q-char-spin-1}. But in contrast to the ordinary spin-1 $q$-character, here contributions for the spin-0 subquotient are not cancelled out due to the $\kappa$ terms (even if $\bar{k}_{\pm}/k_{\pm}=1$). 
 
 For the values $j=3/2$ and $j=2$, see the corresponding refined $q$-characters in Appendix~\ref{app:refchar}.
\end{example}

\begin{example}\label{ex:acute-chi-j}
For ${\bar k}_+/k_+=1$ and ${\bar k}_-/k_-=0$ the first refined q-characters are given by:
\beqa
\acute{\chi}_q(V_v^{(\h)}) &=&\kappa(a)\acute{Y}_{a}\Big|_{z\rightarrow u^{-2}}\ ,\nonumber\\
\acute{\chi}_q(V_v^{(1)}) &=&\kappa(aq)\kappa(aq^{-1})\acute{Y}_{aq}\acute{Y}_{aq^{-1}}\Big|_{z\rightarrow u^{-2}} - 1\ ,\nonumber\\
\acute{\chi}_q(V_v^{(\tha)}) &=&\kappa(aq^2)\kappa(a)\kappa(aq^{-2})\acute{Y}_{aq^2}\acute{Y}_{a}\acute{Y}_{aq^{-2}} - \kappa(aq^2)\acute{Y}_{aq^2}- \kappa(aq^{-2})\acute{Y}_{aq^{-2}}\Big|_{z\rightarrow u^{-2}}\ ,\nonumber\\
\acute{\chi}_q\big(V_v^{(2)}\big)
    &=& \kappa(aq^{-3})\kappa(aq^{-1})\kappa(aq)\kappa(aq^{3})
        {\acute Y}_{aq^{-3}}{\acute Y}_{aq^{-1}}{\acute Y}_{aq}{\acute Y}_{aq^{3}}\Big|_{z\rightarrow u^{-2}}
     \nonumber\\
    &&- \kappa(aq)\kappa(aq^{3})
        {\acute Y}_{aq}{\acute Y}_{aq^{3}}
     - \kappa(aq^{-3})\kappa(aq^{3})
        {\acute Y}_{aq^{-3}}{\acute Y}_{aq^{3}}
    - \kappa(aq^{-3})\kappa(aq^{-1})
        {\acute Y}_{aq^{-3}}{\acute Y}_{aq^{-1}}
    \Big|_{z\rightarrow u^{-2}} + 1\ ,\nonumber
\eeqa
with the identification $a= q^{-2}v^{-2}$.
\end{example}

Let $V$ be a type~1 finite-dimensional representation of $\Loop$.
By Theorem~\ref{thm:homhcute}, and following the same line of reasoning as above that uses a finite Jordan-Hölder series of $V$ in terms of $V_p$ from~\eqref{eq:Vptens} and their multiplicities $m_p$,  we obtain 
 \beq\label{def:refcharV}
 \acute{\chi}_q(V)=\acute{h}_q \circ \acute{t}_V(z)\big|_{\bar\varepsilon_\pm=0}
 = \sum_p m_p\acute{\chi}_q(V_p)\ ,
 \eeq
 where $\acute{\chi}_q(V_p)= \prod_{i=1}^N\acute{\chi}_q(V_{v_i}^{(j_i)})$ for each $V_p$ of the form~\eqref{eq:Vptens},
 and $\acute{\chi}_q(V_{v}^{(j)})$'s are
uniquely determined by~\eqref{eq:chichirecOq}, which is based on~\eqref{eq:chirefhalf}.

\subsection{Relating $\acute{\chi}_q$ to the universal (tensor) K-matrix}\label{sec:defchi} 
Similarly to the standard $q$-characters~\eqref{def:charV} with~\eqref{eq:tV-def} that are given in terms of $\mathfrak{R}$, the actions of the two types of $q$-character maps~\eqref{eq:typechi} on $V$ are expected to admit expressions in terms of a universal K-matrix~\cite[Sec.\,1.6.1]{AV24}. Let us assume the existence of a universal (tensor) K-matrix in an appropriate completion $\widehat{\otimes}$ of the tensor product
\beqa
\mathfrak{K}\ \in O_q \ \widehat{\otimes}\  \Loop\ , \label{eq:univK}
\eeqa
 that respects the axioms in~\cite[Def.\,2.7]{LBG} involving the twist automorphism $\psi$ given by~\eqref{app:psiaut} and with the Drinfeld twist $J=1\otimes1$.
Let us consider the refined $q$-character map $\acute{\chi}_q$ from Definition~\ref{def:chicute}. 
Our goal is to conjecture an expression for $\acute\chi_q(V)$ -- recall the definition~\eqref{def:refcharV} -- that involves the universal (tensor) K-matrix~\eqref{eq:univK}. As a preliminary, we proceed in two steps.

The first step aims to relate our fundamental K-operator $\bold K^{\{r\}}(u)$ 
and its higher spin-$j$ K-operators generalizations to $\mathfrak{K}$. 
We recall from Section~\ref{sec:RmatOq} that  $\bold K^{\{r\}}(u)$ is the image by $\gamma \otimes \id$ of the fundamental K-operator ${\cal K}(u)$ for $\cal A_q$ given by~\eqref{K-Aq}, where $\gamma:{\cal A}_q \rightarrow O_q$  is the surjective homomorphism defined in~\eqref{mapDeldel}. Then applying  the map $\gamma$ to the equations in the conjecture~\cite[Conj.\,3.1]{BGL25}, see also comments above it, we conjecture:\footnote{We note that $\bar\nu(u)$ and ${\bold K}^{\{r\}}(u)$ are, respectively, 
the images via $\gamma$ of $\nu^{(\h)}(u)$ and ${\cal K}^{(\h)}(u)$  in~\cite[Conj.\,3.1]{BGL25}. 
Moreover, the relation~\eqref{eq:relmu1} is the image via $\gamma$ of~\cite[Eq.\,(3.5) \& (3.49)]{BGL25}.}:
\beqa
\bigl(\id \otimes (\pi^{\h}_{1} \circ \phi_{u^{-1}})\bigr)
(\mathfrak{K}) = \bar \nu(u) {\bold K}^{\{r\}}(u)\ ,\label{eq:kunivK12}
\eeqa
where $\pi^j_v$, for $v\in \mathbb{C}^*$, is the spin-$j$ finite-dimensional evaluation representation as in Appendix~\ref{sec:formevrep}, and the gradation map $\phi$ is from~\eqref{grad1}-\eqref{grad3}.
Moreover, the higher spin-$j$ generalizations of ${\bold K}^{\{r\}}(u)$, denoted  ${\bold K}^{\{r\}}_{(j)}(u)$, can be derived recursively using the fusion procedure from~\cite[Def.\,5.6]{LBG} under the map $\gamma \otimes \id$:\footnote{In this subsection we follow conventions on $\Loop$ from~\cite{LBG,BGL25}.}
\begin{equation*}
{\bold K}^{\{r\}}_{(j+\h)}(u) =  \mathcal{F}^{(j+\h)}_\fu \left[{\bold K}^{\{r\}}_{(\h)}(u q^{-j})\right]_1 R^{(\h,j)}(u^2 q^{-j+\h}) \left[{\bold K}^{\{r\}}_{(j)}(u q^{\h})\right]_2 \mathcal{E}^{(j+\h)}_\fu ,
\end{equation*}
where the $\Loop$-intertwiner $\mathcal{E}$ and its pseudo-inverse $\mathcal{F}$ are defined in~\cite[Eq.\,(3.14)\,\&\,Lem.\,3.1]{LBG}, and
the spin-$j$ matrix $R^{(\h,j)}(u)$ is given in~\cite{KR83}, see also explanations in~\cite[Ex.\,2.3]{BGL25}.
 Therefore, following~\cite[Cor.\,3.2]{BGL25}, we get 
\beqa
(\id \otimes (\pi^{j}_{1} \circ \phi_{u^{-1}}))(\mathfrak{K}) &:=& \bar {\bold K}^{\{r\}}_{(j)}(u)\label{eq:barKrj} \\
&=&\bar \nu^{(j)}(u) {\bold K}^{\{r\}}_{(j)}(u)\ ,\label{eq:kunivKj}
\eeqa
where ${\bold K}^{\{r\}}_{(\h)}(u)\equiv {\bold K}^{\{r\}}(u)$, $\bar\nu^{(\h)}\equiv \bar \nu(u)$ and
	\begin{align}
			\bar\nu^{(j)}(u) &= \displaystyle{\left ( \prod_{m=0}^{2j-1}   \bar\nu(uq^{j-\h-m})  \right) } \left( \displaystyle \prod_{k=0}^{2j-2} \prod_{\ell=0}^{2j-k-2} \pi^{\h} (\mu (u^2q^{2j-2-2k-\ell})) \right) \ 
		\end{align}
with $\pi^\h(\mu(u))$ given by~\eqref{eq:defmu}.

The second step aims to express our (normalized) universal boundary tranfer matrix $\acute t_v^{(j)}(u)$ given by~\eqref{eq:tcutej} in terms of $\mathfrak{K}$. For $j=\h$, it is immediate: we insert~\eqref{eq:kunivK12} in~\eqref{eq:acutet12} using~\eqref{eq:Krbar12} and get 
 \beqa
		\acute{t}_v^{(\h)}(u)  &=& 
        \normalfont{\text{tr}}_{{\mathbb C}^2}\bigl(\bar K^{+}(uvq^{\h})(\id \otimes (\pi^{j}_{v^{-1}} \circ \phi_{u^{-1}q^{-\h}}))(\mathfrak{K})\bigr)  
		  	\ ,\label{eq:acutetunivK12} 
		 \eeqa
using that $\pi_v\circ \phi_u = \pi_1\circ \phi_{vu}$.
The higher spin-$j$ generalization of~\eqref{eq:acutet12} that satisfies~\eqref{TTacute-j} is constructed as follows. We consider the (normalized) spin-$j$  dual K-matrix  $\bar K^{+(j)}(u) =\bar\nu_+^{(j)}(u)K^{+(j)}(u)$ that is built recursively using the fusion procedure~\cite[Eq.\,(2.83)]{BGL25}  with the initial condition $\bar K^{+(\h)}(u)\equiv \bar K^{+}(u)$ given by~\eqref{eq:Kpcute} and\footnote{With respect to the normalization factor $\nu_+(u)$ in~\cite[Sec.\,3.3]{BGL25}, we have $\bar\nu_+(u)=\nu_+(u) (u^2q^2-u^{-2}q^{-2})/(q+q^{-1})$. We also note that the relation satisfied by $\bar\nu_+(u)$ given by~\eqref{eq:relmu2} is equivalent to~\cite[Eq.\,(3.24)]{BGL25}.}
\begin{align} \label{exp-nupjh}
			\bar\nu_+^{(j)}(u) &= \displaystyle{\left ( \prod_{m=0}^{2j-1}   \bar\nu_+(uq^{j-\h-m})  \right) } \left( \displaystyle \prod_{k=0}^{2j-2} \prod_{\ell=0}^{2j-k-2} \pi^{\h} (\mu (u^2q^{2j-2-2k-\ell}))^{-1} \right) \ .
		\end{align}
Taking the image by $\gamma$ of~\cite[Eqs.\,(3.21) \& (3.32)]{BGL25} and using~\eqref{eq:tcutej},~\eqref{eq:barKrj}, it follows
\beqa
\acute{t}_v^{(j)}(u) = \normalfont{\text{tr}}_{{\mathbb C}^{2j+1}}\bigl(\bar K^{+(j)}(uvq^{\h})\bar{\bold  K}^{\{r\}}_{(j)}(uvq^{\h})\bigr) \ ,
\eeqa
which, by \cite[Eq.\,(3.32)]{BGL25}, satisfies the TT-relation~\eqref{TTacute-j}. According to Corollary~\ref{cor:degtcute} and Definition~\ref{def:chicute}, we apply $\acute{h}_q$ to the latter expression for $\bar\varepsilon_\pm=0$ together with~\eqref{eq:kunivKj}, and we get 
\beqa
\acute{\chi}_q(V_v^{(j)}) =  \acute{h}_q \circ
    \normalfont{\text{tr}}_{{\mathbb C}^{2j+1}}\Big(\bar K^{+(j)}(uvq^\h)\big|_{\bar\varepsilon_\pm=0}
		  		\bigl(\id \otimes (\pi^{j}_{v^{-1}}\circ \phi_{u^{-1}q^{-\h}})\bigr){\mathfrak K}
		  		\Big) \ ,
\eeqa
which satisfies~\eqref{eq:chichirecOq}. Then, we note that $\pi^j_{v^{-1}} = \pi^j_{v} \circ \psi$, using the automorphism $\psi$ in~\eqref{app:psiaut} and~\eqref{app:defev}. All together, this leads us to the following conjecture.

       \begin{conj}\label{conj:chiunivK}
          Let $V$ be an arbitrary type~1 finite-dimensional representation of $\Loop$. There exists a matrix 
          $\bar K^{+V}(u)\in \End(V)$   that satisfies the dual reflection equation~\cite[Eq.\,(13)]{Skly88} such that  the refined $q$-character map of Definition~\ref{def:chicute} acts as
        \begin{equation}
		  	\acute{\chi}_q(V) =  \acute{h}_q \circ
    \normalfont{\text{tr}}_{V}\Big(\bar K^{+V}(u)
		  		\big(\id \otimes (\pi^{V^\psi}\circ \phi_{u^{-1}q^{-\h}})\big){\mathfrak K}
		  		\Big)  \quad \in \quad O_q\widetilde{H}[[u^{-2}]] \ 
		  	\ \label{tguniv} 
		  \end{equation}
	provided
      \beqa
\deg^{O_q}\Bigr(\normalfont{\text{tr}}_{V}\bigl(\bar K^{+V}(u)
		  		\big(\id \otimes (\pi^{V^\psi}\circ \phi_{u^{-1}q^{-\h}})\big){\mathfrak K}\Bigl) =0\ ,
      \eeqa
  and  $\bar K^{+V}(u)$ at the specialization $V=V^{(\h)}_v$  equals $\bar K^{+}(uvq^{\h})\big|_{\bar\varepsilon_\pm=0}$ with~\eqref{eq:Kpcute}. It therefore depends on the parameters $\bar{k}_\pm$.
  Here, $\pi^{V^\psi}(x) =\pi^{V}\circ \psi(x)$, where the automorphism $\psi$ of $\Loop$ is given by~\eqref{app:psiaut}. 
          \end{conj}

To conclude this section, let us discuss the interpretation of our refined $q$-characters from the perspective of representation theory of $\Loop$.
As recalled in Remark~\ref{rem:Aphi},  the standard $q$-characters $\chi_q(V)$ keep track of the vector space decomposition with respect to the generalized eigenvalues of $\Psig(z)$ and $\Phig(z)$. 
However, it is a very non-trivial task to recover relations in the Grothendieck ring, as~\eqref{eq:chichirec}, from explicit expression of the $q$-characters.
This is easily seen for the $j=1$ expression  in Example~\ref{ex:q-char-spin-1}, and it gets more complicated for higher spins. 
  Instead, our two-parameter family of refined $q$-characters~\eqref{def:refcharV} keeps track of the information on the relations  in the Grothendieck ring $[\mathrm{Rep}\Loop]$.
   Indeed, this is due to the dependence on 
   the two ratios $\ov{k}_\pm/k_\pm$ and  normalization factors $\kappa(a)$ with shifted arguments in the final expression for $\acute\chi_q(V)$. 
  If one can calculate such refined $q$-character from the universal K-matrix, according to Conjecture~\ref{conj:chiunivK}, they should take the form as in~\eqref{eq:chirefhalf} and~\eqref{eq:acutchi-spin1}, and~\eqref{eq:chi32} and~\eqref{eq:chi2}.
  In this case, one can easily recover the relations.
  For instance, the Grothendieck ring relation~\eqref{eq:chichirecOq} at $j=1$ remains explicitly apparent in~\eqref{eq:acutchi-spin1}: the overall contribution with coefficient $\kappa(aq)\kappa(aq^{-1})$ keeps track of the  product $\bigl[V_{vq^{-1/2}}^{(1/2)}\bigr]\cdot \bigl[V_{vq^{1/2}}^{(1/2)}\bigr]$ occuring in the r.h.s.\ of~\eqref{eq:chichirecOq} for $j=1$, whereas the last term $1$ keeps track of $V^{(0)}$. More generally, we recover relations in $[\mathrm{Rep}\Loop]$ by simply counting the $\kappa$ terms along the following algorithm: 
  \begin{enumerate}
      \item Set $\bar{k}_+=1$ and $\bar{k}_-=0$, and expand the refined $q$-character~\eqref{tguniv} on a representation $V$ in terms of the generating functions $\acute{Y}_a$ from~\eqref{eq:Ya};
      
      \item Group terms according to the \textit{$\kappa$-string} of length $k$: 
      \beq\label{eq:kappa-string}
      Z\kappa(aq^{n_1})\kappa(aq^{n_2})\ldots \kappa(aq^{n_k})\ , \qquad Z\in\mathbb{Z}\ ,
      \eeq
       which is the coefficient of the corresponding $\acute{Y}$-string $\acute{Y}_{aq^{n_1}}\ldots \acute{Y}_{aq^{n_k}}$, see e.g.\ in~\eqref{ex:acute-chi-j} for $\acute{\chi}_q\big(V_v^{(2)}\big)$;
      
      \item Replace in the $\kappa$-string each $\kappa(aq^{n})$ by the element $\bigl[V_{vq^{-n/2}}^{(1/2)}\bigr]$ in $[\mathrm{Rep}\Loop]$, and multiply them according to the $\kappa$-string~\eqref{eq:kappa-string};

      \item Delete $\acute{Y}_a$ terms;

      \item And finally sum up such products of spin-$\h$ elements $\bigl[V_{\bullet}^{(1/2)}\bigr]$, where $\bullet$ is an appropriate evaluation parameter from the previous step, multiplied by $Z\in\mathbb{Z}$.

      \item The final result gives a relation in $[\mathrm{Rep}\Loop]$ expressing $[V]$ as a polynomial in $\bigl[V_{\bullet}^{(1/2)}\bigr]$'s with coefficients in $\mathbb{Z}$.
  \end{enumerate}
         
\vspace{5mm}

{\bf Acknowledgments:} We thank Guillaume Lemarthe for collaborations at the early stage of this work. We thank A.\,Appel who pointed to us the reference \cite{P24}, and discussions in February 2024, which partly motivated this work. We also thank S.\,Belliard for discussions around the Drinfeld type presentation of $O_q$, P.\,Terwilliger and B.\,Vlaar for discussions. P.B. thanks T.\,Prze\'zdziecki for discussions and the organizers of the  conference
`Algebraic structures in the Yang-Baxter equation' at Heriot-Watt University, June
2024, where most of the results of this paper -- except Section~\ref{sec:refqchar} -- were presented.
P.B. and A.M.G.  are supported by C.N.R.S.  The work of A.M.G. was also partially supported by l'Agence nationale de la recherche
(ANR grant ``New algebraic structures in quantum integrability: towards 3D" NASQI3D ANR-24-CE40-7252).


\begin{appendix}
	
\section{Drinfeld-Jimbo presentation of $U_q(\widehat{sl}_2)$ and its root vectors}\label{ap:Loop}
	Define the extended Cartan matrix $\{a_{ij}\}$ with $a_{ii}=2$ and $a_{ij}=-2$ for $i\neq j$. The quantum affine algebra $U_q(\widehat{sl}_2)$ is generated by the elements $\{\tk^{\pm\h}_j,\te_j,\tf_j\}$, $j\in \{0,1\}$ which satisfy the defining relations
	\beqa \tk_i\tk_j=\tk_j\tk_i\ , \quad \tk^\h_i \te_j=q^{\h a_{ij}}\te_j\tk_i^\h\ , \quad
	\tk_i^\h \tf_j=q^{-\h a_{ij}}\tf_j\tk_i^\h\ ,\quad
	[\te_i,\tf_j]=\delta_{ij}\frac{\tk_i-\tk_i^{-1}}{q-q^{-1}}\ \label{def:Uq1}
	\eeqa
	and (the $q$-Serre relations)
	\beqa [\te_i,[\te_i,[\te_i,\te_j]_{q}]_{q^{-1}}]=0\ \qquad \mbox{and}\qquad
	[\tf_i,[\tf_i,[\tf_i,\tf_j]_{q}]_{q^{-1}}]=0\ . \label{def:Uq2}\eeqa
	The product $\tk_0\tk_1=\tk_1\tk_0=C$ is the central element of the algebra. 
		The \textit{quantum loop} algebra $\Loop$ is defined as the quotient
	of $U_q(\widehat{sl}_2)$  by the relation: 
	\beqa
	C=1 \qquad \Longrightarrow \qquad \tk_1 = \tk_0^{-1} \ .\label{eq:k1kO}
	\eeqa

	We recall that the Hopf algebra structure of $\Loop$ is given by the
	coproduct 
    $\Delta: \Loop\to \Loop\otimes \Loop$, the antipode $S: \Loop\to \Loop$ 	and the counit ${\epsilon}: \Loop\to {\mathbb C}$ defined by
	\beqa\label{eq:coprod}
    \Delta(\te_i)&=&\te_i\otimes 1 +
	\tk_i^{-1}\otimes \te_i\ ,\nonumber \\
	\Delta(\tf_i)&=&\tf_i\otimes \tk_i + 1 \otimes \tf_i\ ,\label{coprod}\\
	\Delta(\tk_i)&=& \tk_i \otimes \tk_i\ ,\nonumber
	\eeqa
	\beqa S(\te_i)=-\tk_i\te_i\ ,\quad  S(\tf_i)=-\tf_ik_i^{-1}\ ,\quad  S(\tk_i)=\tk_i^{-1} 
	\label{antipode}\eeqa
	and
	\beqa {\epsilon}(\te_i)={\epsilon}(\tf_i)=0 \qquad \mbox{and}\qquad {\epsilon}(\tk_i)=1\ .
	\label{counit}\eeqa

Following~\cite[Ex.\,2.13]{LBG}, we also introduce the twist automorphism of $\Loop$, here denoted by $\psi$:
\beqa
&&\psi(\te_0)=\tf_1\ ,\quad \psi(\te_1)=\tf_0\ , \label{app:psiaut}\quad\psi(\tk_0)=\tk_1^{-1}\ ,\nonumber\\
&& \psi(\tf_1)=\te_0\ ,\quad \psi(\tf_0)=\te_1\ , \quad\psi(\tk_1)=\tk_0^{-1}\ .\nonumber
\eeqa
	
	\subsection{Gradation}\label{sec:grad}
We recall that	$\Loop$ can be regarded as a graded algebra with respect to the so-called {\it principal} or {\it homogeneous} gradation.   Define the algebra homomorphism ${\rm D}^{(\sfr_0,\sfr_1)}_z$:~
	$\Loop \rightarrow \Loop\otimes {\mathbb C}[z,z^{-1}]$ 
	such that (see e.g. \cite{grad})
	\begin{equation}
		{\rm D}^{(\sfr_0,\sfr_1)}_z(\te_i)=z^{\sfr_i} \te_i\,,~~~{\rm D}^{(\sfr_0,\sfr_1)}_z(\tf_i)=z^{-\sfr_i}\tf_i\,,
		~~~{\rm D}^{(\sfr_0,\sfr_1)}_z(\tk_i)=\tk_i\,\label{general-g1}\ ,
	\end{equation}
	where $\sfr_i$ are arbitrary real numbers.
	The {\it homogeneous} gradation is given by $\sfr_0=0,\sfr_1=1$. The {\it principal} gradation is given by $\sfr_0=\sfr_1=-1$. In the main text, we denote the homomorphism associated with the principal gradation by 
	\beqa
	\phi_{z}\equiv {\rm D}^{(\sfr_0=-1,\sfr_1=-1)}_z\ . \label{gradp}
	\eeqa

	\subsection{Root vectors}\label{ap:rootvecdef}
	Recall the root system of $\widehat{sl}_2$ defined by $\cal R$, see (\ref{rootsys}). It is the disjoint union of the positive root system 
\beqa
\roots_+=\{n\delta+\alpha_0,n\delta+\alpha_1, m\delta\,|\,n\in {\mathbb N}, m\in {\mathbb N}^*\} \label{rootsysplus}
\eeqa
and the negative root system $\roots_-=-\roots_+$. Denote by $\te_\gamma$ (resp.  $\tf_{-\gamma}$) for $\gamma \in {\cal R}_+$ (resp. $\gamma \in {\cal R}_-)$ the following system  of root vectors. The fundamental root vectors of $\Loop$ are given by: 
\begin{equation} \label{app:firstRV}
	\te_{\alpha_1} \equiv \te_1 \ , \qquad \te_{\alpha_0} \equiv \te_0  \ , \qquad 	\tf_{\alpha_1}\equiv \tf_1 \ , \qquad \tf_{\alpha_0}\equiv \tf_0  \ .
\end{equation}
Other root vectors are defined\footnote{Other systems of root vectors have been considered in the literature \cite{Dam93,Beck,BCP}. For a discussion around the various choices of definitions with respect to the braid group action, see \cite[Sec.\,3.1]{Dam93}.} through the following set of recursion relations \cite[Eq.\,(4.13)-(4.20)]{Boos2012}. For $k \in \mathbb{N}^*$, they read:
\begin{equation} \label{app:rootvect}
	\begin{aligned}
		\te'_{k\delta}&= q^{-1} [\te_{\alpha_1+(k-1)\delta}, \te_{\alpha_0}]_{q} \ ,\\
		\te_{\alpha_1+k \delta}&=[2]_q^{-1} [\te_{\alpha_1 +(k-1)\delta},\te'_\delta] \ ,\\
		\te_{\alpha_0+k \delta} &= [2]_q^{-1} [ \te'_\delta,\te_{\alpha_0+(k-1)\delta}]\ ,\\
		\tf'_{k \delta} &= q [\tf_{\alpha_0},\tf_{\alpha_1+(k-1)\delta}]_{q^{-1}} \ ,\\
		\tf_{\alpha_1+k\delta} &= [2]_q^{-1} [\tf'_\delta, \tf_{\alpha_1 +(k-1) \delta} ] \ ,\\
		\tf_{\alpha_0+k\delta} &= [2]_q^{-1} [\tf_{\alpha_0+(k-1)\delta}, \tf'_\delta] \, ,
	\end{aligned}
\end{equation}
and the imaginary root vectors $e_{k \delta}$, $f_{k\delta}$ are defined via the generating functions ($\lambda = q-q^{-1}$):
\begin{align*}
 1 + \lambda \displaystyle{\sum_{k=1}^{\infty} } \te'_{k\delta} z^{k} &= 	\exp \left(\lambda\displaystyle{\sum_{k=1}^{\infty}} \te_{k\delta} z^{k} \right) \ , \\
  1-\lambda \displaystyle{\sum_{k=1}^{\infty} } \tf'_{k\delta} z^{-k}  &=	\exp\left(- \lambda \displaystyle{\sum_{k=1}^{\infty}} \tf_{k\delta} z^{-k}\right) \ .
\end{align*}
In particular, one has $\te'_{\delta}=\te_{\delta}$, $\tf'_\delta= \tf_{\delta}$. For the imaginary root vectors of $\Loop$, one has:
\beqa
 [\te_{k\delta},\te_{\ell\delta}]=[\tf_{k\delta},\tf_{\ell\delta}]=0  \qquad  \mbox{and} \qquad [\te_{k\delta},\tf_{\ell\delta}]=0 \qquad  \mbox{for any $k,\ell$} . \label{comim}
\eeqa
The relations between the Cartan elements $\tk_i$ and the root vectors follow from (\ref{app:rootvect}) and (\ref{def:Uq1}), (\ref{def:Uq2}):

\begin{equation} \label{app:Cartroot}
	\begin{aligned}
		[\tk_i,\te'_{k\delta}]&=0=[\tk_i,\te_{k\delta}] \qquad\qquad \qquad
		[\tk_i,\tf'_{k\delta}]=0=[\tk_i,\tf_{k\delta}]\ ,\\
		\tk_i\te_{\alpha_j+k\delta}\tk_i^{-1}&=q^{a_{ij}}\te_{\alpha_j+k\delta}\ ,\qquad \qquad  \tk_i\tf_{\alpha_j+k\delta}\tk_i^{-1}=q^{-a_{ij}}\tf_{\alpha_j+k\delta} \ .
		\end{aligned}
\end{equation}

Several examples of PBW bases for $\Loop$ can be considered. For instance, a PBW basis is obtained for the  uniquely determined total ordering $<$ on ${\cal R}_+$ defined by 
\begin{equation}\label{eq:rootord}
	k\delta + \alpha_1 <  (k+1)\delta + \alpha_1 < m\delta < (m+1)\delta <  (\ell+1)\delta + \alpha_0 < \ell\delta + \alpha_0 \ , \quad k,\ \ell\in{\mathbb N}\ , \ m\in {\mathbb N}^*\ ,
\end{equation}

	\subsection{Evaluation representation  $\pi_{v}^{j}$ and its Drinfeld polynomial}\label{sec:formevrep} For $\Loop$, any finite-dimensional representation is either of type $(1,1)$ or $(-1,-1)$ in the terminology of~\cite{CP91}. Here, we consider type $(1,1)$, while type $(-1,-1)$ is obtained using a twist automorphism, see~\cite[Sec.\,3.3]{CP91}. 

  We first recall that the quantum algebra $U_qsl_2$ is a Hopf subalgebra of $\Loop$ generated by $\te_1,\tf_1$ and $\tk_1$, and for brevity we denote them by $E$, $F$ and $K$,
respectively. Evaluation representations $(V_{v}^{(j)},\pi_{v}^j)$, $j\in\h{\mathbb N}^*$, $v\in{\mathbb C}^*$ of $\Loop$ are obtained by pulling back representations of $U_qsl_2$ using the evaluation homomorphism ${\mathsf{ev}}_{v}: \Loop \rightarrow U_qsl_2$ with evaluation parameter $v\in {\mathbb C}^*$ combined with  $\pi^j: U_qsl_2 \rightarrow \End({\mathbb C}^{2j+1})$, the spin-$j$ representation of $U_qsl_2$.  In the present paper, we essentially use the conventions from~\cite[Sec.\,4.3]{Boos2012}. 
Firstly, we introduce ${\mathsf{ev}}_v = \varphi \circ \phi_v$ using the (principal) gradation map~\eqref{gradp} and Jimbo's homomorphism $\varphi$ such that
    \beqa
    \varphi(\tk_1) = K\ ,\quad \varphi(\tk_0) = K^{-1}\ ,\quad \varphi(\te_1) = \varphi(\tf_0) = E\ ,\quad \varphi(\te_0) = \varphi(\tf_1) = F \ .
    \eeqa
    We then compute the image via  ${\mathsf{ev}}_v$ of the root vectors~\eqref{app:rootvect}. We get:
    \beqa
    {\mathsf{ev}}_{v}(\te_{\alpha_1+k\delta}) &=& (-1)^k v^{-2k-1}K^{-k}E \ ,\nonumber\\
     {\mathsf{ev}}_{v}(\te_{\alpha_0+k\delta}) &=& (-1)^k v^{-2k-1}FK^{-k} \ ,\nonumber\\ 
     {\mathsf{ev}}_{v}(\te'_{k\delta}) &=& (-1)^{k-1} v^{-2k}q^{-k}[E,F]_{q^{k}}K^{-k+1} \ ,\nonumber\\
     {\mathsf{ev}}_{v}(\te_{k\delta}) &=& (-1)^{k-1} \frac{v^{-2k}q^{-k}}{(q-q^{-1})k}\Big( C_k - (q^k + q^{-k})K^{-k} \Big)\ ,\label{app:evek}\\
     {\mathsf{ev}}_{v}(\tf_{\alpha_1+k\delta}) &=& (-1)^k v^{2k+1}FK^{k} \ ,\nonumber\\
     {\mathsf{ev}}_{v}(\tf_{\alpha_0+k\delta}) &=& (-1)^k v^{2k+1}K^{k}E \ ,\nonumber\\
      {\mathsf{ev}}_{v}(\tf'_{k\delta}) &=& (-1)^{k-1} v^{2k}q^k[E,F]_{q^{-k}}K^{k-1} \ ,\nonumber\\
       {\mathsf{ev}}_{v}(\tf_{k\delta}) &=& -(-1)^{k-1} \frac{v^{2k}q^{k}}{(q-q^{-1})k}\Big( C_k - (q^k + q^{-k})K^{k} \Big)\ ,\nonumber
    \eeqa
    where $C_k$ are central in $U_qsl_2$. They are defined using the generating function in the formal variable $y$:
    \beqa
    \sum_{k=1}^\infty (-1)^{k-1} C_k \frac{y^k}{k} = \log\big( 1 + Cy + y^2\big)
    \eeqa
    in terms of the Casimir element $C = q^{-1}K + qK^{-1} +(q-q^{-1})^2 EF$.

   Secondly, let $\{{\mathsf v}_0,{\mathsf v}_1,...,{\mathsf v}_{2j}\}$ be a basis  of ${\mathbb C}^{2j+1}$ for the spin-$j$ finite-dimensional representation of $U_qsl_2$ such that:  
    \beqa
    \pi^j(K){\mathsf v}_\ell = q^{2j-2\ell} {\mathsf v}_\ell\ ,\quad \pi^j(E){\mathsf v}_\ell = [\ell]_q [2j-\ell+1]_q{\mathsf v}_{\ell-1}\ ,\quad 
    \pi^j(F){\mathsf v}_\ell =  {\mathsf v}_{\ell+1}\ .\label{app:actKEF}
    \eeqa
    We recall that~\cite[Eq.\,(4.59)]{Boos2012}
    \beqa
    \pi^j(C_k)=\bigl(q^{(2j+1)k} + q^{-(2j+1)k}\bigr) I\!\!I_{2j+1}\ .\label{App:Ckvalj}
    \eeqa
    
    Altogether, the spin-$j$ evaluation representation $(V_{v}^{(j)},\pi_{v}^{j})$ of  $\Loop$ is obtained by
    \beqa
    \pi_{v}^j = \pi^j \circ {\mathsf{ev}}_{v} \ \label{app:defev}.
    \eeqa
	In particular, the spin-$\h$ evaluation representation of $\Loop$  on the root vectors takes the following form:
		\begin{align*}
		\pi_{v}^{\frac 12}(\te_{\alpha_1+k \delta})&=  (-1)^k v^{-2k -1}q^{-k} E_{12}  \ ,\\
		\pi_{v}^{\frac 12}(\te_{\alpha_0+k \delta})&=(-1)^k v^{-2k-1} q^{-k} E_{21}   \ ,\\
		\pi_{v}^{\frac 12}(\tf_{\alpha_1+k\delta})&= (-1)^k v^{2k+1} q^{k} E_{21}  \ ,\\
		\pi_{v}^{\frac 12}(\tf_{\alpha_0+k \delta})&=  (-1)^k v^{2k+1} q^{k} E_{12} \ ,\\  \label{app:piroot} \numberthis
		\pi_{v}^{\frac 12}(\te'_{k \delta})&= (-1)^{k-1} v^{-2k} q^{-k+1}(E_{11} - q^{-2} E_{22})  \ ,\\
		\pi_{v}^{\frac 12}(\te_{k \delta})&= (-1)^{k-1} v^{-2k}\frac{[k]_q}{k}  (E_{11} - q^{-2k}E_{22}) \ ,\\
		\pi_{v}^{\frac 12}(\tf'_{k \delta})&= (-1)^{k-1} v^{2k}q^{k-1} ( E_{11} -q^2 E_{22}) \ ,\\
		\pi_{v}^{\frac 12}(\tf_{k\delta})&=  (-1)^{k-1} v^{2k}\frac{[k]_q}{k} (E_{11} - q^{2k} E_{22}) \ ,\\
        \pi_{\alpha}^{\frac 12}(q^{\h h_1}) &= q^{\h}E_{11} + q^{-\h}E_{22} \ , \quad \tk_1=q^{h_1}\ ,
	\end{align*}
	where  $E_{ab}$ is the matrix with zero everywhere except 1 in the entry $(a,b)$.
    
The Drinfeld polynomial associated with $V_v^{(j)}$ is computed as follows (see \cite[Cor.\,4.2]{CP91} in the case of homogeneous gradation). Recall the Drinfeld generating function $\Psig(z)$ given by~\eqref{def:psi}. Firstly, we have:
\beqa
\pi_v^{j}\big(\Psig(-u^{-2})\big)\cdot {\mathsf v}_0 &=& \pi^j \circ {\mathsf{ev}}_{v} \circ\Omega ( \tk_1^{-1}) \times \exp\left((q-q^{-1}) \sum_{k=1}^{\infty} \pi^j \circ {\mathsf{ev}}_{v} \circ \Omega (\te_{k\delta}) (-u^{-2})^{k} \right)\cdot {\mathsf v}_0\ \label{app:DrPol1}\\
&=& q^{2j}\frac{P_j(q^{-1}v^{-2};q^{-2}u^{-2})}{P_j(q^{-1}v^{-2};u^{-2})}\cdot {\mathsf v}_0 \nonumber
\eeqa
with Drinfeld polynomial 
\beqa
P_j(x;z) = (1-q^{2j-1}xz)(1-q^{2j-3}xz)\cdots (1-q^{-2j+1}xz)\ ,\label{app:Drpoly}
\eeqa
where we used Table~\ref{tab:Omx},~\eqref{app:defev}, \eqref{app:evek}-\eqref{App:Ckvalj}
  and the fact that the evaluation map commutes with the anti-automorphism $\Omega$.
   More generally, on any vector ${\mathsf v}_\ell$ we obtain:
  \beqa
  \pi_v^{j}\big(\Psig(-u^{-2})\big)\cdot {\mathsf v}_\ell &=&  q^{2j-2\ell} \frac{R_j(q^{-1}u^{-2})}{R_j(qu^{-2})}\frac{Q_\ell(qu^{-2})}{Q_\ell(q^{-1}u^{-2})}\cdot {\mathsf v}_\ell \label{app:psivn}
  \eeqa
  where $R_j(z)$, $Q_\ell(z)$ are polynomials of degree $2j$ and $2\ell$, respectively, given by:
  \beqa
  R_j(z)&=& \prod_{k=1}^{2j}(1-aq^{2j-2k+1}z)\ ,\label{app:R}\\
  Q_\ell(z)&=& \prod_{k=1}^{\ell}(1-aq^{2j-2k+3}z)(1-aq^{2j-2k+1}z)\ \quad \mbox{with}\quad a=q^{-2}v^{-2}\ .\label{app:Q}
  \eeqa
  Here, the r.h.s.\,of~\eqref{app:psivn} is considered as a power series in $u^{-2}$. Let us note that the action of $\Phig(-u^{-2})$ given by~\eqref{def:phi}  on ${\mathsf v}_\ell$ also produces the r.h.s.\,of~\eqref{app:psivn}, but now considered as a power series in $u^2$.
\begin{rem}\label{app:remVV}
Note that our choice of coproduct and evaluation map in the present paper slightly differs from the one in~\cite{LBG,BGL25} by an automorphism of $\Loop$ denoted $\nu$, see~\cite[Eqs.\,(2.26)-(2.27)]{LBG}. The precise relationship between the evaluation representation ${\cal V}^{(j)}_{u}$, for $u\in{\mathbb C}^*$, for the conventions in~\cite{LBG,BGL25} and  $V^{(j)}_{v}$ in the present paper is established by comparing their Drinfeld polynomials.
It gives:
\beqa
V^{(j)}_{v} = {\cal V}^{(j)}_{vq^\h}\ .\label{app:VV}
\eeqa
\end{rem}

	\section{The algebra ${\cal A}_q$}\label{ap:AqRmat}
	In this appendix, the alternating central extension of the $q$-Onsager algebra ${\cal A}_q$, its PBW basis, R-matrix presentation and center ${\cal Z}$ are briefly reviewed.   We refer the reader to \cite{BSh1,BasBel,Ter21,Ter21d,LBG,L23,BGL25,Zh26a} for further details.

\subsection{Defining relations and PBW basis}
	\begin{defn}[\cite{BSh1}] ${\cal A}_q$ \label{DefAq} is the associative algebra with parameter $\rho\in{\mathbb C}^*$ and generators $\{{\tW}_{-k},{\tW}_{k+1},$
		${\tG}_{k+1},{\tilde{\tG}}_{k+1}|k\in {\mathbb N}\}$ satisfying the following relations:
		\begin{align}
			\big[{\tW}_0,{\tW}_{k+1}\big]=\big[{\tW}_{-k},{\tW}_{1}\big]=\frac{1}{(q+q^{-1})}\big({\tilde{\tG}_{k+1} } - {{\tG}_{k+1}}\big)\ ,\label{qo1}\\
			\big[{\tW}_0,{\tG}_{k+1}\big]_q=\big[{\tilde{\tG}}_{k+1},{\tW}_{0}\big]_q=\rho{\tW}_{-k-1}-\rho{\tW}_{k+1}\ ,\label{qo2}\\
			\big[{\tG}_{k+1},{\tW}_{1}\big]_q=\big[{\tW}_{1},{\tilde{\tG}}_{k+1}\big]_q=\rho{\tW}_{k+2}-\rho{\tW}_{-k}\ ,\label{qo3}\\
			\big[{\tW}_{-k},{\tW}_{-l}\big]=0\ ,\quad 
			\big[{\tW}_{k+1},{\tW}_{l+1}\big]=0\ ,\label{qo4}\quad \\
			\big[{\tW}_{-k},{\tW}_{l+1}\big]
			+\big[{{\tW}}_{k+1},{\tW}_{-l}\big]=0\ ,\label{qo5}\\
			\big[{\tW}_{-k},{\tG}_{l+1}\big]
			+\big[{{\tG}}_{k+1},{\tW}_{-l}\big]=0\ ,\label{qo6}\\
			\big[{\tW}_{-k},{\tilde{\tG}}_{l+1}\big]
			+\big[{\tilde{\tG}}_{k+1},{\tW}_{-l}\big]=0\ ,\label{qo7}\\
			\big[{\tW}_{k+1},{\tG}_{l+1}\big]
			+\big[{{\tG}}_{k+1},{\tW}_{l+1}\big]=0\ ,\label{qo8}\\
			\big[{\tW}_{k+1},{\tilde{\tG}}_{l+1}\big]
			+\big[{\tilde{\tG}}_{k+1},{\tW}_{l+1}\big]=0\ ,\label{qo9}\\
			\big[{\tG}_{k+1},{\tG}_{l+1}\big]=0\ ,\quad   \big[{\tilde{\tG}}_{k+1},\tilde{{\tG}}_{l+1}\big]=0\ ,\label{qo10}\\
			\big[{\tilde{\tG}}_{k+1},{\tG}_{l+1}\big]
			+\big[{{\tG}}_{k+1},\tilde{{\tG}}_{l+1}\big]=0\ .\label{qo11}
		\end{align}
	\end{defn}
	In the literature, the generators ${\tW}_{-k},{\tW}_{k+1},{\tG}_{k+1},{\tilde{\tG}}_{k+1}$ are called the alternating generators of ${\cal A}_q$. Note that another presentation of ${\cal A}_q$ has been introduced and studied, known as the `compact' presentation of the alternating central extension of the $q$-Onsager algebra;  See \cite[Thm.\,1.1, Thm.\,1.2]{Ter21c}.
	\begin{rem} \label{rem:auto}
		There exists an automorphism  $\sigma$ of ${\cal A}_q$ with $\sigma^2= \mathrm{id}$ such that \cite{BSh1}
	\beqa
		\sigma({\tW}_{-k})={\tW}_{k+1}\ , \quad \sigma({\tW}_{k+1})={\tW}_{-k}\ , \quad \sigma({\tG}_{k+1})=\tilde{\tG}_{k+1}\ , \quad \sigma(\tilde{\tG}_{k+1})={\tG}_{k+1} \ .\label{remomega}  
		\eeqa
		There exists an anti-automorphism  $\dagger$ of ${\cal A}_q$ with $\dagger^2= \mathrm{id}$ such that  \cite{Ter18}
		\beqa
		\dagger({\tW}_{-k})={\tW}_{-k}\ , \quad \dagger({\tW}_{k+1})={\tW}_{k+1}\ , \quad \dagger({\tG}_{k+1})=\tilde{\tG}_{k+1}\ , \quad \dagger(\tilde{\tG}_{k+1})={\tG}_{k+1} \ .\label{remdagger} 
		\eeqa
		 Define $\tau = \dagger \circ \sigma$. The anti-automorphism  $\tau$ with $\tau^2= \mathrm{id}$ is such that \cite{Ter21b}
		\beqa
		\tau({\tW}_{-k})={\tW}_{k+1}\ , \quad \tau({\tW}_{k+1})={\tW}_{-k}\ , \quad \tau({\tG}_{k+1})={\tG}_{k+1}\ , \quad \tau(\tilde{\tG}_{k+1})={\tilde\tG}_{k+1} \ .  \label{remtau}
		\eeqa
		\end{rem}
	
\begin{rem} For $q\rightarrow 1$, ${\cal A}_q$ specializes to the alternating central extension of the Onsager algebra \cite{Ter21e}. 
\end{rem}

A PBW basis for ${\cal A}_q$ is obtained by its alternating generators in any linear order $<$ such that	\cite[Thm.\,6.1]{Ter21}
	\begin{align}
		\tG_{i+1} < \tW_{-j} < \tW_{k+1} < \tilde{\tG}_{\ell+1}\qquad \qquad i,j,k, \ell \in \mathbb N.
		\label{eq:order}
	\end{align}

\subsection{R-matrix presentation and center $\cal Z$}\label{sec:RpresAq}	
		A R-matrix presentation for ${\cal A}_q$
		is known~\cite{BSh1}, taking  the form
	of a reflection algebra~\cite{Skly88}. Recall the R-matrix~\eqref{R}.  
	Introduce the convenient parametrization:
	\begin{equation}
		\rho=k_+k_-(q+q^{-1})^2\   \quad \text{with} \, \ \ k_\pm\in \mathbb{C}^*.  
		\label{rho}
	\end{equation}
	\begin{thm}[\cite{BSh1}]\label{def:Aq0} $\cal A_q$ admits a presentation in the form of a reflection algebra. Introduce the generating functions:
		\begin{align}
			{\tW}_+(u)=\sum_{k\in {\mathbb N}}{\normalfont \tW}_{-k}U^{-k-1} \ , \quad {\tW}_-(u)=\sum_{k\in  {\mathbb N}}{\normalfont \tW}_{k+1}U^{-k-1} \ ,\label{c1}\\
			\quad {\tG}_+(u)=\sum_{k\in {\mathbb N}}{\normalfont \tG}_{k+1}U^{-k-1} \ , \; \quad {\tG}_-(u)=\sum_{k\in {\mathbb N}}{\normalfont \tilde{{\tG}}_{k+1}}U^{-k-1} \ ,\label{c2} \; 
		\end{align}
		where the  shorthand notation  $U=(qu^2+q^{-1}u^{-2})/(q+q^{-1})$ is used. The defining relations are given by:
		\begin{align} 
			R(u/v)\ {\cal K}_1(u)\ R(uv)\ {\cal K}_2(v)\
			= \ {\cal K}_2(v)\  R(uv)\ {\cal K}_1(u)\ R(u/v)\ 
			\label{RE} 
		\end{align}
		with 		%
		\begin{equation}
			{\cal K}(u)=\begin{pmatrix} 
				uq \tW_+(u)-u^{-1}q^{-1}\tW_-(u) &\frac{1}{k_-(q+q^{-1})}\tG_+(u)+\frac{k_+(q+q^{-1})}{q-q^{-1}} \\ 
				\frac{1}{k_+(q+q^{-1})}\tG_-(u) +\frac{k_-(q+q^{-1})}{q-q^{-1}}& uq \tW_-(u) -u^{-1}q^{-1}\tW_+(u) 
			\end{pmatrix}  \label{K-Aq} \ .
		\end{equation}
		Recall that in our convention  the rational functions $U^{-k}$ are considered here as power series in $u^{-2}$.
	\end{thm}
	We call $\mathcal{K}(u)$ the fundamental K-operator for $\mathcal{A}_q$. Let us note that $U^{-1}$ can be written as a power series in $u^{-2}$, and~\eqref{c1} and~\eqref{c2} have no constant terms. It follows:
	\beqa
	\mathcal{K}(u) \in \mathcal{A}_q[[u^{-1}]]\otimes \End(\mathbb{C}^2)\ .\label{K-Aqexp}
	\eeqa

	Note that the defining relations~\eqref{RE}  are equivalent to the set of exchange relations (in the `symmetric' indeterminates $U,V$, symmetric under the exchange $ u \leftrightarrow q^{-1}u^{-1}$) between the generating functions (\ref{c1}), and (\ref{c2}):
	\begin{align}\label{eq1-Aq} \left [ \tW_\pm(u),\tW_\pm(v) \right ]&=0 \ ,  \\ 
		\label{eq2-Aq}\left [ \tW_+(u),\tW_-(v) \right] &+ \left [ \tW_-(u),\tW_+(v) \right ] =0 \ ,
	\end{align}
	\begin{align}
		\label{eq3-Aq}&	\left [ \tG_\epsilon(u),\tW_\pm(v) \right ] + \left [ \tW_\pm(u),\tG_\epsilon(v) \right ] =0, \quad  \epsilon=\pm \ , \\
		\label{eq4-Aq}&	\left [ \tG_\pm(u),\tG_\pm(v) \right ] =0\ , \\
		\label{eq5-Aq}&	\left [ \tG_+(u),\tG_-(v) \right ] + \left [ \tG_-(u), \tG_+(v) \right ]=0\ ,
	\end{align}
	\vspace{-5mm}
	\begin{align}
		\label{eq6-Aq}&(U-V) \left [ \tW_\pm(u),\tW_\mp(v) \right ]= \frac{(q-q^{-1})}{\rho(q+q^{-1})} \left ( \tG_\pm(u) \tG_\mp(v)-\tG_\pm(v) \tG_\mp(u) \right )\\ \nonumber
		&\hspace{4cm}	+ \frac{1}{(q+q^{-1})} \left ( \tG_\pm(u) -\tG_\mp(u)+\tG_\mp(v)-\tG_\pm(v) \right )\ ,\\
		\label{eq7-Aq}	&\tW_\pm(u)\tW_\pm(v)-\tW_\mp(u)\tW_\mp(v)+\frac{1}{\rho(q^2-q^{-2})}\left [ \tG_\pm(u),\tG_\mp(v)\right] \\ \nonumber
		&\hspace{4cm}	+ \frac{1- UV}{U-V}\left ( \tW_\pm(u)\tW_\mp(v)-\tW_\pm(v)\tW_\mp(u)\right) =0 \ , \\
		\label{eq8-Aq}&	U \left [ \tG_\mp(v), \tW_\pm(u) \right ]_q -V\left [ \tG_\mp(u),\tW_\pm(v)\right ]_q -(q-q^{-1})\left ( \tW_\mp(u)\tG_\mp(v) - \tW_\mp(v) \tG_\mp(u) \right ) \\ \nonumber
		&\hspace{4cm}	+ \rho \left ( U \tW_\pm(u) - V \tW_\pm(v) - \tW_\mp(u)+\tW_\mp(v) \right )=0 \ , \\ 
		\label{eq9-Aq}&	U \left [ \tW_\mp(u),\tG_\mp(v)\right ]_q - V \left [ \tW_\mp(v),\tG_\mp(u) \right ]_q -(q-q^{-1}) \left(\tW_\pm(u) \tG_\mp(v) - \tW_\pm(v) \tG_\mp(u) \right ) \\ \nonumber
		&\hspace{4cm}	+ \rho \left( U \tW_\mp(u)-V\tW_\mp(v) - \tW_\pm(u) +\tW_\pm(v) \right ) =0 \ .
	\end{align}
	\smallskip
	
In the R-matrix formalism, the center of $\mathcal{A}_q$, denoted $\cal Z$, is generated by the so-called quantum determinant ~\cite{Skly88} associated with the K-operator (\ref{K-Aq}), as we now recall. We define the following generating function
	\begin{equation}
		\Gamma(u)=\normalfont{\text{tr}}_{12}\big({\cal P}^{-}_{12} {\cal K}_1(u)  R(qu^2)  {\cal K}_2(u q) \big) \ .  \label{gammaform}
	\end{equation}

	\begin{prop}[\cite{BasBel}]  The quantum determinant of  the fundamental K-operator
		\begin{equation}
			\Gamma(u)=
			\frac{(u^2q^2-u^{-2}q^{-2})}{2(q-q^{-1})}\left( \Delta(u) - \frac{2\rho}{q-q^{-1}}\right) \ ,\label{gamma}
		\end{equation}
		with
		\begin{align} 
			\Delta(u)&= -(q-q^{-1})(q^2+q^{-2})\Big(\tW_+(u)\tW_+(uq) + \tW_-(u)\tW_-(uq)\Big) \label{deltau}  \\ \nonumber
			&+(q-q^{-1})(u^2q^2+u^{-2}q^{-2})\Big(\tW_+(u)\tW_-(uq) + \tW_-(u)\tW_+(uq)\Big)\\ 
			&- \frac{(q-q^{-1})}{\rho} \Big(\tG_+(u)\tG_-(uq) + \tG_-(u)\tG_+(uq)\Big) - \ \tG_+(u) - \tG_+(uq) - \tG_-(u) - \tG_-(uq)  \  ,\nonumber
		\end{align}
		is  such that  $\big[\Gamma(u),{\cal K}(v)_{mn}\big]=0$ where ${\cal K}(v)_{mn}$ denote the entries of the K-operator (\ref{K-Aq}). 
	\end{prop}
\begin{rem}
	Using the set of relations  that follow from (\ref{eq1-Aq})-(\ref{eq9-Aq}) -- see~\cite[Lem.\,B.1]{LBG} -- various equivalent expressions for $\Delta(u)$ can be derived. 
\end{rem}
	
In particular, expanding $\Delta(u)$ as a formal power series  in  $u^{-2}$, central elements are extracted:
\begin{equation}
	\Delta(u)  = \displaystyle{\sum_{k=0}^\infty}  u^{-2k-2} c_{k+1} \Delta_{k+1} \quad \mbox{with} \quad  c_{k} = - \frac{ (q+q^{-1})^{k} (q^{k} + q^{-k})} {q^{2k}} \ . \label{eq:ckp1}
\end{equation}
By the above Proposition, $[\Delta_{k+1},a]=0$ $\forall a\in {\cal A_q}$.
A formula expressing the central elements $\Delta_{k+1}$ in terms of the alternating generators is given in \cite[Lem.\,2.1]{BasBel}. 
By \cite[Thm.\,10.2(i)]{Ter21}, it follows:
\begin{prop}\label{DefZ} The center ${\cal Z}$ of ${\cal A}_q$ is generated by the elements $\{\Delta_{k+1}\}_{k=0}^{\infty}$ defined by~\eqref{eq:ckp1},~\eqref{deltau}.
\end{prop}

	\begin{rem}\label{rem:convention} In \cite{Ter21b}, the central elements $\{{\cal Z}^\vee_k\}_{k\in {\mathbb N}}$ are introduced. With respect to the conventions in \cite[Defs.\,8.1,\,8.3]{Ter21b} (using also \cite[Lem.\,8.8]{Ter21}), one has the correspondence:   \beqa
	&&q^{-1}u^{-2} \rightarrow t ,\qquad \qquad \qquad \qquad\quad  \ \  \rho \rightarrow -(q^2-q^{-2})^2 \ ,\nonumber\\ 
     &&{\tW}_\mp(u q^{\h}) \rightarrow S \mathcal{W}^\pm(S)\ ,\qquad  \qquad \quad{\tW}_\mp(u q^{-\h}) \rightarrow T \mathcal{W}^\pm(T) \nonumber\\
     && {\tG}_+(u q^{\h}) + \rho/(q-q^{-1}) \rightarrow \mathcal{G}(S) \ ,\quad {\tG}_-(u q^{-\h}) + \rho/(q-q^{-1}) \rightarrow \tilde{\mathcal{G}}(T) \nonumber
    \eeqa
    with
    \beqa
    S=\frac{q+q^{-1}}{q^{-1}t + qt^{-1}} \ ,\quad T=\frac{q+q^{-1}}{qt + q^{-1}t^{-1}} \ ,\nonumber
    \eeqa
    and
	\begin{equation}
		\frac{\Gamma(u q^{-\h})}{(q+q^{-1})^2(qu^2-q^{-1}u^{-2}) }  \rightarrow {\cal Z}^\vee(t)    \ 	\quad \mbox{where}\quad  {\cal Z}^\vee(t) = \sum_{k\in{\mathbb N}} {\cal Z}^\vee_k t^k\ .\label{gamZvee}
	\end{equation}
 Using (\ref{gamma}), it follows  
	\beqa
	\frac{(q-q^{-1})(q+q^{-1})^k(q^k+q^{-k})}{\rho}  \frac{\Delta_k}{2} \rightarrow {\cal Z}^\vee_k   \quad \mbox{and} \quad -\frac{\rho}{(q^2-q^{-2})^2} \rightarrow {\cal Z}^\vee_0 =1 \ .
	\eeqa
\end{rem}

Let us point out that the constant term of $\Delta(u) - 2\rho/(q-q^{-1})$ is 
	$- 2\rho/(q-q^{-1})$, so $\Gamma(u)$ invertible. It follows that $\mathcal{K}(u)$ is invertible too and its inverse is given by~\cite[Eq.\,(5.24)]{LBG}:
		\begin{equation}
			\left [	\mathcal{K}(u) \right ]^{-1} = \frac{(u^{-2}-u^2)}{\Gamma(u q^{-1})} \, \mathcal{K}(u^{-1}) \ .\label{KinvAq}
		\end{equation}
Note that when we write ${\tW}_\pm(u^{-1})$ and ${\tG}_\pm(u^{-1})$ it means replacing $u$ by $u^{-1}$ in $U$ and developing in power series in $u^{-2}$, this way we get 
\begin{equation}\label{eq:u-inv-convention}				
	{\tW}_\pm(u^{-1})= {\tW}_\pm(u q^{-1})\qquad \text{and} \qquad
	{\tG}_\pm(u^{-1})= {\tG}_\pm(u q^{-1})\ .
\end{equation}		
In particular, this shows that $\mathcal{K}(u^{-1})$ is equally in $ \mathcal{A}_q[[u^{-1}]]\otimes \End(\mathbb{C}^2)$.

\section{The generating function ${\bar B}(z)$}\label{ap:Bbarsol}
In this appendix, we determine the generating function ${\bar B}(z)$ in~\eqref{Bbargen}  for the standard case $\varepsilon_\pm=0$ and the non-standard case.
 
\subsection{The standard case}\label{sec:stTer}

Let us determine the auxiliary generating function ${\bar B}(z)$ for the special case  $\varepsilon_\pm=0$. For this choice, ~\eqref{BinGGred} simplifies to
\beqa
B(u^{-2}) = \xi^{-1}  {\bar B}^{(0)}(u^{-2})^{-1} {\bar B}^{(0)}(q^2u^{-2})^{-1} 
 \ ,\label{BinBbarBbar}
\eeqa
where ${\bar B}(z)|_{\varepsilon_\pm=0}$ is denoted ${\bar B}^{(0)}(z)$ with
\beqa
\bar B^{(0)}(z)= \sum_{n=0}^\infty{\bar B}^{(0)}_n z^n \ .\label{Bbarst}
\eeqa
The coefficients $\{{\bar B}^{(0)}_n\}_{n\in {\mathbb N}} \in O_q$ can be computed recursively.  Namely, using~\eqref{BinBbarBbar} and the same equation with the shift $u \rightarrow q u$, as well as \eqref{eq:comBbarBdel}, we get respectively:
	\beqa
	B(u^{-2}){\bar B}^{(0)}(q^2u^{-2}) = \xi^{-1}{\bar B}^{(0)}(u^{-2})^{-1} \quad \mbox{and} \quad 	B(q^{-2}u^{-2}){\bar B}^{(0)}(q^{-2}u^{-2}) = \xi^{-1}{\bar B}^{(0)}(u^{-2})^{-1}\ ,
	\eeqa
which gives
\beqa
B(u^{-2}){\bar B}^{(0)}(q^2u^{-2}) = B(q^{-2}u^{-2}){\bar B}^{(0)}(q^{-2}u^{-2})\ .\label{BbarBeqst}
\eeqa
Expanding both sides of the equation using (\ref{rvgen1}) and (\ref{Bbarst}) yields to the recurrence relations
\beqa
\sum_{k=0}^p (q^{2k}-q^{-2p}) B_{(p-k)\delta}{\bar B}^{(0)}_{k}=0\ \quad \mbox{for all \ $p\in {\mathbb N}$} \ ,\label{eq:recBbarBdel}
\eeqa
that determine uniquely  $\{{\bar B}^{(0)}_n\}_{n\in {\mathbb N}}$ in terms of $\{{B}_{n\delta}\}_{n\in {\mathbb N}}$. For instance, one gets for the first coefficients ${\bar B}^{(0)}_0=\frac{\rho}{q-q^{-1}}$, with the sign fixed by Definition~\ref{def:Bbar}, and
\beqa
{\bar B}^{(0)}_1&=&
	-(q+q^{-1}) B_{\delta} \ ,\nonumber\\
    {\bar B}_2^{(0)}&=&  -q^{-1}\frac{(q+q^{-1})^2}{(q^2+q^{-2})} B_{2\delta} - \frac{(q^3-q^{-3})}{(q^2+q^{-2})(q-q^{-1})^2}B_{\delta}^2 \nonumber \ .
\eeqa 

\subsection{The non-standard case}\label{sec:nstTer}
We now consider the generating function $\bar B(z)$ that solves~\eqref{BinGGred} for $\varepsilon_\pm\neq 0$.
We first define the invertible generating function 
\beqa
{\eta}(u^{-2}) = \sum_{n=0}^{\infty} \eta_k u^{-2k}\ , \quad \eta_k\in {\mathbb C}^*\ ,\label{etagenAp}
\eeqa
such that	
\beqa
\eta(u^{-2})^{-1}\eta(q^{2}u^{-2})^{-1} &=&\frac{ d(uq^{-1})^{(\epsp,\epsm)}_{\rho}}{(q+q^{-1})^2(u^2-u^{-2})}\ \label{etaetaAp} \\
&=&  \left(1+ \frac{ \epsp^2 + \epsm^2 + \epsp\epsm(u^2+u^{-2})}{(u^2-u^{-2})^2} \right) \ ,\nonumber
\eeqa
using $d(u)_\rho^{(\epsp,\epsm)}$ given in~(\ref{eq:durho}).
Our first task is to determine the scalars $\{\eta_k\}_{k\in{\mathbb N}}$ in terms of $\{\epsp,\epsm\}$:

\begin{lem}\label{lem:eta-recur}
The generating function $\eta(u)$  given by~\eqref{etagenAp} satisfying \eqref{etaetaAp} is uniquely (up to an overall sign) determined by the recursion
\beqa\label{eq:eta-recur}
\sum_{k=0}^p (q^{2k}-q^{-2p}) f_{p-k}\eta_{k}=0\ \quad \mbox{for all \ $p\in {\mathbb N}$} \ ,
\eeqa
where  $f(u) = \sum_{k=0}^\infty f_k u^{-2k}$ is  the r.h.s.\,of (\ref{etaetaAp}). 
Moreover,  the inverse generating function
\beqa
{\eta}(u^{-2})^{-1} = \sum_{n=0}^{\infty} \eta^{inv}_k u^{-2k}\ , \quad \eta^{inv}_k\in {\mathbb C}^*\ ,\label{etainvAp}
\eeqa
is uniquely determined by the recursion
\beqa
\sum_{\ell=0}^{p} \eta_{p-\ell}\eta_\ell^{inv} =0 \ , \qquad p\in {\mathbb N}^*\ .\label{recetainv}
\eeqa 
\end{lem}

\begin{proof}
    First, using~\eqref{eq:durho}, one finds:
\beqa
f_0=1 \ ,\qquad f_{2p+2} &=& (p+1)(\epsp^2+\epsm^2)
\ ,  \label{fpair}\\
f_{2p+1} &=&  (2p+1)\epsp\epsm \ ,\quad p\in {\mathbb N}\ .\label{fimpair}
\eeqa
To determine the scalars $\{\eta_k\}_{k\in{\mathbb N}}$, we proceed as follows. According to~\eqref{etaetaAp} and $f_0=1$, we can set $\eta_0=\pm 1$.  Then, for $k\geq 1$ observe that (\ref{etaetaAp}) is equivalent to:
\beqa
\eta(q^{2}u^{-2})f(u)= \eta(q^{-2}u^{-2})f(uq)\ , 
\eeqa
which is  equivalent to the recurrence relations~\eqref{eq:eta-recur}.
Those relations determine uniquely $\{\eta_k\}_{k\in {\mathbb N}^*}$ in terms of $\{f_{k}\}_{k\in {\mathbb N}^*}$ with (\ref{fpair}), (\ref{fimpair}), and with an overall factor $\eta_0=\pm1$.

Finally, using $\eta(u^{-2})\eta(u^{-2})^{-1}=1$, we get $\eta^{inv}_0=\eta_0$ and the recurrence relations~\eqref{recetainv}.
These relations determine uniquely  $\{\eta^{inv}_k\}_{k\in {\mathbb N}^*}$  in terms of  $\{\eta_k\}_{k\in {\mathbb N}^*}$.
\end{proof}

\begin{example}\label{exC1} Let us fix $\eta_0=\eta^{inv}_0=1$, the choice we use elsewhere. We then have:
	\beqa
	 \eta_1 = -\frac{\epsp\epsm}{q^2+1}\ ,\quad \eta_2 = \frac{q^{-2}}{q^4-q^{-4}}\left(  \left(\frac{q^3-q^{-3}}{q-q^{-1}}\right)\epsp^2\epsm^2 -(q^2-q^{-2})(\epsp^2+\epsm^2)  \right)\ .
	\eeqa

	\beqa
	 \eta^{inv}_1 = \frac{\epsp\epsm}{q^2+1}\ ,\quad \eta^{inv}_2 =
	 \frac{q^{-2}}{(q+q^{-1})^2(q^2+q^{-2})}\left(  (q+q^{-1})^2(\epsp^2+\epsm^2)  - \epsp^2\epsm^2 \right)\ .\label{coeffeta0inv}
	\eeqa
\end{example}

\medskip

Finally, the generating function  ${\bar B}(z)$ associated with the non-standard case is expressed in terms of the generating function associated with the standard case, previously denoted ${\bar B}^{(0)}(z)$.
Starting from~\eqref{BinGGred} with~\eqref{etaetaAp}, we obtain the non-standard analog of (\ref{BbarBeqst}) given by:
\beqa
B(u^{-2})\eta(q^2u^{-2}) {\bar B}(q^2u^{-2})= B(q^{-2}u^{-2})\eta(q^{-2}u^{-2}){\bar B}(q^{-2}u^{-2})\ .\label{BbarBeqncAp}
\eeqa
\begin{prop}\label{prop:Bbar}
The generating function $\bar B(z)$ that satisfies~\eqref{BinGGred} is given by
\beqa
{\bar B}(z) =\eta(z)^{-1} {\bar B}^{(0)}(z)\ .\label{eq:Bbargen}
\eeqa
where $\eta(z)^{-1}$ and ${\bar B}^{(0)}(z)$ are the power series~\eqref{etainvAp} and~\eqref{Bbarst} with coefficients that solve~\eqref{recetainv} and ~\eqref{eq:recBbarBdel}, respectively.
\end{prop}
\begin{proof} By direct comparison of~\eqref{BbarBeqncAp}
with~\eqref{BbarBeqst}, and because the latter equation has a unique solution.
\end{proof}

\section{Examples of refined $q$-characters for $j=3/2$ and $j=2$}\label{app:refchar}
The following expressions for the refined $q$-characters are computed using~\eqref{eq:chichirecOq},~\eqref{eq:chirefhalf} and~\eqref{eq:acutchi-spin1}
Below, we identify $\alpha=\frac{\bar k_+}{k_+}$, $\beta=\frac{\bar k_-}{k_-}$ and $a=q^{-2}v^{-2}$.

{\small
	\begin{align}
		\acute{\chi}_q\big(V_v^{(3/2)}\big)
		&= \kappa(aq^{-2})\,\kappa(a)\,\kappa(aq^{2})\Big(
		\alpha^{3}\, {\acute Y}_{aq^{-2}}\,{\acute Y}_{a}\,{\acute Y}_{aq^{2}}
		+ \alpha^{2}\beta\, {\acute Y}_{aq^{-2}}\,{\acute Y}_{a}\,{\acute Y}_{aq^{4}}^{-1}\nonumber \\
		&\qquad\qquad\qquad\qquad\qquad
		+ \alpha^{2}\beta\, {\acute Y}_{aq^{2}}
		+ \alpha\beta^{2}\, {\acute Y}_{aq^{4}}^{-1}
		+ \alpha^{2}\beta\, {\acute Y}_{aq^{-2}} \nonumber\\
		&\qquad \qquad\qquad\qquad\qquad
		+ \alpha\beta^{2}\, {\acute Y}_{aq^{-2}}\,{\acute Y}_{aq^{2}}^{-1}\,{\acute Y}_{aq^{4}}^{-1}
		+ \alpha\beta^{2}\, {\acute Y}_{a}^{-1}
		+ \beta^{3}\, {\acute Y}_{a}^{-1}\,{\acute Y}_{aq^{2}}^{-1}\,{\acute Y}_{aq^{4}}^{-1}
		\Big) \label{eq:chi32} \\
		&\quad \qquad\qquad\qquad\qquad\quad
		- \kappa(aq^{2})\big(\alpha\, {\acute Y}_{aq^{2}} + \beta\, {\acute Y}_{aq^{4}}^{-1}\big)
		- \kappa(aq^{-2})\big(\alpha\, {\acute Y}_{aq^{-2}} + \beta\, {\acute Y}_{a}^{-1}\big)\Bigg|_{z\rightarrow u^{-2}}\ ,\nonumber
	\end{align}
}
	
	{\small 
	\begin{align}
	\acute{\chi}_q\big(V_v^{(2)}\big)
		&= \kappa(aq^{-3})\,\kappa(aq^{-1})\,\kappa(aq)\,\kappa(aq^{3})\Big(
		\alpha^{4}\, {\acute Y}_{aq^{-3}}\,{\acute Y}_{aq^{-1}}\,{\acute Y}_{aq}\,{\acute Y}_{aq^{3}} \nonumber\\
		& \qquad \qquad\qquad\qquad\qquad\qquad\qquad
		+ \alpha^{3}\beta\, {\acute Y}_{aq^{-3}}\,{\acute Y}_{aq^{-1}}\,{\acute Y}_{aq}\,{\acute Y}_{aq^{5}}^{-1}
		+ \alpha^{2}\beta^{2}\, {\acute Y}_{aq^{-3}}\,{\acute Y}_{aq^{-1}}\,{\acute Y}_{aq^{3}}^{-1}\,{\acute Y}_{aq^{5}}^{-1} \nonumber\\
		&\qquad \qquad\qquad\qquad\qquad\qquad\qquad
		+ \alpha\beta^{3}\, {\acute Y}_{aq^{-3}}\,{\acute Y}_{aq}^{-1}\,{\acute Y}_{aq^{3}}^{-1}\,{\acute Y}_{aq^{5}}^{-1}
		+ \beta^{4}\, {\acute Y}_{aq^{-1}}^{-1}\,{\acute Y}_{aq}^{-1}\,{\acute Y}_{aq^{3}}^{-1}\,{\acute Y}_{aq^{5}}^{-1} \nonumber\\
		&\qquad \qquad\qquad\qquad\qquad\qquad\qquad
		+ \alpha^{3}\beta\, {\acute Y}_{aq}\,{\acute Y}_{aq^{3}}
		+ \alpha^{3}\beta\, {\acute Y}_{aq^{-3}}\,{\acute Y}_{aq^{3}}
		+ \alpha^{2}\beta^{2}\, {\acute Y}_{aq^{-1}}^{-1}\,{\acute Y}_{aq^{3}}\nonumber \\
		&\qquad \qquad\qquad\qquad\qquad\qquad\qquad
		+ \alpha^{3}\beta\, {\acute Y}_{aq^{-3}}\,{\acute Y}_{aq^{-1}}
		+ \alpha^{2}\beta^{2}\, {\acute Y}_{aq^{-3}}\,{\acute Y}_{aq}^{-1}
		+ \alpha\beta^{3}\, {\acute Y}_{aq^{-1}}^{-1}\,{\acute Y}_{aq}^{-1} \nonumber\\
		&\qquad \qquad\qquad\qquad\qquad\qquad\qquad
		+ \alpha^{2}\beta^{2}\, {\acute Y}_{aq}\,{\acute Y}_{aq^{5}}^{-1}
		+ \alpha^{2}\beta^{2}\, {\acute Y}_{aq^{-3}}\,{\acute Y}_{aq^{5}}^{-1}
		+ \alpha\beta^{3}\, {\acute Y}_{aq^{-1}}^{-1}\,{\acute Y}_{aq^{5}}^{-1} \nonumber\\
		&\qquad \qquad\qquad\qquad\qquad\qquad\qquad
		+ \alpha\beta^{3}\, {\acute Y}_{aq^{3}}^{-1}\,{\acute Y}_{aq^{5}}^{-1}
		+ \alpha^{2}\beta^{2}
		\Big) \nonumber\\
		&\qquad \qquad\quad
		- \kappa(aq)\,\kappa(aq^{3})\Big(
		\alpha^{2}\, {\acute Y}_{aq}\,{\acute Y}_{aq^{3}}
		+ \alpha\beta\, {\acute Y}_{aq}\,{\acute Y}_{aq^{5}}^{-1}
		+ \beta^{2}\, {\acute Y}_{aq^{3}}^{-1}\,{\acute Y}_{aq^{5}}^{-1}
		+ \alpha\beta
		\Big) \label{eq:chi2}\\
		&\qquad \qquad\quad
		- \kappa(aq^{-3})\,\kappa(aq^{3})\Big(
		\alpha^{2}\, {\acute Y}_{aq^{-3}}\,{\acute Y}_{aq^{3}}
		+ \alpha\beta\, {\acute Y}_{aq^{-1}}^{-1}\,{\acute Y}_{aq^{3}}
		+ \alpha\beta\, {\acute Y}_{aq^{-3}}\,{\acute Y}_{aq^{5}}^{-1} + \beta^{2}\, {\acute Y}_{aq^{-1}}^{-1}\,{\acute Y}_{aq^{5}}^{-1}
		\Big)\nonumber\\
		&\qquad\qquad\quad
		- \kappa(aq^{-3})\,\kappa(aq^{-1})\Big(
		\alpha^{2}\, {\acute Y}_{aq^{-3}}\,{\acute Y}_{aq^{-1}}
		+ \alpha\beta\, {\acute Y}_{aq^{-3}}\,{\acute Y}_{aq}^{-1}
		+ \beta^{2}\, {\acute Y}_{aq^{-1}}^{-1}\,{\acute Y}_{aq}^{-1}
		+ \alpha\beta
		\Big)\Bigg|_{z\rightarrow u^{-2}} \nonumber\\
		&\qquad \qquad\quad+ 1 \ .\nonumber
	\end{align}
}
\normalsize

\end{appendix}

\end{document}